\documentclass{article}
\usepackage[english]{babel}
\usepackage[utf8]{inputenc}
\usepackage[T1]{fontenc}
\usepackage{graphicx}
\usepackage{bm} 
\usepackage{mathrsfs}
\usepackage{mathtools}
\usepackage{units}
\usepackage{tensor}
\usepackage{accents}
\usepackage{xspace}
\usepackage{enumitem}
\usepackage{dsfont}
\usepackage{hyperref}
\usepackage{bookmark}
\usepackage{cleveref}
\usepackage{amssymb}
\usepackage{amsmath}
\usepackage{helvet}
\usepackage{mathrsfs}
\usepackage{pgf, tikz, pgfplots}
\usepackage[margin=1.3in]{geometry}
\pgfplotsset{compat=1.15}
\usetikzlibrary{shapes.geometric, arrows,backgrounds,fit,positioning}

\newcommand{\R}{\mathbb R}

\newcommand{\llbracket}{[\![}
\newcommand{\rrbracket}{]\!]}

\newcommand{\T}{\mathcal T}    
\newcommand{\E}{\mathcal E}    
\newcommand{\mN}{\mathcal N}    
\newcommand{\Ebw}{{\eta_{\mathrm{bw}}}}        
\newcommand{\jump}[1]{\left\llbracket#1\right\rrbracket} 

\DeclareMathOperator{\im}{Im}    
\DeclareMathOperator{\di}{div}    
\DeclareMathOperator{\Id}{Id}    

\DeclarePairedDelimiter{\norm}{\lVert}{\rVert} 
\DeclarePairedDelimiter{\abs}{\lvert}{\rvert}    
\newcommand{\dx}{\text{ d$x$}}    
\newcommand{\de}[1]{\text{ d$#1$}} 
\renewcommand{\phi}{\varphi}

\renewcommand{\epsilon}{\varepsilon}

\newcommand{\Pcal}{\mathcal P}    

\newcommand{\scal}[2]{\left(#1,#2\right)}    
\newcommand{\tr}[1]{#1^{\text{T}}}   

\newcommand{\comment}[1]{}  
\newcommand{\bw}{\mathrm{bw}}

\tikzset{block/.style={rectangle, draw, text width=4.2em, text centered, minimum height=3em, rounded corners=5pt},
title/.style={font=\bfseries},
line/.style={draw, -latex'},
int/.style={draw, minimum size=2em, font=\upshape}}

\title{Hierarchical a posteriori error estimation of Bank--Weiser type in the
FEniCS Project\thanks{R.B.\ would like to acknowledge the support of the
ASSIST research project of the University of Luxembourg. This publication
has been prepared in the framework of the DRIVEN project funded by the
European Union's Horizon 2020 Research and Innovation programme under Grant
Agreement No.\ 811099. F.C.'s work is partially supported by the I-Site BFC
project NAANoD and the EIPHI Graduate School (contract ANR-17-EURE-0002).}}

\author{Raphaël Bulle\thanks{Institute of Computational Engineering, University of
Luxembourg, 6 Avenue de la Fonte, 4362 Esch-sur-Alzette,
Luxembourg (\texttt{raphael.bulle@uni.lu}, \texttt{stephane.bordas@uni.lu}, \texttt{jack.hale@uni.lu})}
\and Jack S. Hale\footnotemark[2]
\and Alexei Lozinski\thanks{Laboratoire de Mathématiques de Besançon,
UMR CNRS 6623, Université de Bourgogne Franche-Comté, 16 route de Gray, 25030
Besançon Cedex, France (\texttt{alexei.lozinski@univ-fcomte.fr})} 
\and Stéphane P. A. Bordas\footnotemark[2]
\and Franz Chouly\thanks{Université de Bourgogne Franche-Comté, Institut de
Mathématiques de Bourgogne, 21078 Dijon, France
(\texttt{franz.chouly@u-bourgogne.fr})}
}

\begin{document}
    \maketitle
    \begin{center}
        \begin{minipage}{0.8\textwidth}
        \footnotesize
        \textbf{Abstract.} In the seminal paper of Bank and Weiser [\emph{Math.\ Comp.}, 44 (1985), pp.\
        283--301] a new a posteriori estimator was introduced.
        This estimator requires the solution of a local Neumann problem on every cell of the
        finite element mesh.
        Despite the promise of Bank--Weiser type estimators, namely locality,
        computational efficiency, and asymptotic sharpness, they have seen little use
        in practical computational problems.
        The focus of this contribution is to describe a novel implementation of
        hierarchical estimators of the Bank--Weiser type in a modern high-level finite
        element software with automatic code generation capabilities.
        We show how to use the estimator to drive (goal-oriented) adaptive mesh
        refinement and to mixed approximations of the nearly-incompressible
        elasticity problems.
        We provide comparisons with various other used estimators.
        An open source implementation based on the FEniCS Project finite element
        software is provided as supplementary material.
    \end{minipage}
    \end{center}

    \section{Introduction}
    A posteriori error estimation\ \cite{ainsworth_posteriori_2011} is the de
    facto tool for assessing the discretization error of finite element method (FEM)
    simulations, and iteratively reducing that error using adaptive mesh
    refinement strategies\ \cite{nochetto_theory_2009}.

    This paper is concerned with the description and justification of an
    implementation of an error estimator introduced in the seminal paper of
    Bank and Weiser\ \cite[Section 6]{bank_posteriori_1985}.
    In that paper an error estimate was derived involving the solution of local
    Neumann problems on a special finite element built on nested or
    hierarchical spaces. 
    Despite its excellent performance and low computational cost, this estimator
    has seen relatively sparse use in practical computational problems.
    The overarching goal of this contribution is to provide access to an
    efficient, generic and extensible implementation of Bank--Weiser type
    estimators in a modern and widely used finite element software,
    specifically, the FEniCS Project\ \cite{alnaes_fenics_2015}.

    \subsection{Background}
    The literature on a posteriori error estimation and adaptive finite
    element methods is vast, so we focus on articles on practical software
    implementations of adaptive finite element methods and comparative
    performance studies.

    The T-IFISS\ \cite{bespalov_t-ifiss_2019} software package, based on the
    existing IFISS\ \cite{elman_ifiss_2014} package, is a finite element
    software written in MATLAB/Octave with a focus on a posteriori error
    estimation and adaptive finite element methods.
    Recently\ \cite{bespalov_goal-oriented_2019}, T-IFSS has been
    extended to solve adaptive stochastic Galerkin finite element methods.
    The stated emphasis of T-IFISS\ \cite{bespalov_t-ifiss_2019} is on being a
    laboratory for experimentation and exploration, and also to enable the
    rapid prototyping and testing of new adaptive finite element methods.
    A number of estimation and marking strategies are implemented in T-IFISS,
    although not the Bank--Weiser estimator we consider in this paper.
    T-IFISS only works for two-dimensional problems and it was never intended
    to be a high-performance code suitable for large-scale computations e.g.\
    high-performance computing systems using the Message Passing Interface
    (MPI).

    The PLTMG package\ \cite{bank_pltmg_1998} is one of the oldest open finite
    element softwares for solving elliptic problems that is still under active
    maintenance, and includes many advanced features such as $hp$-adaptive
    refinement, a posteriori error estimation, domain decomposition and
    multigrid preconditioning.
    The a posteriori error estimation is based on a superconvergent patch
    recovery estimation technique introduced
    in\ \cite{bank_superconvergent_2007}.
    PLTMG only works in two dimensions and is naturally limited from a
    usability perspective due to the programming tools available at its
    inception (Fortran and ANSI C).

    In\ \cite{funken_efficient_2011} an adaptive first-order polynomial finite
    element method was implemented in a code called p1afem using MATLAB.
    The primary goal was to show how the basic finite element algorithm could
    be implemented efficiently using MATLAB's vectorization capabilities.
    A standard residual estimator\ \cite{bubuska_feedback_1984} is used to
    drive an adaptive mesh refinement algorithm.
    Again, like T-IFISS, p1afem only works in two dimensions.

    In\ \cite{rognes_automated_2013} a novel methodology for automatically
    deriving adaptive finite element methods from the high-level specification
    of the goal functional and (potentially non-linear) residual equation was
    implemented in the FEniCS Project.
    The emphasis of the paper\ \cite{rognes_automated_2013}, in contrast with
    the T-IFISS toolbox\ \cite{bespalov_t-ifiss_2019}, is on the
    \emph{automatic} construction of goal-oriented adaptive finite element
    methods, without much knowledge required on the part of the user.
    The implicit residual problems are automatically localised using bubble
    functions living on the interior and facets of the cell, and the dual
    problem~\cite{giles_adjoint_2002} is derived and solved automatically on
    the same finite element space as the primal problem, before being
    extrapolated to a higher-order finite element space using a patch-wise
    extrapolation operator.
    In practice the automatically derived estimators seem to be able to
    effectively drive adaptive mesh refinement for a range of different PDEs.

    Explicit residual estimators are also commonly employed by users of
    high-level finite element software packages as they can usually be
    expressed straightforwardly in a high-level form language, e.g.\
    \cite{alnaes_fenics_2015,prudhomme_feel_2012}.
    For example,\ \cite{hale_simple_2018} used the FEniCS Project to implement
    an explicit residual error estimator for the Reissner-Mindlin plate
    problem from\ \cite{beirao_da_veiga_priori_2008}.
    The authors of\ \cite{duprez_quantifying_2020} used the FEniCS Project to
    implement an explicit residual estimator for elasticity problems within a
    dual-weighted residual framework.
    The dual problem is solved on a higher-order finite element space in order
    to ensure that the weighting by the dual residual solution does not vanish
    \cite{rognes_automated_2013}.
    In\ \cite{houston_automatic_2018} the authors use an explicit
    dual-weighted residual strategy for adaptive mesh refinement of
    discontinuous Galerkin finite element methods.
    In addition, as the name suggests, they can be explicitly computed as they
    involve only functions of the known finite element solution and the
    problem data.

    In the present work, aside of the Bank--Weiser estimator we will consider an
    explicit residual estimator\ \cite{babuska_posteriori_1978} named
    \textit{residual estimator} in the following, a flux reconstruction based on
    averaging technique estimator~\cite{zienkiewicz_simple_1987}, referred to as
    \textit{Zienkiewicz--Zhu estimator}, and a variant of the Bank--Weiser
    estimator introduced in\ \cite{verfurth_posteriori_1994} and referred to as
    the \textit{bubble Bank--Weiser estimator}.
    The residual estimator was proved to be both reliable and (locally)
    efficient in\ \cite{verfurth_posteriori_1994} for any finite element order and in
    any dimension.
    The proof of reliability and (local) efficiency of Zienkiewicz--Zhu
    estimator has been derived in\ \cite{rodriguez_remarks_1994}, for linear finite
    elements in dimension two and generalised to any averaging technique in
    any dimension in\ \cite{carstensen_each_2002} and any finite element order
    in\ \cite{bartels_each_2002}.
    The bubble Bank--Weiser estimator was proved to be reliable and locally efficient
    in\ \cite{verfurth_posteriori_1994} for any dimension and any finite element order.

    A proof of the equivalence between the Bank--Weiser estimator and the exact
    error was derived in the original paper\ \cite{bank_posteriori_1985}.
    However, this proof requires a \textit{saturation assumption}\
    \cite{bank_posteriori_1985, dorfler_small_2002, nochetto_removing_1993}
    asking for the best approximation with higher order finite elements to be
    strictly smaller than that of lower order elements and which is known to be
    tricky to assert in practice.
    Some progress has been made in\ \cite{nochetto_removing_1993} removing the saturation
    assumption from the analysis.
    However, this progress was made at the price of restricting the framework to
    linear polynomial finite elements and dimension two only.
    The equivalence proof between Bank--Weiser and residual estimators have been
    improved by the authors in\ \cite{bulle_removing_2020} where it was extended to dimension three.

    \subsection{Contribution}
    We show how robust and cheap hierarchical error estimation strategies can be
    implemented in a high-level finite element framework, e.g.\ the FEniCS Project\
    \cite{alnaes_fenics_2015}, Firedrake\
    \cite{gibson_slate_2019,rathgeber_firedrake_2017}, freefem\texttt{++}\
    \cite{hecht_freefem_2012}, Feel\texttt{++}\ \cite{prudhomme_feel_2012},
    GetFEM\ \cite{RenardGetFEM} or
    Concha\ \cite{concha2010}.
    Specifically, the contribution of our paper to the existing literature is:
    \begin{itemize}
        \item A generic and efficient implementation of the Bank--Weiser
            estimator in the open source FEniCS Project finite element software
            that works for Lagrange finite elements of arbitrary polynomial
            order and in two and three spatial dimensions. We provide
            implementations for the popular but legacy DOLFIN finite
            element solver~\cite{alnaes_fenics_2015}, and the new DOLFINx
            solver\ \cite{Habera2020}. The two versions are functionally identical, although in terms of overall speed and parallel scaling the DOLFINx version is superior due to underlying architectural improvements. Hence we only show parallel scaling results with this new version.
        The code is released under an open source (LGPLv3) license\
        \cite{hale_implementation_2020}.
        Because the code utilises the existing automatic code generation
        capabilities of FEniCS along with a custom finite element
        assembly routine, the packages are very compact (a few hundred lines of
        code, plus documentation and demos). Additionally, the estimators are
        implemented in near mathematical notation using the Unified Form
        Language, see the Appendices for code snippets.
        \item A numerical comparison of the Bank--Weiser estimator with various
        estimators mentioned earlier.
        We examine the relative efficiency, and their performance within an
        adaptive mesh refinement loop on various test problems.
        Unlike\ \cite{carstensen_estimator_2010}, we do not aim at running a
        competition of error estimators but at stressing the potential of
        the Bank--Weiser estimator since, as the authors of\
        \cite{carstensen_estimator_2010} point out, a single error
        estimation strategy is not sufficient to cover the particulars of all
        possible problems.
        \item Relying on results in\ \cite{becker_weighted_2011}, we show a
        goal-oriented adaptive mesh refinement algorithm can be driven by
        weighted sum of estimators, computed separately on primal and dual
        problems discretized on the same finite element space.
        This avoids the extrapolation operation of\ \cite{rognes_automated_2013} 
        or the need to compute the dual solution
        in a higher-order finite element space\ \cite{becker_optimal_2001}.
        \item Using the same basic methodology as for the Poisson problem, we
        extend our approach to estimating errors in mixed approximation of
        nearly incompressible elasticity problems.
        This idea was originally introduced in\
        \cite{ainsworth_posteriori_2011} and is still an active research
        topic, see e.g.\ \cite{khan_robust_2019} for a parameter-robust implicit
        residual estimator for nearly-incompressible elasticity.
    \end{itemize}

    \subsection{Outline}
    An outline of this paper is as follows;
    in\ \cref{subsec:notation} we outline
    the main notation and definitions used in this paper.
    In\ \cref{sec:primal,sec:bw} we show the derivation of the primal problem
    and the Bank--Weiser error estimator.
    In\ \cref{sec:algimp} we derive a new method for computing the Bank--Weiser
    estimator and discuss its implementation in FEniCS.
    In\ \cref{sec:apps} we discuss the use of the approach for various
    applications such as goal-oriented adaptive mesh refinement and for mixed
    approximations of PDEs.
    Then, in\ \cref{sec:results} we show some results on two and three
    dimensional Poisson test problems as well as on linear elasticity
    problems, before concluding in\ \cref{sec:conclusions}.

    \subsection{Notation}\label{subsec:notation}
    In this section we outline the main notations used in the rest of
    the paper.
    Let $\Omega$ be an bounded open domain of $\R^d$ ($d=1,2$ or $3$), with
    polygonal/polyhedral boundary denoted by $\Gamma := \partial \Omega$.
    We consider $\Gamma = \Gamma_D \cup \Gamma_N$ a partition of the boundary.
    We assume $\Gamma_D$ is of positive measure.
    We denote by $n: \Gamma \rightarrow \R^d$ the outward unit normal vector
    along $\Gamma$.
    Let $\omega$ be a subset of $\overline{\Omega}$.
    For $l \in \R$ we denote by $H^l(\omega)$ the Sobolev space of order $l$.
    The space $H^0(\omega) = L^2(\omega)$ is the Lebesgue space of square
    integrable functions over $\omega$.
    The space $H^l(\omega)$ is endowed with the usual inner product
    $\scal{\cdot}{\cdot}_{l,\omega}$ and norm $\norm{\cdot}_{l, \omega}$.
    We omit the subscript $l$ when $l=0$ and subscript $\omega$ when $\omega=
    \Omega$.
    We denote $H^1_D(\Omega)$ the subspace of $H^1(\Omega)$ of functions with
    zero trace on $\Gamma_D$.
    We make use of the notation $\partial_n v := \nabla v \cdot n$ for
    the normal derivative of a smooth enough function $v$.
    For $l \in \R$ and for a $d$-dimensional subset $\omega$ of $\Omega$, we
    also define the following vector fields spaces
    $\bm L^2(\omega) := \left(L^2(\omega)\right)^d$ and $\bm H^l(\omega) :=
    \left( H^l(\omega) \right)^d,$ with respective inner products defined as
    their scalar counterparts, replacing the scalar product by the Euclidean
    inner product or the Frobenius double dot product.
    The space $\bm H^1_D(\Omega)$ is the subspace of $\bm H^1(\Omega)$ of
    functions with zero trace on $\Gamma_D$.
    From now on, the bold font notation will be reserved to vector fields.
    With these notations at hand we can proceed with the rest of the paper.

    \section{Primal problem statement and finite element
    discretization}\label{sec:primal}
    We consider the Poisson problem with mixed Dirichlet and Neumann boundary
    conditions.
    Let $\Gamma = \Gamma_D \cup \Gamma_N$ be a partition of the boundary.
    We apply a Dirichlet boundary condition on $\Gamma_D$ and a Neumann
    boundary condition on $\Gamma_N$.
    Let $f \in L^2(\Omega)$, $u_D \in H^{1/2}(\Gamma_D)$ and $g \in L^2
    (\Gamma_N)$ be known data.
    We seek a function $u$:
    \begin{equation}
        -\Delta u = f \text{ in } \Omega,\qquad u = u_D \text{ on } \Gamma_D,\qquad {\partial_n u} = g \text{ on } \Gamma_N.
        \label{eq:strong_form}
    \end{equation}
    Problem\ \cref{eq:strong_form} can be written in an equivalent
    weak form: Find $u \in H^1(\Omega)$ of trace $u_D$ on $\Gamma_D$ such that
    \begin{equation}
        \scal{\nabla u}{\nabla v} = \scal{f}{v} + \scal{g}{v}_{\Gamma_N},\quad \forall v \in H^1_D(\Omega).
        \label{eq:weak_form}
    \end{equation}
    The weak problem \cref{eq:weak_form} can be discretized
    using the Lagrange finite element method.
    We take a mesh $\T$ of the domain $\Omega$, consisting of cells $\T =
    \left\lbrace T \right\rbrace$, facets $\E = \left\lbrace E \right\rbrace$
    (we call \textit{facets} the edges in dimension two and the faces in
    dimension three), and vertices $\mN = \left\lbrace\chi \right\rbrace$.
    The mesh $\T$ is supposed to be regular in Ciarlet's sense: $h_T/\rho_T
    \leqslant \gamma,\ \forall T \in \T$, where $h_T$ is the diameter of a cell
    $T$, $\rho_T$ the diameter of its inscribed ball, and $\gamma$ is a positive
    constant fixed once and for all.
    The subset of facets in the interior of the mesh (i.e.\ those that are not
    coincident with the boundary $\Gamma$) is denoted $\E_I$.
    The subset of facets lying on $\Gamma_D$ is denoted $\E_D$.
    The subset of facets lying on $\Gamma_N$ is denoted $\E_N$.
    The subset of facets lying on the boundary of the domain $\Gamma$ is
    denoted $\E_B = \E_D \cup \E_N$.
    Until the end of this work we assume that the mesh resolves the boundary
    conditions, in other words for any edge $E \in \Gamma$ then $E \in \Gamma_D$
    or $E \in \Gamma_N$.
    Let $n^{+} \in \R^d$ and $n^{-} \in \R^d$ be the outward unit normals to a
    given edge as seen by two cells $T^{+}$ and $T^{-}$ incident to a common
    edge $E$.
    If we denote $\Pcal_k(T)$ the space of polynomials of order
    $k$ on a cell $T$, the continuous Lagrange finite
    element space of order $k$ on the mesh $\T$ is defined by
    \begin{equation}
        V^k := \left\lbrace v_k \in H^1(\Omega), v_{k|T} \in \Pcal_k(T) \;
        \forall T \in \T \right\rbrace.
    \end{equation}
    We denote $V_D^k$ the finite element space composed of functions of $V^k$
    vanishing on the boundary $\Gamma_D$.
    We consider the finite element problem: Find $u_k \in V^k$ such that $u_k =
    u_{D,k}$ on $\Gamma_D$ and:
    \begin{equation}
        \scal{\nabla u_k}{\nabla v_k} = \scal{f}{v_k} +
        \scal{g}{v_k}_{\Gamma_N},\quad \forall v_k \in V_D^k,
        \label{eq:discrete_weak_form}
    \end{equation}
    and where $u_{D,k}$ is a discretization of $u_D$ on $V^k$ (for example the
    Laplace interpolation or a $L^2$ orthogonal projection).

    \section{The Bank--Weiser estimator}\label{sec:bw}
    In this section we derive the general definition of the Bank--Weiser
    estimator from the equation of the error as it was given in the original
    paper\ \cite{bank_posteriori_1985}.
    We also give a concrete example of the Bank--Weiser estimator for linear finite
    elements.

    \subsection{The global error equation}
    We are interested in estimating the error we commit by approximating
    the solution $u$ by $u_k \in V^k_D$.
    We define this error by the function $e := u - u_k$ and we want to estimate
    its norm $\norm{e}_1$.
    The first step towards this will be to derive a new variational problem for
    which the exact error $e$ is the solution.
    For a cell $T$ of the mesh, we introduce the interior residual as
    \begin{equation}
        r_T := (f + \Delta u_k)_{|T},
        \label{eq:interior_residual}
    \end{equation}
    and for an edge $E$, the edge residual
    \begin{equation}
        J_E = \left \{ \begin{array}{l l}
                0 & \text{if } E \in \E_D,\\
                \jump{\partial_n u_k}_E & \text{if } E \in \E_I,\\
                (g - \partial_n u_k)_{|E} & \text{if } E \in \E_N.
            \end{array}\right.
    \label{eq:poisson_edge_residual}
    \end{equation}
    where the notation $\jump{v}_E := v^+ - v^-$ denotes the jump in the
    value of the function across an interior facet $E \in \E_I$.
    Here, $v^+$ and $v^-$ denote the values of $v$ on the facet $E$ as
    seen by the two incident cells $T^+$ and $T^-$, respectively.
    The error function $e$ satisfies what we call the global error equation
    \begin{equation}
        \scal{\nabla e}{\nabla v} = \sum_{T \in \T} \scal{r_T}{v}_T + \sum_{E
        \in \E_I} \scal{J_E}{v}_{E} + \sum_{E \in \E_N} \scal{J_E}{v}_E, \quad
        \forall v \in H^1_D(\Omega),
        \label{eq:error_equation}
    \end{equation}
    and $e = u_D - u_k$ on the Dirichlet boundary $\Gamma_D$.

    \subsection{The local Bank--Weiser space and the Bank--Weiser
    estimator}\label{subsec:bank-weiser_estimator}
    We introduce now local finite element spaces in order to derive the
    finite element approximation of the error.
    For a cell $T$ of the mesh we define
    \begin{equation}
        V^k_{T,D} := \left \lbrace v_{k,T} \in \mathcal P_k(T),\ v_{k,T} = 0\ \text{in}\ (\Omega \setminus
        \overline T)\cup(\overline T \cap \Gamma_D)\right\rbrace,
    \end{equation}
    as well as
    \begin{equation}
        V^k_T := \left \lbrace v_{k,T} \in \mathcal P_k(T)\right\rbrace.
    \end{equation}
    A key idea in the Bank--Weiser estimator derivation is to introduce an
    appropriate finite element space for the discretization of error.
    This non-standard space has two roles.
    Firstly, for the local problems involving the cells with facets only in
    the interior of the domain or on the Neumann boundary, it should remove
    the constant functions, giving a unique solution.
    Secondly, and as we will notice in \cref{sec:results}, solving the local error equation on
    the finite element space $V^k_{T,D}/\R$ does not necessary lead to an
    accurate estimation of the error.
    However, in some cases, the estimation of the error can be surprisingly
    accurate when the space is judiciously chosen.
    We refer the reader to\ \cite{ainsworth_influence_1996} for a full
    discussion.

    Before introducing this non-standard space, we need some more
    notations.
    Let $k_+$ and $k_-$ be two non-negative integers such that $k_+ > k_-
    \geqslant 0$.
    Let $\widetilde{T}$ be the reference cell fixed once for all (independent
    from the mesh $\T$).
    We denote
    \begin{equation}
        \mathcal L_{\widetilde{T}}:\ V_{\widetilde{T}}^{k_+} \longrightarrow
        V_{\widetilde{T}}^{k_+},\quad \im(\mathcal L_{\widetilde{T}}) =
        V_{\widetilde T}^{k_-},
    \end{equation}
    the Lagrange interpolation operator between the local spaces
    $V_{\widetilde{T}}^{k_+}$ and $V_{\widetilde{T}}^{k_-} \subset
    V_{\widetilde{T}}^{k_+}$.
    Moreover, for any cell $T$ of the mesh, there exists an affine bijection
    \begin{equation}
        \begin{array}{r c c l}
            S : & \widetilde{T} & \longrightarrow & T\\
            & \widetilde{x} & \longmapsto & S(\widetilde{x}) =: x
        \end{array}
    \end{equation}
    mapping $\widetilde{T}$ onto $T$.
    From the mapping $S$ we deduce another mapping given by
    \begin{equation}
        \begin{array}{r c c l}
            \mathcal S: & V_T^{k_+} & \longrightarrow & V_{\widetilde T}^{k_+}\\
            & v(x) & \longmapsto & \mathcal S(v)(\widetilde x) := v(S
            (\widetilde x)).
        \end{array}
    \end{equation}
    If we denote $d_+$ the dimension of $V_{\widetilde T}^{k_+}$ and $d_-$
    the dimension of $V_{\widetilde T}^{k_-}$, given $\mathcal B_{\widetilde
    T}^+ := \{\widetilde \phi_1,\cdots, \widetilde \phi_{d_+}\}$ the
    basis of shape functions of $V_{\widetilde{T}}^{k_+}$ and $\mathcal
    B_T^+ := \{\phi_1,\cdots, \phi_{d_+}\}$ the basis of $V_T^{k_+}$, we
    can always find a mapping $S$ (and a mapping $\mathcal S$) such that
    \begin{equation}
        \mathcal S(\phi_{T,i}) = \widetilde \phi_{T,i},\quad \forall i
        \in \{1,\cdots,d_+\},
        \label{eq:basis_order}
    \end{equation}
    We choose $S$ and $\mathcal S$ so.
    For a given cell $T$ of the mesh, we define the Lagrange interpolation
    operator on $T$ as follows
    \begin{equation}
        \mathcal L_T := \mathcal S^{-1} \circ \mathcal L_{\widetilde T} \circ
        \mathcal S.
        \label{eq:link_operator_cells}
    \end{equation}
    Note, due to \eqref{eq:basis_order}, the matrix of $\mathcal S$ in
    the couple of basis $(\mathcal B_T^+, \mathcal B_{\widetilde T}^+)$ is
    the identity matrix of size $d_+ \times d_+$.
    Consequently, if we denote $G$ the matrix of $\mathcal L_T$ in the basis
    $\mathcal B_T^+$ and $\widetilde G$ the matrix of $\mathcal L_{\widetilde
    T}$ in the basis $\mathcal B_{\widetilde T}^-$, we have
    \begin{equation}
        G = \Id^{-1} \widetilde G \Id = \widetilde G.
        \label{eq:interpolation_independent_cell}
    \end{equation}
    For a cell $T$ of the mesh, the \emph{local Bank--Weiser space} $V_T^{\bw}$
    is defined as the null space of $\mathcal L_T$, in other words
    \begin{equation}
        V_T^{\bw} := \ker(\mathcal L_T) = \left\lbrace v_T^{\bw} \in V_T^{k_+},\ 
        \mathcal L_T v_T^{\bw} = 0 \right\rbrace.
    \end{equation}
    Similarly, we define
    \begin{equation}
        V_{T,D}^{\bw} := \left \{ v_T^{\bw} \in V^{\bw}_T,\ v_T^{\bw}=0\ \text{on
        } \overline{T}\cap \Gamma_D\right\}.
    \end{equation}
    With these new spaces in hands, we can derive a local discrete counterpart of
    equation \cref{eq:error_equation} on any cell $T$:
    Find $e_T \in V_T^{\bw}$ such that:
    \begin{equation}
        \scal{\nabla e_T^{\bw}}{\nabla v_T^{\bw}} = \scal{r_T}{v_T^{\bw}} +
        \frac{1}{2} \sum_{E \in \partial T} \scal{J_E}{v_T^{\bw}}_E, \quad \forall v_T^{\bw}
        \in V_{T,D}^{\bw},
        \label{eq:bw_local_error_equation}
    \end{equation}
    and $e_T = \pi_T^{\bw} \left(u_D - u_k\right)$ on $\Gamma_D$, where
    $\pi_T^{\bw}: L^2(T) \rightarrow V^{\bw}_T$ is a proper projection operator
    (the way this projection is implemented is detailed in \cref{subsec:outline}).

    Note, the definition of the edge residual $J$ takes into account the error
    on the Neumann boundary data approximation.
    The Dirichlet boundary data approximation has to be incorporated to the
    linear system during the solve of \cref{eq:bw_local_error_equation}, as well
    will see later.
    For a detailed discussion on a priori and a posteriori error
    estimation with inhomogeneous Dirichlet boundary conditions see
    \cite{Aurada2013,bartels_inhomogeneous_2004}.

    Finally, on the cell $T$ the local Bank--Weiser estimator $\eta_{\bw,T}$ is
    defined by
    \begin{equation}
        \eta_{\mathrm{bw},T} := \norm{\nabla e_T^{\bw}}_T,
        \label{eq:local_estimate}
    \end{equation}
    where $e_T$ is defined in \cref{eq:bw_local_error_equation} and the global
    Bank--Weiser estimator by the sum of local estimates
    \begin{equation}
        \eta^2_{\mathrm{bw}} := \sum_{T \in \T} \eta_{\mathrm{bw},T}^2.
        \label{eq:global_estimate}
    \end{equation}

    Note, althought it is not shown in this study, it is straightforward to
    generalize the Bank--Weiser estimator for other kind of elliptic operators
    by changing the energy norm in \cref{eq:local_estimate} accordingly.
    \subsection{A particular example}
    If we assume $k=1$ (i.e.\ we solve
    \cref{eq:discrete_weak_form} using linear finite elements) one can define
    the space $V_T^{\bw}$ from the choice of $k_+ = 2$, $k_- = k = 1$.
    This example was the case considered in the numerical tests of the original
    paper \cite{bank_posteriori_1985}.
    The space $V_T^{\bw}$ consists of quadratic polynomial
    functions (in $V_T^2$) vanishing at the degrees of freedom of the standard
    linear finite element functions (in $V_T^1$) i.e.\ the degrees of freedom
    associated with the vertices of $T$.

    \section{Algorithms and implementation details}\label{sec:algimp}
    The linear system corresponding to \cref{eq:bw_local_error_equation} is not
    accessible in FEniCS. This prevent us from directly solving the Bank--Weiser
    equation.
    We propose to bypass the problem by constructing the linear system
    corresponding to \cref{eq:bw_local_error_equation} from another linear
    system derived from finite element spaces that are accessible directly in
    FEniCS.
    \subsection{Method outline}\label{subsec:outline}
    \begin{enumerate}
        \item We consider the following singular value decomposition (SVD) of $G$
            \begin{equation}
                G = U \Sigma \tr{V},
            \end{equation}
            where $\Sigma$ is a diagonal matrix composed of the singular values
            of $G$.
            The columns of the matrix $V$ are singular vectors of $G$,
            associated with singular values.
            The columns associated with singular values zero span the null
            space of $G$.
            We take the submatrix $N$ made of the columns of $V$ spanning the
            null space of $G$.
            Note that, since $G$ does not depend on any cell $T$, the same
            property holds for $N$.
        \item We build the matrix $A_T^+$ and vector $b_T^+$ of the local linear
            system corresponding to the following variational formulation in the
            space $V_T^{k_+}$, available in FEniCS:
        \begin{equation}
            \scal{\nabla e_T^+}{\nabla v_T^+} =
            \scal{r_T}{v_T^+} + \frac{1}{2}\sum_{E \in \partial T}
            \scal{J_E}{v_T^+}_E,\quad \forall v_T^+ \in
            V_T^{k_+}.
            \label{eq:ho_equation}
        \end{equation}
            We integrate the Dirichlet boundary condition directly into $A_T^+$
            and $b_T^+$, by considering the vector associated to $\pi^+_T(u_D - u_k)$, where
            $\pi^+_T$ is the $L^2$ projection onto
            $V_T^{k_+}$.
            More precisely, the rows and columns of $A_T^+$ corresponding to degrees of
            freedom on the Dirichlet boundary are zeroed and the corresponding
            diagonal entries are replaced by ones.
            The entries of $b_T^+$ corresponding to these degrees of freedom
            are replaced by the corresponding entries in the vector of
            $\pi^+_T(u_D - u_k)$.
        \item We construct the matrix $A_T^{\bw}$ and vector $b_T^{\bw}$ as
            follow
        \begin{equation}
            A_T^{\bw} = \tr{N} A_T^+ N \quad \text{and} \quad
            b_T^{\bw} = \tr{N} b_T^+,
            \label{eq:link_N}
        \end{equation}
            where $A_T^{\bw}$ and $b_T^{\bw}$ are the matrix and vector
            which allow to recover the bilinear and linear forms of
        \cref{eq:bw_local_error_equation} in a basis of $V_T^{\bw}$.
        \item We solve the linear system
        \begin{equation}
            A_T^{\bw} x_T^{\bw} = b_T^{\bw},
            \label{eq:bw_linear_system_Vbw}
        \end{equation}
        \item We bring the solution back to $V_T^{k_+}$, considering
            $Nx_T^{\bw}$, in order to post-process it and compute the local
            contribution of the Bank--Weiser estimator
            \cref{eq:local_estimate}.
    \end{enumerate}

    \subsection{Computational details}
    We now give more details specific to our implementation in FEniCS of each one
    of the above steps.

    \begin{figure}
        \begin{center}
            \resizebox{.9\textwidth}{!}{
            \begin{tikzpicture}
                \node[int, text width=2.8cm, align=center] (i1) {Compute\\ $G=G_1G_2$};
                \node[int,below=3mm of i1, text width=2.8cm, align=center] (i2) {SVD of $G$\\ to obtain $V$};
                \node[int,below=3mm of i2, text width=2.8cm, align=center] (i3) {Extract $N$\\ from $V$};
                \node[title, below=3mm of i3, align=center] (title1) {Computation of N};
                \node[block, fit=(i1) (i2) (i3) (title1),draw,minimum width=2.8cm] (Input) {};

                \node[int, above right=-2.7mm and 1.7cm of i1, text width=2.8cm,align=center] (j1) {Compute $A_T^+$\\ and $b_T^+$};
                \coordinate[right=8mm of j1] (j1');
                \node[int,below=3mm of j1, text width=2.8cm, align=center] (j2) {Compute $A_T^{\bw}$\\ and $A_T^{\bw}$};
                \node[int,below=3mm of j2, text width=2.8cm, align=center]
                (j3) {Solve \cref{eq:bw_local_error_equation}};
                \node[int,below=3mm of j3, text width=2.8cm, align=center] (j4) {Compute local\\ BW estimator};
                \coordinate[right=8mm of j4] (j4');
                \node[title, below right=3mm and -3cm of j4, text width=3.5cm, align=center] (title2) {Computation of\\ the local estimators};
                \node[block, fit=(j1) (j1') (j2) (j3) (j4) (j4') (title2),draw,minimum width=4.5cm, right=1cm of Input] (Middle) {};

                \node[title, draw,rounded corners=5pt, right=1cm of Middle, text width=4cm, align=center] (Right) {Computation of\\ the global estimator};

                \path[line] (i1) -- (i2) ;
                \path[line] (i2) -- (i3) ;
                \path[line] (j1) -- (j2) ;
                \path[line] (j2) -- (j3) ;
                \path[line] (j3) -- (j4) ;
                \draw [thick, ->] (Input) -- (Middle) ;
                \draw[thick, -] (j4) -- (j4') ;
                \draw[thick, -] (j4') -- node[left]{\rotatebox{90}{Loop over cells}} (j1') ;
                \draw[thick, ->] (j1') -- (j1) ;
                \draw[thick, ->] (Middle) -- (Right) ;

            \end{tikzpicture}
            }
        \end{center}
        \caption{Overall process of the Bank--Weiser estimator algorithm.}
    \end{figure}
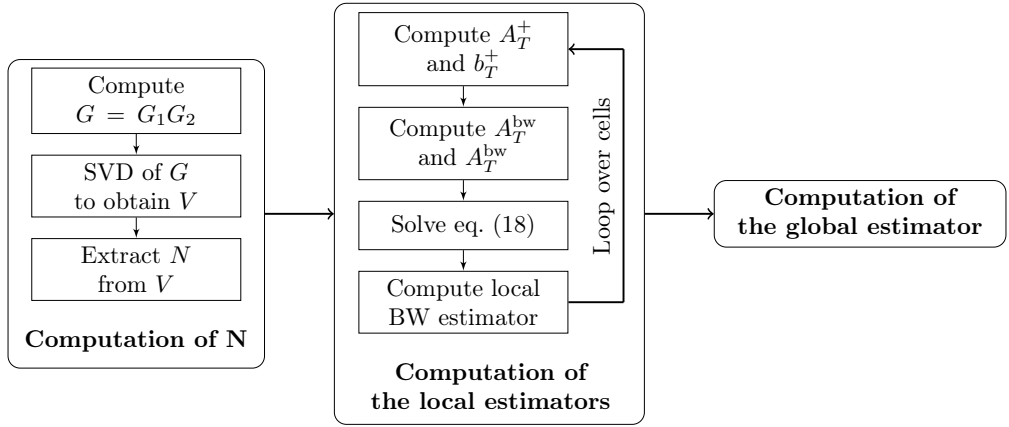

    \begin{enumerate}[leftmargin=*, label=\arabic*.]
        \item \textit{Computation of $N$.} This is the key point of our
        implementation.
        The operator $\mathcal L_T$ can be written as follows:
        \begin{equation}
            \begin{array}{r c c c c l}
                \mathcal L_T: & V_T^{k_+} & \longrightarrow & V_T^{k_-} &
                \longrightarrow & V_T^{k_+}\\
                 & v^+ & \longmapsto & \mathcal G_1(v^+) & \longmapsto &
                \mathcal G_2\big(\mathcal G_1(v^+)\big).
            \end{array}
        \end{equation}
        Then, the matrix $G$ is obtained via the following product
        \begin{equation}
            G = G_2G_1,
            \label{eq:G_decomposition}
        \end{equation}
        where $G_1$ and $G_2$ are respectively the matrix in the couple of
        basis $(\mathcal B_T^+, \mathcal B_T^-)$ of the Lagrange
        interpolation operator from $V_T^{k_+}$ to $V_T^{k_-}$, denoted
        $\mathcal G_1$ and the matrix in the same couple of basis of the
        canonical injection of $V_T^{k_-}$ into $V_T^{k_+}$, denoted $\mathcal
        G_2$.
        The matrices $G_1$ and $G_2$ can be calculated either using the
        Finite Element Automatic Tabulator (FIAT)\ \cite{kirby_algorithm_2004}
        or, as we choose to do, using the interpolator construction functions of
        the DOLFIN/x finite element library\ \cite{logg_dolfin_2010}.
        The next step consists in computing the unitary matrix $V$ of right
        singular vectors of $G$.
        This computation is done using the singular value decomposition (SVD)
        algorithm available in the SciPy library\ \cite{virtanen_scipy_2020}.
        We can write the matrix $V$ as follows,
        \begin{equation}
            V = \big(\xi_1^0\ | \cdots |\ \xi_{d_{\bw}}^0\ |\ \xi_1\ |
            \cdots | \xi_{d_-} \big),
        \end{equation}
        where $\mathcal B_T^{\bw} := \{\xi_1^0, \cdots, \xi_{d_{\bw}}^0\}$ is
        the set of singular vectors of $G$ corresponding to a zero singular
        value, spanning $V_T^{\bw}$ and $\{\xi_1,\cdots, \xi_{d_-}\}$ is
        spanning the supplementary space.
        The matrix $N$ is then chosen as the submatrix of $V$, keeping only the
        columns from $\mathcal B_T^{\bw}$:
        \begin{equation}
            N := \big(\xi_1^0\ | \cdots |\ \xi_{d_{\bw}}^0\big).
        \end{equation}
        The linear algebra operations needed to form the submatrix $N$ from $V$
        are performed using the NumPy library\ \cite{van_der_walt_numpy_2011}.
        \item \textit{Computation of $A_T^+$ and $b_T^+$.}
        The equation \cref{eq:ho_equation} is expressed directly in the Unified
        Form Language (UFL)\ \cite{alnaes_unified_2014} and efficient
        C\texttt{++} code for calculating the cell local tensors
        $A_T^+$ and $b_T^+$ for a given cell $T$ is then
        generated using the FEniCS Form Compiler (FFC)\
        \cite{kirby_compiler_2006,olgaard_automated_2009}.
            If the cell $T$ has an edge on a Dirichlet boundary $\mathcal{E}_D$, the matrix
        $A_T^+$ and vector $b_T^+$ must be modified in
        order to enforce the boundary condition.
        \item \textit{Computation of $A_T^{\bw}$ and $b_T^{\bw}$.} The matrix
        $A_T^{\bw}$ and vector $b_T^{\bw}$ are constructed using
        \cref{eq:link_N}.
    \item \textit{Solution of the linear system
        \eqref{eq:bw_linear_system_Vbw}.} The linear system
    \cref{eq:bw_linear_system_Vbw} is solved using a partial-pivot LU
    decomposition algorithm from the Eigen dense linear algebra
    library\ \cite{guennebaud_eigen_2010} in DOLFIN and xtensor-blas, which
    calls LAPACK's dgesv in DOLFINx.
        \item \textit{Computation of the Bank--Weiser estimator.} Finally, the
        solution $x_T^{\bw}$ is sent back to $V_T^{k_+}$ using $N$ and the
        norm of the corresponding function, giving the local estimator
        \cref{eq:local_estimate} is computed using standard high-level functions
        already available within FEniCS.
        The global estimator \cref{eq:global_estimate} is computed using the
        information of all the local contributions.
    \end{enumerate}
    \subsection{Additional remarks}
    \begin{itemize}
        \item The custom assembler composed of steps 2.-5.\ is performed by
            looping over every cell of the mesh and, by virtue of using the
            abstractions provided by DOLFINx, works in parallel on distributed
            memory computers using the Message Passing Interface (MPI) standard.
            For performance reasons these steps have been written in
            C\texttt{++} and wrapped in Python using the pybind11 library so
            that they are available from the Python interface to DOLFIN/x.
            In contrast, the first step must only be performed once since the
            matrix $N$ is the same for every cell of the mesh.

        \item A posteriori error estimation methods such as the one we are
            considering here assume that the linear system associted with the
            primal problem \cref{eq:weak_form} is solved exactly.
            However for performance reasons, here we use PETSc conjugate
            gradient iterative method.
            Using inexact solutions can have an influence on the total error but
            also on the a posteriori error estimator itself.
            It is a known issue \cite{Arioli2013} and several authors have
            proposed ways to estimate the algebraic error, see e.g.\
            \cite{Anciaux-Sedrakian2020,Papez2018}.
            Since algebraic error estimation is beyond the scope of this work,
            in all our numerical results we set PETSc residual tolerance small
            enough to neglect this part of the error.

        \item Because we use the automatic code generation capabilities of
            FEniCS, our approach can be readily applied to other definitions for
            the spaces $V_T^{k_+}$ and $V_T^{k_-}$, and to vectorial problems
            like linear elasticity, as we will see in the next section.

        \item For large problems the storage of the global higher order space
            $V^{k_+}$ can be an issue since it requires a lot of memory space.
            However we avoid this problem by considering the local higher order
            spaces $V_T^{k_+}$ (and local lower order spaces $V_T^{k_-}$) only.

        \item In the numerical results secton we compare several versions of
            Bank--Weier estimator and especially the one we call bubble
            Bank--Weiser estimator and denote $\eta^b_T$ which can be obtained
            with our method by taking $V_T^+$ as the space $V_T^2 +
            \mathrm{Span}\{\psi_T\}$ (the local space of quadratic functions
            enriched with the space spanned by the interior bubble function) and
            $V_T^{k_-}$ as $V_T^1$.  The resulting space $V^{\bw}_T$ is spanned
            by the interior bubble function and the edges bubbles functions of
            the cell $T$.
    \end{itemize}

    \section{Applications}\label{sec:apps}
    In this section we show a number of applications, including adaptive mesh
    refinement, goal-oriented estimation and extensions to more complex mixed
    finite element formulations for the nearly-incompressible elasticity
    problems.

    \subsection{Adaptive mesh refinement}\label{sec:adaptive}
    As well as simply providing an estimate of the global and local error, the
    estimator can be used to drive an adaptive mesh refinement strategies.
    In the following we compare different refinement strategy all based on
    the following loop:
    \begin{center}
        ... $\rightarrow$ SOLVE $\rightarrow$ ESTIMATE $\rightarrow$ MARK
        $\rightarrow$ REFINE $\rightarrow$ ...
    \end{center}
    The loop can be terminated once a given criterion e.g.\ maximum number of
    iterations, or global error less than a given tolerance, has been reached.
    A detailed discussion on adaptive refinement methods can be found in\
    \cite{nochetto_theory_2009}.
    In the following we expand on the specific algorithms used in our case.

    \subsubsection{Solve} The weak form\ \cref{eq:weak_form} is discretized
    using a standard finite element method implemented within FEniCS.
    The resulting linear systems are solved using the appropriate algorithms
    available within PETSc\ \cite{balay_petsc_2016}, e.g.\ conjugate gradient
    method preconditioned with Hypre BoomerAMG\ \cite{falgout_hypre_2002},
    or direct methods, e.g.\ MUMPS\
    \cite{amestoy_fully_2001,amestoy_hybrid_2006}.

    \subsubsection{Estimate} The Bank--Weiser estimator $\Ebw$ is formulated and
    implemented as described in \cref{sec:algimp}.
    The local contributions of the estimator provide an estimate of the local
    error for each cell in the mesh  and are subsequently used to mark the mesh.
    In addition the global estimator can be used to determine when to stop
    iterating.

    \subsubsection{Mark} We have used two distinct marking strategies throughout
    the results section: the maximum strategy on the three-dimensional test
    cases and Dörfler strategy on the two-dimensional ones.
    We follow the presentation in\ \cite{pfeiler_dorfler_2019}.
    In the maximum marking strategy\ \cite{bubuska_feedback_1984}, a cell is
    marked if its indicator is greater than a fixed fraction of the maximum
    indicator.
    More precisely, given a marking fraction $\theta \in (0, 1]$, the marked
    set $\mathcal{M} \subset \T$ is the subset such that:
    \begin{equation}
        \eta_{\bw, T} \ge \theta \max_{T \in \T} \eta_{\bw, T}, \quad \forall T \in \T.
    \end{equation}

    In the Dörfler marking strategy~\cite{dorfler_convergent_1996} (sometimes
    referred to as the equilibrated marking strategy) enough elements must
    be marked such that the sum of their estimators is larger than a fixed
    fraction of the total error.
    Given a marking fraction $\theta \in (0, 1]$, the marked set
    $\mathcal{M}$ is the subset with minimal cardinality $\#\mathcal{M}$ such
    that
    \begin{equation}
        \label{eq:dorfler}
        \sum_{T \in \mathcal{M}} \eta_{\bw, T}^2 \ge \theta \sum_{T \in \T}
        \eta_{\bw, T}^2.
    \end{equation}
    We implement an $\mathcal{O}( N \log N )$ with $N := \# \T$ complexity
    algorithm for finding the minimum cardinality set by sorting the
    indicators in decreasing order and finding the cutoff point such that
    \cref{eq:dorfler} is satisfied.
    Because of the ordering operation this set is guaranteed to have minimal
    cardinality.
    We note that recent work\ \cite{hoare_algorithm_1961,pfeiler_dorfler_2019}
    proposes a $\mathcal{O}(N)$ complexity algorithm for finding the set with
    minimum cardinality.

    \subsubsection{Refine} We use two-dimensional and three-dimensional
    variants of the algorithm proposed in~\cite{plaza_local_2000}, sometimes
    referred to as the Plaza algorithm.
    This algorithm works by subdividing the facets of each marked triangle or
    tetrahedron cell and then subdividing each triangle or tetrahedral
    cell so that it is compatible with the refinement on the facets.
    The algorithm has $\mathcal{O}(M)$ complexity in the number of added mesh
    vertices $M$.
    This algorithm already exists in DOLFIN~\cite{logg_dolfin_2010} and was
    used for the numerical results in~\cite{rognes_automated_2013}.

    \subsection{Goal-oriented adaptive mesh refinement}\label{sec:wgo}
    In many practical applications it is desirable to control the error in a
    specific quantity of interest, rather than the (global, i.e.\ across the
    entire domain $\Omega$) energy norm\ \cite{becker_optimal_2001}.
    In this section we show how the basic Bank--Weiser estimator can be used
    to control error in a goal functional, rather than in the natural norm.
    To do this, we use a weighted marking strategy proposed
    in\ \cite{becker_weighted_2011}.

    Let $\mathcal{J}: L^2(\Omega) \to \R$ be a given linear functional.
    Associated with $\mathcal J(u)$ and the primal problem
    \cref{eq:weak_form} is  the \emph{dual} or \emph{adjoint} problem:
    Find the dual solution $z \in H^1_D(\Omega)$ such that
    \begin{equation}
        \scal{\nabla v}{\nabla z} = \mathcal{J}(v),\quad \forall v \in H^1_D(\Omega).
    \end{equation}
    The dual problem, like the primal problem, can also be approximated using
    the finite element method.
    Find $z_k \in V^k$ such that
    \begin{equation}
        \scal{\nabla v_k}{\nabla z_k} = \mathcal{J}(v_k) = \scal{c}{v_k} +
        (h, v_k)_{\Gamma},\quad \forall v_k \in V^k.
    \end{equation}
    Using Galerkin orthogonality and Cauchy-Schwarz, it follows that
    \begin{align}
        \abs{\mathcal J(u) - \mathcal J(u_k)} &= \abs{\scal{\nabla(u - u_k)
        }{\nabla z}} \\
        &= \abs{\scal{\nabla(u - u_k)}{\nabla(z - z_k)}} \\
        &\le \norm{ \nabla(u - u_k)} \norm{ \nabla(z - z_k)},
        \label{eq:goal_oriented_inequality}
    \end{align}
    where the inequality holds due to Galerkin orthogonality.

    Approximating the primal and dual errors $\norm{\nabla(u - u_k)}$ and
    $\norm{\nabla(z - z_k)}$ with any estimators $\eta_u$ and $\eta_z$
    respectively, gives us an estimator for the error in the goal functional
    $|J(u) - J(u_k)|$ as the product of $\eta_u$ and $\eta_z$, thanks to
    \cref{eq:goal_oriented_inequality}:
    \begin{equation}
        \eta_w := \eta_u \eta_z
    \end{equation}
    In addition, if $\eta_u$ and $\eta_z$ are reliable estimators i.e. if there
    exist two constants $C_u$ and $C_z$ only depending on the mesh regularity
    such that
    \begin{equation}
        \norm{\nabla(u - u_k)} \leq C_u \eta_u,\quad \text{and}\quad
        \norm{\nabla(z - z_k)} \leq C_z \eta_z,
    \end{equation}
    then, $\eta_w$ is reliable as well
    \begin{equation}
        |J(u) - J(u_k)| \leqslant C_u C_z \eta_w.
    \end{equation}
    Note that because the error in the goal functional is bounded by the
    product of two estimates, the element marking strategy must incorporate
    information from local indicators for both approximations to reduce the
    error on refinement.
    There are multiple strategies for doing this in the literature, see
    e.g.\ \cite{mommer_goal-oriented_2009}.
    We have chosen to implement the weighted goal-oriented (WGO) marking
    strategy from\ \cite{becker_weighted_2011}.
    The local WGO estimator is then defined as
    \begin{equation}
        \eta_{w, T}^2 := \frac{\eta_z^2}{\eta_u^2 + \eta_z^2} \eta^2_{u, T} +
        \frac{\eta_u^2}{\eta_u^2 + \eta_z^2} \eta^2_{z, T}, \quad \forall T \in
        \T.
    \end{equation}
    The marking and refinement using $\eta_{w, T}^2$ then follows in
    exactly the same manner as in the standard adaptive refinement strategy.

    \subsection{Extension to linear elasticity problems}
    Our implementation of the Bank--Weiser estimator can be directly applied to
    mixed formulations of (nearly-incompressible) linear elasticity problems
    using the results in\ \cite{khan_robust_2019}.
    In\ \cite{ainsworth_posteriori_1997} a new a posteriori error estimator is introduced
    for mixed formulations of Stokes problems consisting in solving a local
    Poisson problem based on the local residuals on each cell.
    This estimator has been proved to be reliable and efficient in\
    \cite{ainsworth_posteriori_1997} under a saturation assumption.
    This assumption has been later removed in\ \cite{Liao_simple_2012}.
    The reliability and efficiency of the estimator for mixed formulations of
    linear elasticity is proved in\ \cite{khan_robust_2019} without the
    need of a saturation assumption.
    In addition, they show that the estimator is robust in the incompressible
    limit.
    
    \subsubsection{Nearly-incompressible elasticity}\label{sec:elasticity}
    We consider the problem of linear deformation of an isotropic elastic solid
    $\Omega$ using the Herrmann mixed formulation.
    We consider the stress tensor $\bm{\sigma} : \Omega \rightarrow \R^{d \times
    d}$, the strain tensor $\bm{\epsilon} : \Omega \rightarrow \R^{d \times d}$,
    the load $\bm{f}: \Omega \rightarrow \R^d$ which belongs to
    $\left(L^2(\Omega) \right)^d$, the Dirichlet boundary data $\bm u_D$ in
    $\left(H^{1/2}(\Gamma_D)\right)^d$, the Neumann boundary condition
    (traction) data $\bm g \in \left(L^2(\Gamma_N)\right)^d$ and displacement
    field $\bm{u}: \Omega \rightarrow \R^d$.
    The stress and strain tensors are defined by\\
    \begin{subequations}
        \begin{minipage}{.5\textwidth}
            \begin{align}
                \bm \sigma := 2 \mu \bm \epsilon(\bm u) - p \Id,
                \label{eq:stress_tensor}
            \end{align}
        \end{minipage}
        \begin{minipage}{.45\textwidth}
            \begin{align}
                \bm \epsilon(\bm{u}) := \frac{1}{2}\left( \nabla \bm u + \tr{(\nabla
                \bm u)} \right).
                \label{eq:strain_tensor}
            \end{align}
        \end{minipage}\\
    \end{subequations}
    where $\Id$ is the $d\times d$ identity matrix and $\mu$ and $\lambda$
    are the Lamé coefficients.
    The weak form of this linear elasticity problem reads: find $\bm{u}$ in
    $\bm H^1(\Omega)$ of trace $\bm u_D$ on $\Gamma_D$ and $p \in L^2(\Omega)$ such that
    \begin{subequations}
        \begin{align}
            \label{eq:elasticity:weak:equilibrium}
            2\mu \scal{\bm \epsilon(\bm u)}{\bm \epsilon(\bm v)} -
            \scal{p}{\di(\bm v)} & = \scal{\bm{f}}{\bm{v}} + \scal{\bm g}{\bm
            v}_{\Gamma_N},\quad \forall \bm v
            \in \bm H^1_D(\Omega),\\
            \label{eq:elasticity:weak:pressure}
            \scal{q}{\di(\bm u)} + \frac{1}{\lambda} \scal{p}{q} & = 0,\quad
            \forall q \in L^2(\Omega).
        \end{align}
    \end{subequations}
    The problem given by
    \cref{eq:elasticity:weak:equilibrium,eq:elasticity:weak:pressure} admits a
    unique solution (see e.g.\ \cite{khan_robust_2019}).
    We introduce the finite element spaces $X_D \subset \bm H^1_D
    (\Omega)$ and $M \subset L^2(\Omega)$ such that
    \begin{equation}
        X_D := \left(V^2_D\right)^d,
    \end{equation}
    and $M := V^1$.
    Let $\bm w$ be a discretization of $\bm u_D \in X$.
    Considering the stable Taylor--Hood method of discretization, the mixed
    finite element approximation of
    \cref{eq:elasticity:weak:equilibrium,eq:elasticity:weak:pressure} reads:
    find $\bm u_2 \in X_D$ with $\bm u_2 = \bm w$ on $\Gamma_D$ and $p_1 \in M$ such that
    \begin{subequations}
        \begin{align}
            \label{eq:elasticity:weak:discrete:equilibrium}
            2\mu \scal{\bm \epsilon(\bm u_2)}{\bm \epsilon(\bm v_2)} -
            \scal{p_1}{\di(\bm v_2)} & = \scal{\bm{f}}{\bm{v_2}} + \scal{\bm
            g}{\bm v_2},\quad \forall
            \bm v_2 \in X_D,\\
            \label{eq:elasticity:weak:discrete:pressure}
            \scal{q_1}{\di(\bm u_2)} + \frac{1}{\lambda} \scal{p_1}{q_1} & =
            0,\quad \forall q_1 \in M.
        \end{align}
    \end{subequations}
    Similarly to
    \cref{eq:elasticity:weak:equilibrium,eq:elasticity:weak:pressure}
    transposed to the discrete context,
    \cref{eq:elasticity:weak:discrete:equilibrium,eq:elasticity:weak:discrete:pressure}
    have a unique solution.
    If we denote $\bm e := \bm u - \bm u_2$ and
    $\epsilon := p - p_1$ the discretization error is measured by
    $2\mu \norm{\nabla \bm e_T} + \norm{r_T}$.

    For a cell $T$ and an edge $E$ the residuals are defined by\\
    \begin{subequations}
        \begin{minipage}{.56\textwidth}
            \begin{equation}\label{eq:elasticity_vectorial_interior_res}
                \bm R_T := \left(f +
                \di\left(2\mu \bm \epsilon(\bm u_2)
                \right) - \nabla p_1\right)_{|T},
            \end{equation}
        \end{minipage}
        \begin{minipage}{.43\textwidth}
            \begin{equation}\label{eq:elasticity_scalar_interior_res}
                r_T := \big( \di(\bm u_2) +
                \frac{1}{\lambda} p_1 \big)_{|T},
            \end{equation}
        \end{minipage}\\
        \begin{equation}\label{eq:elasticity_edge_res}
        \bm R_E = \left \{ \begin{array}{l l}
                \frac{1}{2} \jump{\left(p_1 \Id - 2 \mu \bm \epsilon(\bm
            u_2)\right) \bm n} & \text{if } E \in \E_I,\\
            0 & \text{if } E \in \E_D,\\
        \bm g - (p_1 \Id - 2\mu \epsilon(\bm u_2)) \bm n & \text{if } E \in \E_N,
                \end{array} \right.
        \end{equation}
    \end{subequations}
    Here, once again we derive the a posteriori error estimator from these
    residuals and a local Poisson problem, following\ \cite{khan_robust_2019}.
    Let $T$ be a cell of the mesh, the local Poisson problem read: find $\bm
    e_T \in \bm V_T^{\bw}$ such that
    \begin{equation}
            \label{eq:elasticity:error:poisson_equilibrium}
            2 \mu \scal{\nabla \bm e_T}{\nabla \bm v_T}_T =
            \scal{\bm R_T}{\bm v_T}_T - \sum_{E \in \partial T}
            \scal{\bm R_E}{\bm v_T}_E,\quad \forall \bm v_T \in \bm V_T^{\bw}.
    \end{equation}
    The Poisson estimator is then defined by\\
    \begin{subequations}
        \begin{minipage}{.4\textwidth}
            \begin{equation}
                \eta_{\mathrm p}^2 := \sum_{T \in \T} \eta_{\mathrm p, T}^2,
            \end{equation}
        \end{minipage}
        \begin{minipage}{.5\textwidth}
            \begin{equation}
                \eta_{\mathrm p, T}^2 := 2\mu \norm{\nabla \bm e_T}_T^2 + \norm{r_T}_T^2.
            \end{equation}
        \end{minipage}
    \end{subequations}

    This estimator has been proved to be reliable and locally efficient in\
    \cite{khan_robust_2019} as well as robust in the incompressible limit.

    \section{Results}\label{sec:results}
    We illustrate our implementation first on several two dimensional problems as
    Poisson problems with solutions of different regularities and with different
    boundary conditions.
    Then, we also look at examples of linear elasticity, and goal-oriented
    problems.
    We now treat a three dimensional example: a linear elasticity problem
    on a mesh inspired by a human femur bone.
    One can find another example of three dimensional application in\
    \cite{bulle_removing_2020}.
    
    All the numerical results were produced within DOLFIN except the strong
    scaling tests in \cref{sec:strongscaling} which were performed using the
    DOLFINx version of our code.

    We apply different adaptive refinement methods as presented in
    \cref{sec:adaptive}.
    For each method we perform the estimation step with a
    different estimator among the following: $\eta_{\mathrm{res}}$ the residual
    estimator, defined in \cref{app:res}, $\eta_{\mathrm{zz}}$ the
    Zienkiewicz--Zhu estimator, defined in \cref{app:zz}.
    Note that we use the most basic version of the Zienkiewicz--Zhu estimator
    which is not defined for quadratic or cubic
    finite elements nor for linear elasticity problems, and consequently will
    be absent from the comparison in these cases (It is possible to
    extend the idea of the Zienkiewicz--Zhu estimator to higher-order
    polynomials via the definition of the Scott-Zhang interpolator, see\
    \cite{carstensen_axioms_2014, scott_finite_1990}).
    In addition we compare several versions of the Bank--Weiser estimator: the
    bubble Bank--Weiser estimator $\eta_{\mathrm{bw}}^{\mathrm{b}}$ defined from
    the enriched bubble functions space and
    $\eta_{\mathrm{bw}}^{k_+, k_-}$ for multiple choices of the fine and
    coarse spaces orders $k_+$ and $k_-$.

    For each one of the following test cases we will first give a comparison
    of all the refinement strategies by giving the efficiency of the a
    posteriori error estimator on the last mesh of the hierarchy, where the
    efficiency of an estimator $\eta$ is defined as follows:
    \begin{equation}
        \mathrm{eff} := \frac{\eta}{\epsilon_{\mathrm{err}}},
    \end{equation}
    where $\epsilon_{\mathrm{err}}$ is a higher order approximation of the exact
    error computed either from the knowledge of the analytical solution or from
    a higher-order finite element method on a fine mesh.

    \section{Poisson problems}\label{sec:poisson}
    \subsection{L-shaped domain}\label{sec:lshaped}
    We consider a 2D L-shaped domain $\Omega = (-1, 1)^2
    \setminus [-1,0]^2$.
    We solve\ \cref{eq:strong_form} with $f = 0$, $\Gamma_D = \Gamma$, $u_D$
    given by the analytical solution defined below and $\Gamma_N = \varnothing$.
    In polar coordinates, the exact solution is given by
    $u_{\mathrm{exact}} (r, \theta) =
    r^{2/3} \sin \big(2/3 (\theta + \pi/2)\big)$.
    The exact solution belongs to $H^{5/3 - \epsilon}(\Omega)$ for any
    $\epsilon>0$ and its gradient admits a singularity at the vertex of the
    reentrant corner\ \cite[Chapter 5]{grisvard_elliptic}.
    L-shaped domains are widely used to test adaptive mesh refinement
    procedures\ \cite{mitchell_collection_2013}.
    In both linear and quadratic finite elements all the estimators reach
    an expected convergence rate ($\approx -0.5$ in the number of degrees of
    freedom for linear elements and $\approx -1$ for quadratic elements). The
    choice of a posteriori error estimator is not critical for mesh
    refinement purposes, every estimator leading to a hierarchy of meshes on
    which the corresponding errors $\epsilon_{\mathrm{err}}$ are similar. For
    brevity we have not included the convergence plots of these results.\\

\begin{figure}
        \begin{tabular}{c c c}
            \includegraphics[width=0.3\textwidth]{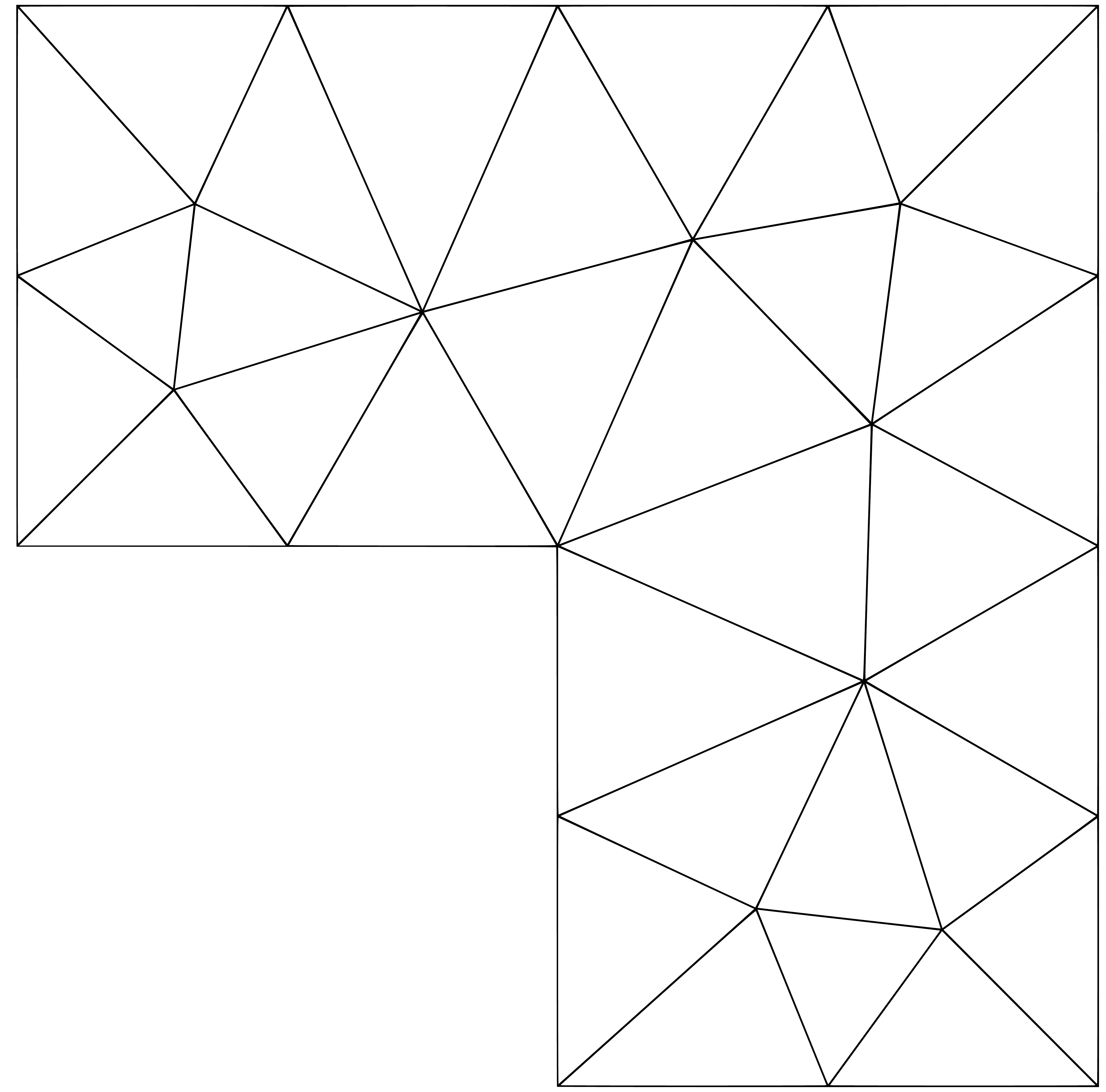}
            &
            \includegraphics[width=0.3\textwidth]{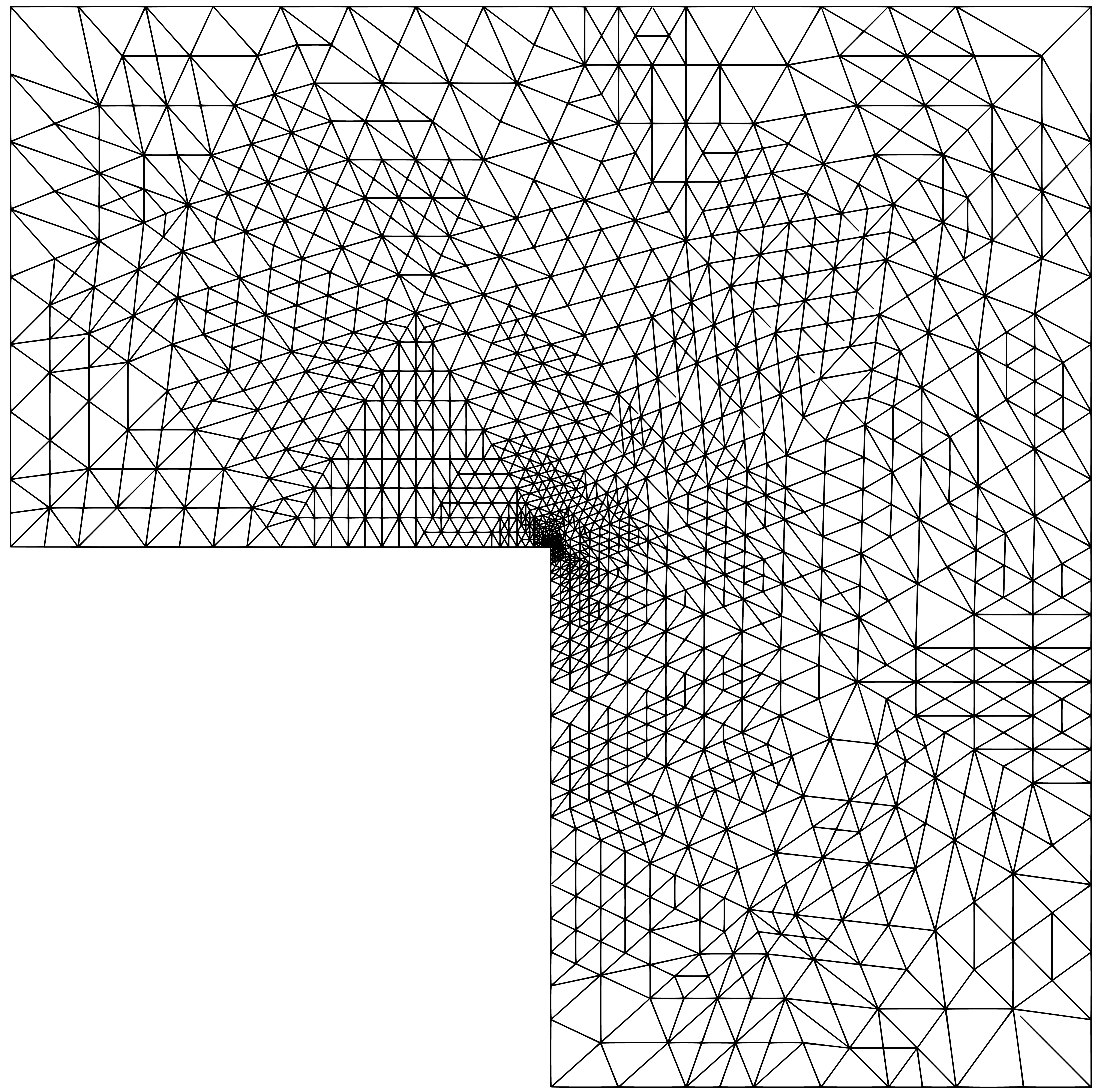}
            &
            \includegraphics[width=0.3\textwidth]{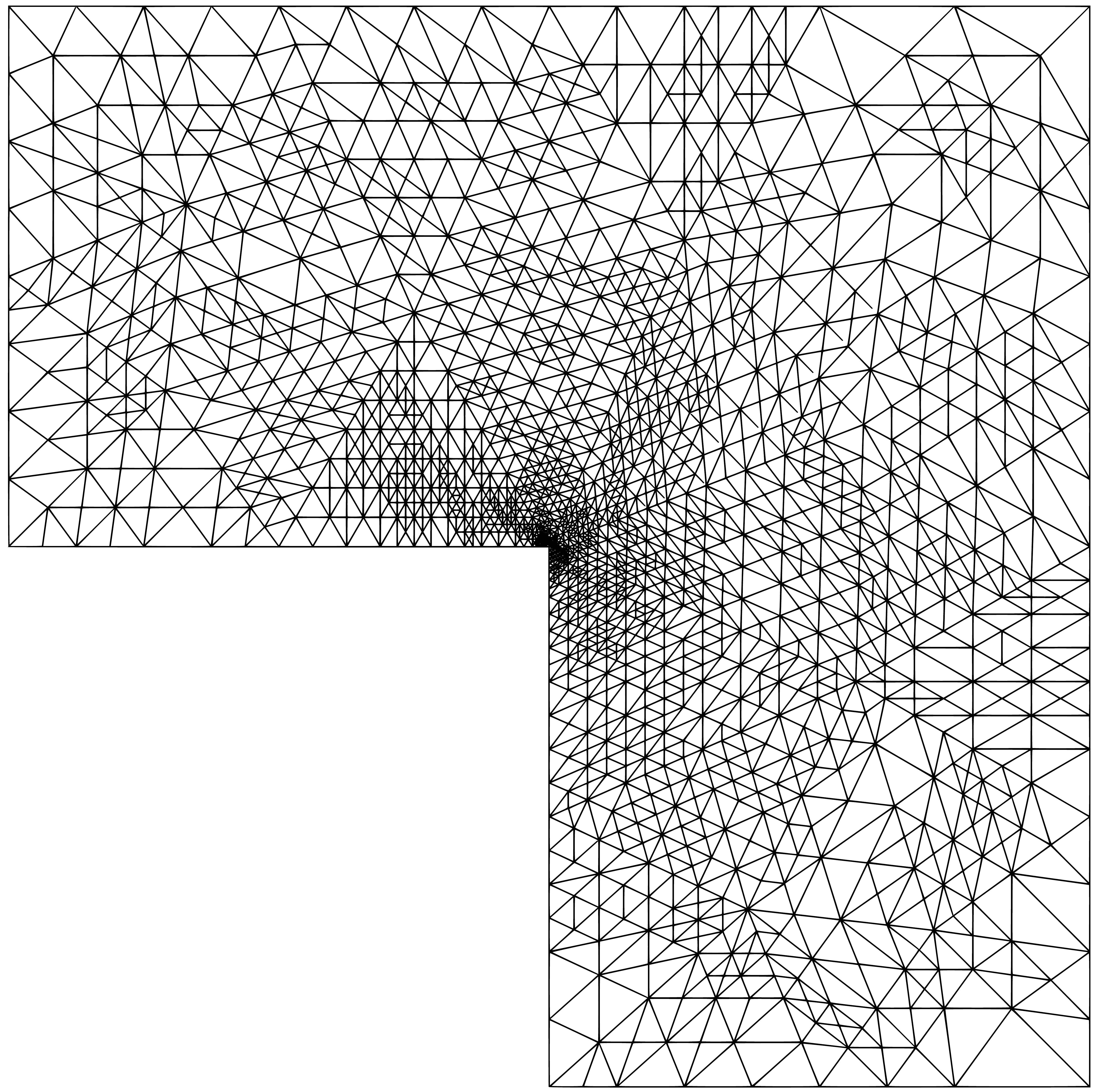}
            \\
            \includegraphics[width=0.3\textwidth]{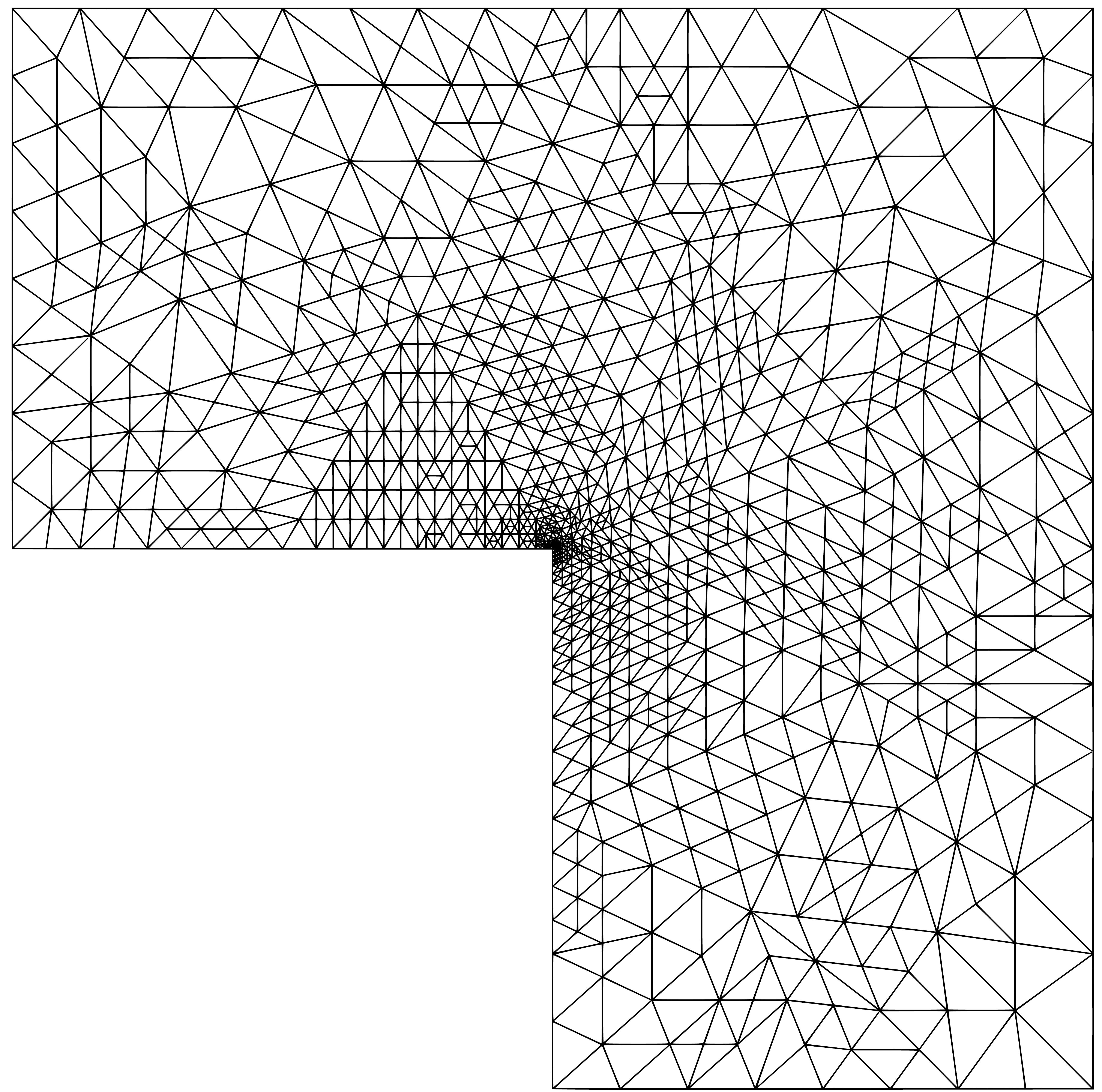}
            &
            \includegraphics[width=0.3\textwidth]{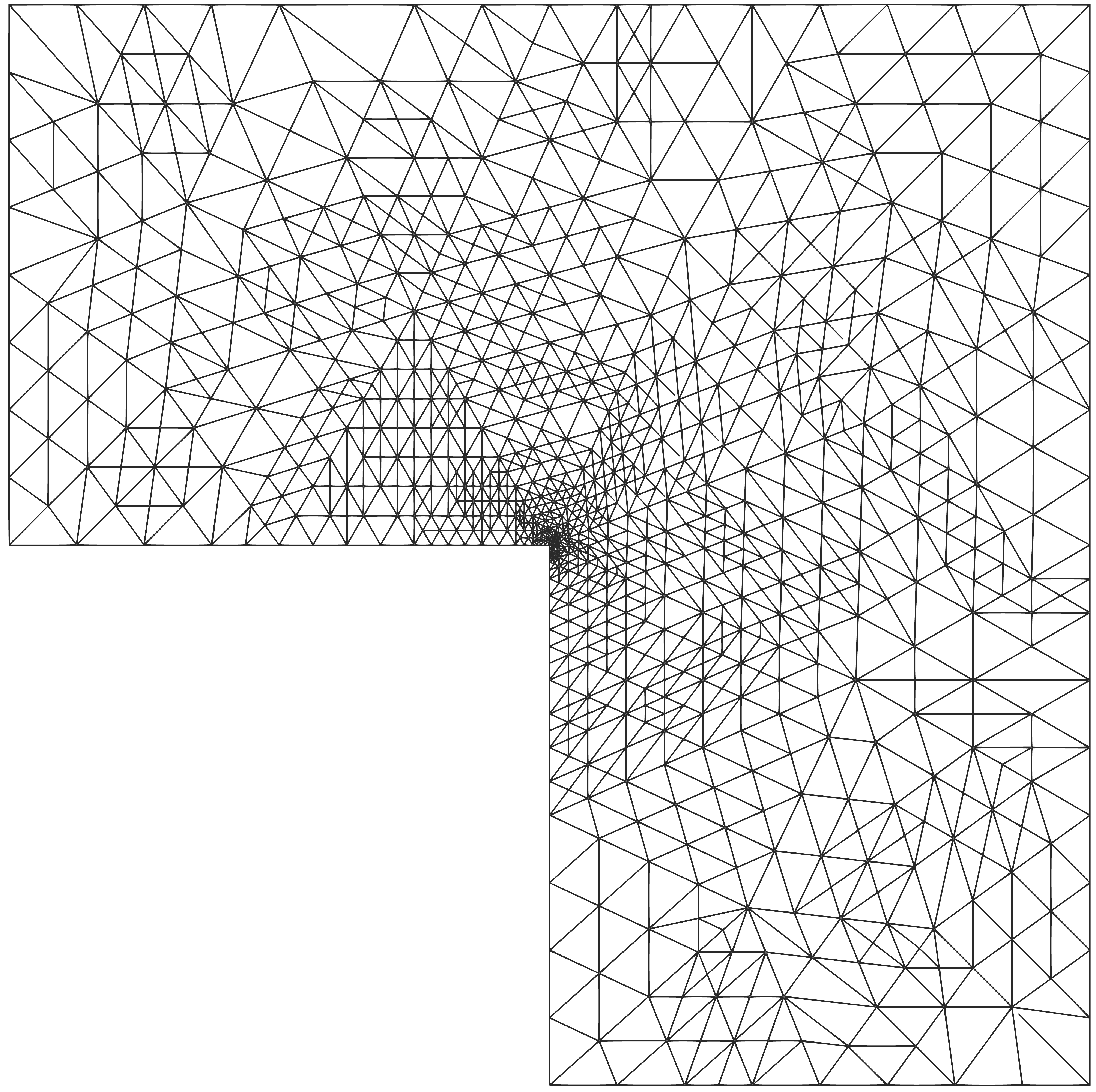}
            &
            \includegraphics[width=0.3\textwidth]{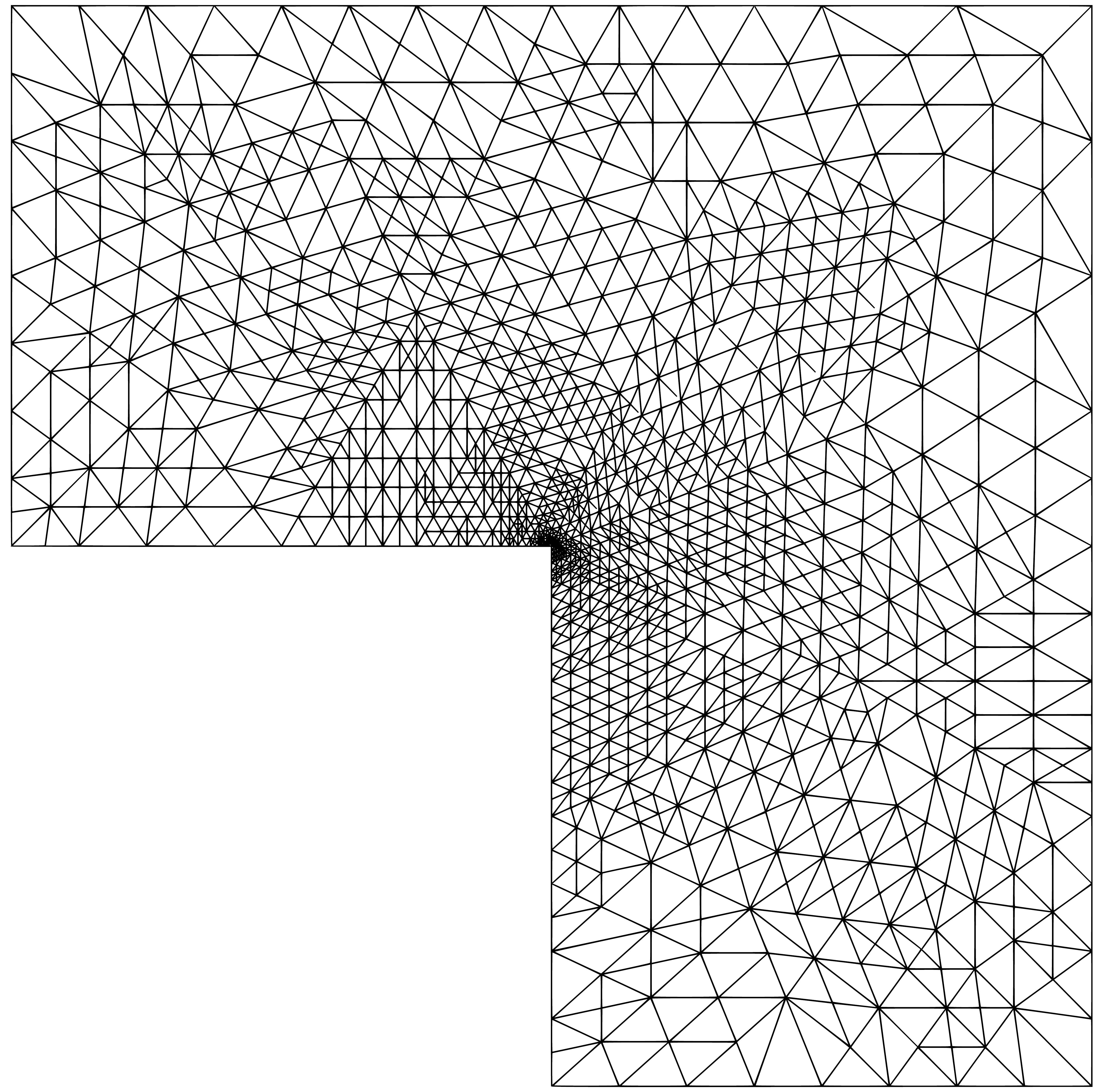}
        \end{tabular}
        \caption{L-shaped Poisson problem with linear elements: On top left the initial mesh used to start all the adaptive
            strategies. From top middle to bottom right, the adaptive meshes
        obtained after seven iterations of refinement strategies steered
        respectively by
        $\eta_{\mathrm{res}}$, $\eta_{\mathrm{bw}}^{\mathrm{b}}$,
        $\eta_{\mathrm{zz}}$, $\eta_{\mathrm{bw}}^{2,1}$ and
        $\eta_{\mathrm{bw}}^{4,2}$}
        \label{fig:meshes}
    \end{figure}

    \textbf{Linear elements.}
    On \cref{fig:meshes} we can see the initial mesh (top left) used to start
    the adaptive refinement strategies. Then, we can see the different refined
    meshes we obtain after seven refinement iterations \comment{, leading to an approximate
    error and number of degrees of freedom given in \cref{fig:table_lshaped}.}\\
    \comment{
    \begin{figure}
        \footnotesize
        \input{./text/L-shaped/tex_table.tex}
        \normalsize
        \caption{Approximate error $\epsilon_{\mathrm{err}}$ and number of degrees
        of freedom for some steps of the adaptive refinement strategies.}
        \label{fig:table_lshaped}
    \end{figure}}
    As we can see on \cref{fig:table_bws_lshaped_P1} the Zienkiewicz--Zhu
    estimator $\eta_{\mathrm{zz}}$ seems to perform the best in terms of
    efficiency while the second best estimator is $\eta_{\bw}^{2,1}$.
    The bubble Bank--Weiser estimator $\eta_{\bw}^{\mathrm b}$ is outperformed
    by almost all the other Bank--Weiser estimators.
    The residual estimator $\eta_{\mathrm{res}}$ largely overestimates the error
    while the estimators $\eta_{\bw}^{k_+, k_-}$ for $k_- > 1$ largely
    underestimates it, leading to poor error approximations.
    Among the poor estimators, $\eta_{\bw}^{3,2}$ is surprisingly off for linear
    elements on this test case.
    This behavior seems to be specific to the L-shaped test cases with linear
    finite elements as we will see below.\\

    \textbf{Quadratic elements.}
    As shown on \cref{fig:table_bws_lshaped_P2}, the best estimator in terms of
    efficiency is $\eta_{\bw}^{2, 0}$ which nearly perfectly matches the error
    $\epsilon_{\mathrm{err}}$.
    We can also notice the very good efficiencies of $\eta_{\bw}^{4,2}$ and
    $\eta_{\bw}^{3,2}$.
    Once again the Bank--Weiser estimators with $k_- > 2$ drastically
    underestimate the error.
    We can notice that the residual estimator is less efficient as the finite
    element degree increases.

    \begin{figure}
    \begin{center}
    \begin{tikzpicture}[x=0.8cm,y=0.6cm] 
    \def\ab{1.34};
    \def\ac{1.53};
    \def\bc{1.22};
    \def\ad{1.53};
    \def\bd{1.53};
    \def\cd{0.0};
    \def\ae{1.59};
    \def\be{1.72};
    \def\ce{0.7};
    \def\de{0.29};
    \def\oa{3.56};
    \def\ob{1.78};
    \def\oc{0.99};
    \draw (0,-1) -- (5, -1) -- (5, 4) -- (0, 4) -- (0, -1);
    \draw (0, 1) -- (5, 1);
    \draw (0, 2) -- (5, 2);
    \draw (0, 3) -- (5, 3);
    \draw (1, -1) -- (1, 4);
    \draw (2, -1) -- (2, 4);
    \draw (3, -1) -- (3, 4);
    \draw (4, -1) -- (4, 4);
    \draw (0, 4) -- (1, 3);
    \draw (0, 0) -- (5, 0);
    \draw (6.01, -1) -- (8.01, -1) -- (8.01, 2) -- (6.01, 2) -- (6.01, -1);
    \draw (6.01, -1) -- (8.01, -1);
    \draw (6.01, 0) -- (8.01, 0);
    \draw (6.01, 1) -- (8.01, 1);
    \draw (7.01, -1) -- (7.01, 2);

    \node at (0.33, 3.25) [text width = 3cm, align=center] {$k_-$};
    \node at (0.75, 3.66) [text width = 3cm, align=center] {$k_+$};
    \node at (0.5, -0.5) [text width = 3cm, align=center] {$3$};
    \node at (0.5, 0.5) [text width = 3cm, align=center] {$2$};
    \node at (0.5, 1.5) [text width = 3cm, align=center] {$1$};
    \node at (0.5, 2.5) [text width = 3cm, align=center] {$0$};

    \node at (1.5, 3.5) [text width = 3cm, align=center] {$1$};
    \node at (2.5, 3.5) [text width = 3cm, align=center] {$2$};
    \node at (3.5, 3.5) [text width = 3cm, align=center] {$3$};
    \node at (4.5, 3.5) [text width = 3cm, align=center] {$4$};

    \node at (1.5, 2.5) [align=center] {\ab};
    \node at (2.5, 2.5) [align=center] {\ac};
    \node at (3.5, 2.5) [align=center] {\ad};
    \node at (4.5, 2.5) [align=center] {\ae};

    \node at (1.5, 1.5) [align=center] {$\varnothing$};
    \node at (2.5, 1.5) [align=center] {\bc};
    \node at (3.5, 1.5) [align=center] {\bd};
    \node at (4.5, 1.5) [align=center] {\be};

    \node at (1.5, 0.5) [align=center] {$\varnothing$};
    \node at (2.5, 0.5) [align=center] {$\varnothing$};
    \node at (3.5, 0.5) [align=center] {\cd};
    \node at (4.5, 0.5) [align=center] {\ce};

    \node at (1.5, -0.5) [align=center] {$\varnothing$};
    \node at (2.5, -0.5) [align=center] {$\varnothing$};
    \node at (3.5, -0.5) [align=center] {$\varnothing$};
    \node at (4.5, -0.5) [align=center] {\de};

    \node at (6.51, 1.5) [text width = 3cm, align=center] {$\eta_{\mathrm{res}}$};
    \node at (6.51, 0.5) [text width = 3cm, align=center] {$\eta^{\mathrm{b}}_{\bw}$};
    \node at (6.51, -0.5) [text width = 3cm, align=center] {$\eta_{\mathrm{zz}}$};
    \node at (7.51, 1.5) [text width = 3cm, align=center] {\oa};
    \node at (7.51, 0.5) [text width = 3cm, align=center] {\ob};
    \node at (7.51, -0.5) [text width = 3cm, align=center] {\oc};
    \end{tikzpicture}
    \end{center}
    \caption{L-shaped Poisson problem with linear elements: efficiencies of $\eta_{\bw}^{k_+,k_-}$ and other estimators on the
    last mesh of an adaptively refined hierarchy.}
    \label{fig:table_bws_lshaped_P1}
    \end{figure}

    \begin{figure}
    \begin{center}
    \begin{tikzpicture}[x=0.8cm,y=0.6cm] 
    \def\ab{0.66};
    \def\ac{1.0};
    \def\bc{1.61};
    \def\ad{1.12};
    \def\bd{2.1};
    \def\cd{0.92};
    \def\ae{1.27};
    \def\be{2.28};
    \def\ce{1.07};
    \def\de{0.31};
    \def\oa{8.67};
    \def\ob{1.84};
    \draw (0,-1) -- (5, -1) -- (5, 4) -- (0, 4) -- (0, -1);
    \draw (0, 1) -- (5, 1);
    \draw (0, 2) -- (5, 2);
    \draw (0, 3) -- (5, 3);
    \draw (1, -1) -- (1, 4);
    \draw (2, -1) -- (2, 4);
    \draw (3, -1) -- (3, 4);
    \draw (4, -1) -- (4, 4);
    \draw (0, 4) -- (1, 3);
    \draw (0, 0) -- (5, 0);
    \draw (6.01, -1) -- (8.01, -1) -- (8.01, 2) -- (6.01, 2) -- (6.01, -1);
    \draw (6.01, -1) -- (8.01, -1);
    \draw (6.01, 0) -- (8.01, 0);
    \draw (6.01, 1) -- (8.01, 1);
    \draw (7.01, -1) -- (7.01, 2);

    \node at (0.33, 3.25) [text width = 3cm, align=center] {$k_-$};
    \node at (0.75, 3.66) [text width = 3cm, align=center] {$k_+$};
    \node at (0.5, -0.5) [text width = 3cm, align=center] {$3$};
    \node at (0.5, 0.5) [text width = 3cm, align=center] {$2$};
    \node at (0.5, 1.5) [text width = 3cm, align=center] {$1$};
    \node at (0.5, 2.5) [text width = 3cm, align=center] {$0$};

    \node at (1.5, 3.5) [text width = 3cm, align=center] {$1$};
    \node at (2.5, 3.5) [text width = 3cm, align=center] {$2$};
    \node at (3.5, 3.5) [text width = 3cm, align=center] {$3$};
    \node at (4.5, 3.5) [text width = 3cm, align=center] {$4$};

    \node at (1.5, 2.5) [align=center] {\ab};
    \node at (2.5, 2.5) [align=center] {\ac};
    \node at (3.5, 2.5) [align=center] {\ad};
    \node at (4.5, 2.5) [align=center] {\ae};

    \node at (1.5, 1.5) [align=center] {$\varnothing$};
    \node at (2.5, 1.5) [align=center] {\bc};
    \node at (3.5, 1.5) [align=center] {\bd};
    \node at (4.5, 1.5) [align=center] {\be};

    \node at (1.5, 0.5) [align=center] {$\varnothing$};
    \node at (2.5, 0.5) [align=center] {$\varnothing$};
    \node at (3.5, 0.5) [align=center] {\cd};
    \node at (4.5, 0.5) [align=center] {\ce};

    \node at (1.5, -0.5) [align=center] {$\varnothing$};
    \node at (2.5, -0.5) [align=center] {$\varnothing$};
    \node at (3.5, -0.5) [align=center] {$\varnothing$};
    \node at (4.5, -0.5) [align=center] {\de};

    \node at (6.51, 1.5) [text width = 3cm, align=center] {$\eta_{\mathrm{res}}$};
    \node at (6.51, 0.5) [text width = 3cm, align=center] {$\eta^{\mathrm{b}}_{\bw}$};
    \node at (6.51, -0.5) [text width = 3cm, align=center] {$\eta_{\mathrm{zz}}$};
    \node at (7.51, 1.5) [text width = 3cm, align=center] {\oa};
    \node at (7.51, 0.5) [text width = 3cm, align=center] {\ob};
    \node at (7.51, -0.5) [text width = 3cm, align=center] {$\varnothing$};
    \end{tikzpicture}
    \end{center}
    \caption{L-shaped Poisson problem with quadratic elements: efficiencies of $\eta_{\bw}^{k_+,k_-}$ and other estimators on the
    last mesh of an adaptively refined hierarchy.}
    \label{fig:table_bws_lshaped_P2}
    \end{figure}

    \subsection{Mixed boundary conditions L-shaped
    domain}\label{sec:mix_bc_lshaped} We solve \cref{eq:strong_form} on the same
    two-dimensional L-shaped boundary domain as in \cref{sec:lshaped} but with
    different boundary conditions.
    We consider $f = 0$, $\Gamma_N = \{(x,y) \in \R^2,\ x<0,\ y=0\}$ and
    $\Gamma_D = \Gamma \setminus \Gamma_N$.
    The boundary data are given by $g = 0$ and $u_D = u_{\mathrm{exact}} =
    r^{1/3} \sin \big(1/3(\theta + \pi/2)\big)$. 
    The exact solution belongs to $H^{4/3 - \epsilon}(\Omega)$ for any $\epsilon
    > 0$ and its gradient has a singularity located at the reentrant corner of
    $\Gamma$ (see\ \cite[Chapter 5]{grisvard_elliptic}).
    As before, each estimator is leading to a convergence rate close to the
    expected one ($\approx -0.5$ for linear elements, $\approx -1$ for quadratic
    elements) and the choice of the estimator does not impact the quality of the
    mesh hierarchy.\\

    \textbf{Linear elements.} 
    First thing we can notice from \cref{fig:table_bws_mixed_bc_lshaped_P1} is
    that the estimators efficiencies are quite different from those in
    \cref{fig:table_bws_lshaped_P1}.
    Most of the Bank--Weiser estimator efficiencies have improved, except when
    $k_- > 1$.
    The Zienkiewicz--Zhu estimator $\eta_{\mathrm{zz}}$ is no longer the most
    efficient and has been outperformed by $\eta_{\bw}^{2,0}$,
    $\eta_{\bw}^{2,1}$ and $\eta_{\bw}^{3,0}$.
    The Bank--Weiser estimator $\eta_{\bw}^{3,2}$ still performs poorly as in
    \cref{fig:table_bws_lshaped_P1}, while the residual estimator
    $\eta_{\mathrm{res}}$ once again largely overestimates the error.\\

    \textbf{Quadratic elements.}
    As for linear elements, the efficiencies in
    \cref{fig:table_bws_mixed_bc_lshaped_P2} are very different from
    \cref{fig:table_bws_lshaped_P2}, many Bank--Weiser estimators are now
    underestimating the error.
    The most efficient estimator is $\eta_{\bw}^{2,1}$ closely followed by the
    bubble Bank--Weiser estimator $\eta_{\bw}^{\mathrm{b}}$.
    As for the previous test cases, the Bank--Weiser estimators with $k_- > 2$
    are largely underestimating the error.

    \begin{figure}
    \begin{center}
    \begin{tikzpicture}[x=0.8cm,y=0.6cm] 
    \def\ab{0.83};
    \def\ac{1.06};
    \def\bc{0.94};
    \def\ad{1.08};
    \def\bd{1.21};
    \def\cd{0.0};
    \def\ae{1.14};
    \def\be{1.34};
    \def\ce{0.55};
    \def\de{0.23};
    \def\oa{2.84};
    \def\ob{1.24};
    \def\oc{0.91};
    \draw (0,-1) -- (5, -1) -- (5, 4) -- (0, 4) -- (0, -1);
    \draw (0, 1) -- (5, 1);
    \draw (0, 2) -- (5, 2);
    \draw (0, 3) -- (5, 3);
    \draw (1, -1) -- (1, 4);
    \draw (2, -1) -- (2, 4);
    \draw (3, -1) -- (3, 4);
    \draw (4, -1) -- (4, 4);
    \draw (0, 4) -- (1, 3);
    \draw (0, 0) -- (5, 0);
    \draw (6.01, -1) -- (8.01, -1) -- (8.01, 2) -- (6.01, 2) -- (6.01, -1);
    \draw (6.01, -1) -- (8.01, -1);
    \draw (6.01, 0) -- (8.01, 0);
    \draw (6.01, 1) -- (8.01, 1);
    \draw (7.01, -1) -- (7.01, 2);

    \node at (0.33, 3.25) [text width = 3cm, align=center] {$k_-$};
    \node at (0.75, 3.66) [text width = 3cm, align=center] {$k_+$};
    \node at (0.5, -0.5) [text width = 3cm, align=center] {$3$};
    \node at (0.5, 0.5) [text width = 3cm, align=center] {$2$};
    \node at (0.5, 1.5) [text width = 3cm, align=center] {$1$};
    \node at (0.5, 2.5) [text width = 3cm, align=center] {$0$};

    \node at (1.5, 3.5) [text width = 3cm, align=center] {$1$};
    \node at (2.5, 3.5) [text width = 3cm, align=center] {$2$};
    \node at (3.5, 3.5) [text width = 3cm, align=center] {$3$};
    \node at (4.5, 3.5) [text width = 3cm, align=center] {$4$};

    \node at (1.5, 2.5) [align=center] {\ab};
    \node at (2.5, 2.5) [align=center] {\ac};
    \node at (3.5, 2.5) [align=center] {\ad};
    \node at (4.5, 2.5) [align=center] {\ae};

    \node at (1.5, 1.5) [align=center] {$\varnothing$};
    \node at (2.5, 1.5) [align=center] {\bc};
    \node at (3.5, 1.5) [align=center] {\bd};
    \node at (4.5, 1.5) [align=center] {\be};

    \node at (1.5, 0.5) [align=center] {$\varnothing$};
    \node at (2.5, 0.5) [align=center] {$\varnothing$};
    \node at (3.5, 0.5) [align=center] {\cd};
    \node at (4.5, 0.5) [align=center] {\ce};

    \node at (1.5, -0.5) [align=center] {$\varnothing$};
    \node at (2.5, -0.5) [align=center] {$\varnothing$};
    \node at (3.5, -0.5) [align=center] {$\varnothing$};
    \node at (4.5, -0.5) [align=center] {\de};

    \node at (6.51, 1.5) [text width = 3cm, align=center] {$\eta_{\mathrm{res}}$};
    \node at (6.51, 0.5) [text width = 3cm, align=center] {$\eta^{\mathrm{b}}_{\bw}$};
    \node at (6.51, -0.5) [text width = 3cm, align=center] {$\eta_{\mathrm{zz}}$};
    \node at (7.51, 1.5) [text width = 3cm, align=center] {\oa};
    \node at (7.51, 0.5) [text width = 3cm, align=center] {\ob};
    \node at (7.51, -0.5) [text width = 3cm, align=center] {\oc};
    \end{tikzpicture}
    \end{center}
    \caption{Mixed boundary conditions L-shaped Poisson problem with linear
        elements: efficiencies of $\eta_{\bw}^{k_+,k_-}$ and other estimators on the
    last mesh of an adaptively refined hierarchy.}
    \label{fig:table_bws_mixed_bc_lshaped_P1}
    \end{figure}

    \begin{figure}
    \begin{center}
    \begin{tikzpicture}[x=0.8cm,y=0.6cm] 
    \def\ab{0.57};
    \def\ac{0.86};
    \def\bc{1.05};
    \def\ad{0.83};
    \def\bd{1.37};
    \def\cd{0.62};
    \def\ae{0.97};
    \def\be{1.43};
    \def\ce{0.78};
    \def\de{0.3};
    \def\oa{5.91};
    \def\ob{1.17};
    \draw (0,-1) -- (5, -1) -- (5, 4) -- (0, 4) -- (0, -1);
    \draw (0, 1) -- (5, 1);
    \draw (0, 2) -- (5, 2);
    \draw (0, 3) -- (5, 3);
    \draw (1, -1) -- (1, 4);
    \draw (2, -1) -- (2, 4);
    \draw (3, -1) -- (3, 4);
    \draw (4, -1) -- (4, 4);
    \draw (0, 4) -- (1, 3);
    \draw (0, 0) -- (5, 0);
    \draw (6.01, -1) -- (8.01, -1) -- (8.01, 2) -- (6.01, 2) -- (6.01, -1);
    \draw (6.01, -1) -- (8.01, -1);
    \draw (6.01, 0) -- (8.01, 0);
    \draw (6.01, 1) -- (8.01, 1);
    \draw (7.01, -1) -- (7.01, 2);

    \node at (0.33, 3.25) [text width = 3cm, align=center] {$k_-$};
    \node at (0.75, 3.66) [text width = 3cm, align=center] {$k_+$};
    \node at (0.5, -0.5) [text width = 3cm, align=center] {$3$};
    \node at (0.5, 0.5) [text width = 3cm, align=center] {$2$};
    \node at (0.5, 1.5) [text width = 3cm, align=center] {$1$};
    \node at (0.5, 2.5) [text width = 3cm, align=center] {$0$};

    \node at (1.5, 3.5) [text width = 3cm, align=center] {$1$};
    \node at (2.5, 3.5) [text width = 3cm, align=center] {$2$};
    \node at (3.5, 3.5) [text width = 3cm, align=center] {$3$};
    \node at (4.5, 3.5) [text width = 3cm, align=center] {$4$};

    \node at (1.5, 2.5) [align=center] {\ab};
    \node at (2.5, 2.5) [align=center] {\ac};
    \node at (3.5, 2.5) [align=center] {\ad};
    \node at (4.5, 2.5) [align=center] {\ae};

    \node at (1.5, 1.5) [align=center] {$\varnothing$};
    \node at (2.5, 1.5) [align=center] {\bc};
    \node at (3.5, 1.5) [align=center] {\bd};
    \node at (4.5, 1.5) [align=center] {\be};

    \node at (1.5, 0.5) [align=center] {$\varnothing$};
    \node at (2.5, 0.5) [align=center] {$\varnothing$};
    \node at (3.5, 0.5) [align=center] {\cd};
    \node at (4.5, 0.5) [align=center] {\ce};

    \node at (1.5, -0.5) [align=center] {$\varnothing$};
    \node at (2.5, -0.5) [align=center] {$\varnothing$};
    \node at (3.5, -0.5) [align=center] {$\varnothing$};
    \node at (4.5, -0.5) [align=center] {\de};

    \node at (6.51, 1.5) [text width = 3cm, align=center] {$\eta_{\mathrm{res}}$};
    \node at (6.51, 0.5) [text width = 3cm, align=center] {$\eta^{\mathrm{b}}_{\bw}$};
    \node at (6.51, -0.5) [text width = 3cm, align=center] {$\eta_{\mathrm{zz}}$};
    \node at (7.51, 1.5) [text width = 3cm, align=center] {\oa};
    \node at (7.51, 0.5) [text width = 3cm, align=center] {\ob};
    \node at (7.51, -0.5) [text width = 3cm, align=center] {$\varnothing$};
    \end{tikzpicture}
    \end{center}
    \caption{Mixed boundary conditions L-shaped Poisson problem with quadratic
        elements: efficiencies of $\eta_{\bw}^{k_+,k_-}$ and other estimators on the
    last mesh of an adaptively refined hierarchy.}
    \label{fig:table_bws_mixed_bc_lshaped_P2}
    \end{figure}

    \subsection{Boundary singularity}\label{sec:boundary_singularity} We
    solve \cref{eq:strong_form} on a two-dimensional unit square domain $\Omega
    = (0,1)^2$ with $u = u_{\mathrm{exact}}$ on $\Gamma_D = \Gamma$, ($\Gamma_N
    = \varnothing$) and $f$ chosen in order to have $u(x,y) =
    u_{\mathrm{exact}}(x,y) = x^{\alpha},$ with $\alpha \geqslant 0 .5$.
    In the following results we chose $\alpha = 0.7$.
    The gradient of the exact solution $u$ admits a singularity along the left
    boundary of $\Omega$ (for $x = 0$). 
    The solution $u$ belongs to $H^{6/5-\epsilon}$ for all $\epsilon >0$\
    \cite{houston_sobolev,mitchell_collection_2017}.
    Consequently, the value of $\alpha$ determines the strength of the
    singularity and the regularity of $u$.
   
    Due to the presence of the edge singularity, all the estimators are achieving a
    convergence rate close to $-0.2$ for linear elements.
    Moreover, this rate does not improve for higher-order elements (for brevity,
    the results for higher-order elements are not shown here).
    The low convergence rate shows how computationally challenging such a
    problem can be.
    Once again the choice of estimator is not critical for mesh
    refinement purposes.\\

    \textbf{Linear elements.}
    The best estimator in terms of efficiency is $\eta_{\bw}^{2,1}$ which
    slightly overestimates the error, closely followed by $\eta_{\bw}^{4,2}$
    underestimating the error as we can see on
    \cref{fig:table_bws_boundary_singularity_P1}.
    Unlike the previous test case, here the Zienkiewicz--Zhu estimator
    $\eta_{\mathrm{zz}}$ grandly underestimates the error.
    The worst estimator is the residual estimator $\eta_{\mathrm{res}}$ which
    gives no precise information about the error.
    We can notice that the poor performance of the estimator
    $\eta_{\bw}^{3,2}$ on the L-shaped test case does not reproduce here.\\

    \textbf{Quadratic elements.}
    Again, \cref{fig:table_bws_boundary_singularity_P2} shows that the best
    estimator is $\eta_{\bw}^{2,1}$ closely followed by $\eta_{\bw}^{2,0}$ and
    the bubble estimator $\eta_{\bw}^{\mathrm{b}}$.
    The residual estimator is getting worse as the finite element degree
    increases.

    \begin{figure}
    \begin{center}
    \begin{tikzpicture}[x=0.8cm,y=0.6cm] 
    \def\ab{0.74};
    \def\ac{1.15};
    \def\bc{1.06};
    \def\ad{0.94};
    \def\bd{1.27};
    \def\cd{0.72};
    \def\ae{1.1};
    \def\be{1.41};
    \def\ce{0.95};
    \def\de{0.62};
    \def\oa{17.02};
    \def\ob{1.23};
    \def\oc{0.6};
    \draw (0,-1) -- (5, -1) -- (5, 4) -- (0, 4) -- (0, -1);
    \draw (0, 1) -- (5, 1);
    \draw (0, 2) -- (5, 2);
    \draw (0, 3) -- (5, 3);
    \draw (1, -1) -- (1, 4);
    \draw (2, -1) -- (2, 4);
    \draw (3, -1) -- (3, 4);
    \draw (4, -1) -- (4, 4);
    \draw (0, 4) -- (1, 3);
    \draw (0, 0) -- (5, 0);
    \draw (6.01, -1) -- (8.01, -1) -- (8.01, 2) -- (6.01, 2) -- (6.01, -1);
    \draw (6.01, -1) -- (8.01, -1);
    \draw (6.01, 0) -- (8.01, 0);
    \draw (6.01, 1) -- (8.01, 1);
    \draw (7.01, -1) -- (7.01, 2);

    \node at (0.33, 3.25) [text width = 3cm, align=center] {$k_-$};
    \node at (0.75, 3.66) [text width = 3cm, align=center] {$k_+$};
    \node at (0.5, -0.5) [text width = 3cm, align=center] {$3$};
    \node at (0.5, 0.5) [text width = 3cm, align=center] {$2$};
    \node at (0.5, 1.5) [text width = 3cm, align=center] {$1$};
    \node at (0.5, 2.5) [text width = 3cm, align=center] {$0$};

    \node at (1.5, 3.5) [text width = 3cm, align=center] {$1$};
    \node at (2.5, 3.5) [text width = 3cm, align=center] {$2$};
    \node at (3.5, 3.5) [text width = 3cm, align=center] {$3$};
    \node at (4.5, 3.5) [text width = 3cm, align=center] {$4$};

    \node at (1.5, 2.5) [align=center] {\ab};
    \node at (2.5, 2.5) [align=center] {\ac};
    \node at (3.5, 2.5) [align=center] {\ad};
    \node at (4.5, 2.5) [align=center] {\ae};

    \node at (1.5, 1.5) [align=center] {$\varnothing$};
    \node at (2.5, 1.5) [align=center] {\bc};
    \node at (3.5, 1.5) [align=center] {\bd};
    \node at (4.5, 1.5) [align=center] {\be};

    \node at (1.5, 0.5) [align=center] {$\varnothing$};
    \node at (2.5, 0.5) [align=center] {$\varnothing$};
    \node at (3.5, 0.5) [align=center] {\cd};
    \node at (4.5, 0.5) [align=center] {\ce};

    \node at (1.5, -0.5) [align=center] {$\varnothing$};
    \node at (2.5, -0.5) [align=center] {$\varnothing$};
    \node at (3.5, -0.5) [align=center] {$\varnothing$};
    \node at (4.5, -0.5) [align=center] {\de};

    \node at (6.51, 1.5) [text width = 3cm, align=center] {$\eta_{\mathrm{res}}$};
    \node at (6.51, 0.5) [text width = 3cm, align=center] {$\eta^{\mathrm{b}}_{\bw}$};
    \node at (6.51, -0.5) [text width = 3cm, align=center] {$\eta_{\mathrm{zz}}$};
    \node at (7.51, 1.5) [text width = 3cm, align=center] {\oa};
    \node at (7.51, 0.5) [text width = 3cm, align=center] {\ob};
    \node at (7.51, -0.5) [text width = 3cm, align=center] {\oc};
    \end{tikzpicture}
    \end{center}
    \caption{Boundary singularity Poisson problem with linear elements: efficiencies of $\eta_{\bw}^{k_+,k_-}$ and other estimators on the
    last mesh of an adaptively refined hierarchy.}
    \label{fig:table_bws_boundary_singularity_P1}
    \end{figure}

    \begin{figure}
    \begin{center}
    \begin{tikzpicture}[x=0.8cm,y=0.6cm] 
    \def\ab{0.46};
    \def\ac{0.91};
    \def\bc{0.96};
    \def\ad{1.29};
    \def\bd{1.54};
    \def\cd{1.2};
    \def\ae{1.51};
    \def\be{1.7};
    \def\ce{1.4};
    \def\de{1.16};
    \def\oa{37.61};
    \def\ob{1.13};
    \draw (0,-1) -- (5, -1) -- (5, 4) -- (0, 4) -- (0, -1);
    \draw (0, 1) -- (5, 1);
    \draw (0, 2) -- (5, 2);
    \draw (0, 3) -- (5, 3);
    \draw (1, -1) -- (1, 4);
    \draw (2, -1) -- (2, 4);
    \draw (3, -1) -- (3, 4);
    \draw (4, -1) -- (4, 4);
    \draw (0, 4) -- (1, 3);
    \draw (0, 0) -- (5, 0);
    \draw (6.01, -1) -- (8.01, -1) -- (8.01, 2) -- (6.01, 2) -- (6.01, -1);
    \draw (6.01, -1) -- (8.01, -1);
    \draw (6.01, 0) -- (8.01, 0);
    \draw (6.01, 1) -- (8.01, 1);
    \draw (7.01, -1) -- (7.01, 2);

    \node at (0.33, 3.25) [text width = 3cm, align=center] {$k_-$};
    \node at (0.75, 3.66) [text width = 3cm, align=center] {$k_+$};
    \node at (0.5, -0.5) [text width = 3cm, align=center] {$3$};
    \node at (0.5, 0.5) [text width = 3cm, align=center] {$2$};
    \node at (0.5, 1.5) [text width = 3cm, align=center] {$1$};
    \node at (0.5, 2.5) [text width = 3cm, align=center] {$0$};

    \node at (1.5, 3.5) [text width = 3cm, align=center] {$1$};
    \node at (2.5, 3.5) [text width = 3cm, align=center] {$2$};
    \node at (3.5, 3.5) [text width = 3cm, align=center] {$3$};
    \node at (4.5, 3.5) [text width = 3cm, align=center] {$4$};

    \node at (1.5, 2.5) [align=center] {\ab};
    \node at (2.5, 2.5) [align=center] {\ac};
    \node at (3.5, 2.5) [align=center] {\ad};
    \node at (4.5, 2.5) [align=center] {\ae};

    \node at (1.5, 1.5) [align=center] {$\varnothing$};
    \node at (2.5, 1.5) [align=center] {\bc};
    \node at (3.5, 1.5) [align=center] {\bd};
    \node at (4.5, 1.5) [align=center] {\be};

    \node at (1.5, 0.5) [align=center] {$\varnothing$};
    \node at (2.5, 0.5) [align=center] {$\varnothing$};
    \node at (3.5, 0.5) [align=center] {\cd};
    \node at (4.5, 0.5) [align=center] {\ce};

    \node at (1.5, -0.5) [align=center] {$\varnothing$};
    \node at (2.5, -0.5) [align=center] {$\varnothing$};
    \node at (3.5, -0.5) [align=center] {$\varnothing$};
    \node at (4.5, -0.5) [align=center] {\de};

    \node at (6.51, 1.5) [text width = 3cm, align=center] {$\eta_{\mathrm{res}}$};
    \node at (6.51, 0.5) [text width = 3cm, align=center] {$\eta^{\mathrm{b}}_{\bw}$};
    \node at (6.51, -0.5) [text width = 3cm, align=center] {$\eta_{\mathrm{zz}}$};
    \node at (7.51, 1.5) [text width = 3cm, align=center] {\oa};
    \node at (7.51, 0.5) [text width = 3cm, align=center] {\ob};
    \node at (7.51, -0.5) [text width = 3cm, align=center] {$\varnothing$};
    \end{tikzpicture}
    \end{center}
    \caption{Boundary singularity Poisson problem with quadratic elements: efficiencies of $\eta_{\bw}^{k_+,k_-}$ and other estimators on the
    last mesh of an adaptively refined hierarchy.}
    \label{fig:table_bws_boundary_singularity_P2}
    \end{figure}

    \subsection{Goal-oriented adaptive refinement using linear elements}\label{sec:go_adaptive}
    We solve the L-shaped domain problem as described in \cref{sec:lshaped}
    but instead of controlling the error in the natural norm, we aim to control
    the error in the goal functional $J(u) = (c, u)$ with $c$ a smooth bump function
    \begin{equation}
        c(\bar{r}) :=
        \left\lbrace
        \begin{array}{l l}
            \epsilon^{-2} \exp \left( - \frac{1}{\bar{r}^2} \right) & 0 \le \bar{r}^2 < 1, \\
            0 & \bar{r}^2 \ge 1.
        \end{array}\right.
    \end{equation}
    where $\bar{r}^2 = \left( \left( x - \bar{x} \right)/\epsilon \right)^2 +
        \left( \left( y - \bar{y}\right) /\epsilon \right)^2$, with $\epsilon
        \in \R$ a parameter that controls the size of the bump
    function, and $\bar{x} \in \R$ and $\bar{y} \in R$ the position of the
    bumps function's center.
    We set $\epsilon = 0.35$ and $\bar{x} = \bar{y} = 0.2$.
    With these parameters the goal functional is isolated to a region close
    to the re-entrant corner.

    We use the goal-oriented adaptive mesh refinement methodology outlined in
    \cref{sec:wgo}.
    We use a first-order polynomial finite element method for the primal and
    dual problem, and the Bank--Weiser error estimation procedure to
    calculate both $\eta_u$ and $\eta_z$.

    The `exact' value of the functional $J(u)$ was calculated on a very fine
    mesh using a fourth-order polynomial finite element space and was used to
    compute higher-order approximate errors for each refinement strategy.

    The weighted goal-oriented strategy refines both the re-entrant corner and
    the broader region of interest defined by the goal functional.
    Relatively less refinement occurs in the regions far away from either of
    these important areas.
    \begin{figure}
        \begin{tabular}{c c c}
            \includegraphics[width=0.3\textwidth]{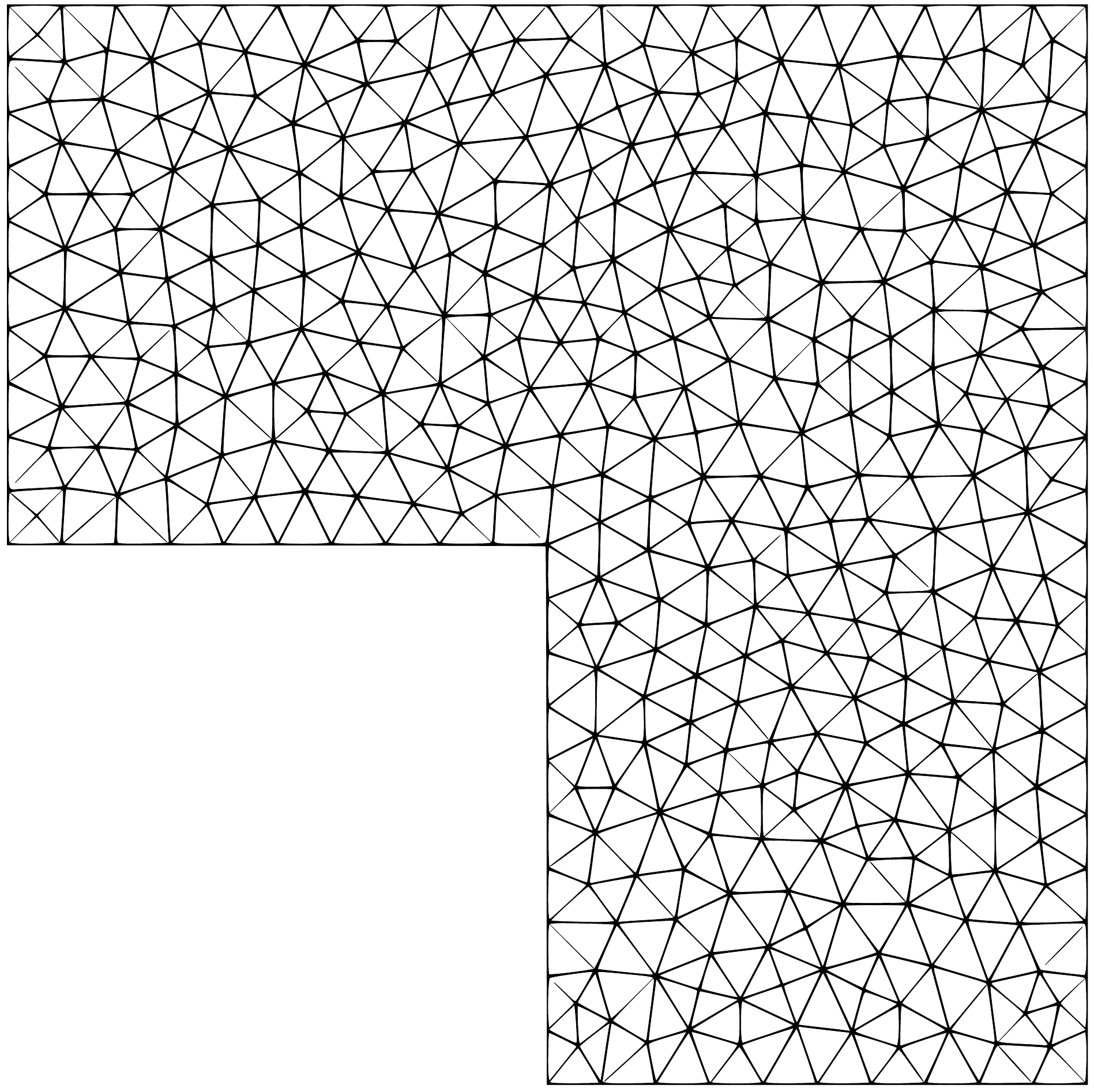}
            &
            \includegraphics[width=0.3\textwidth]{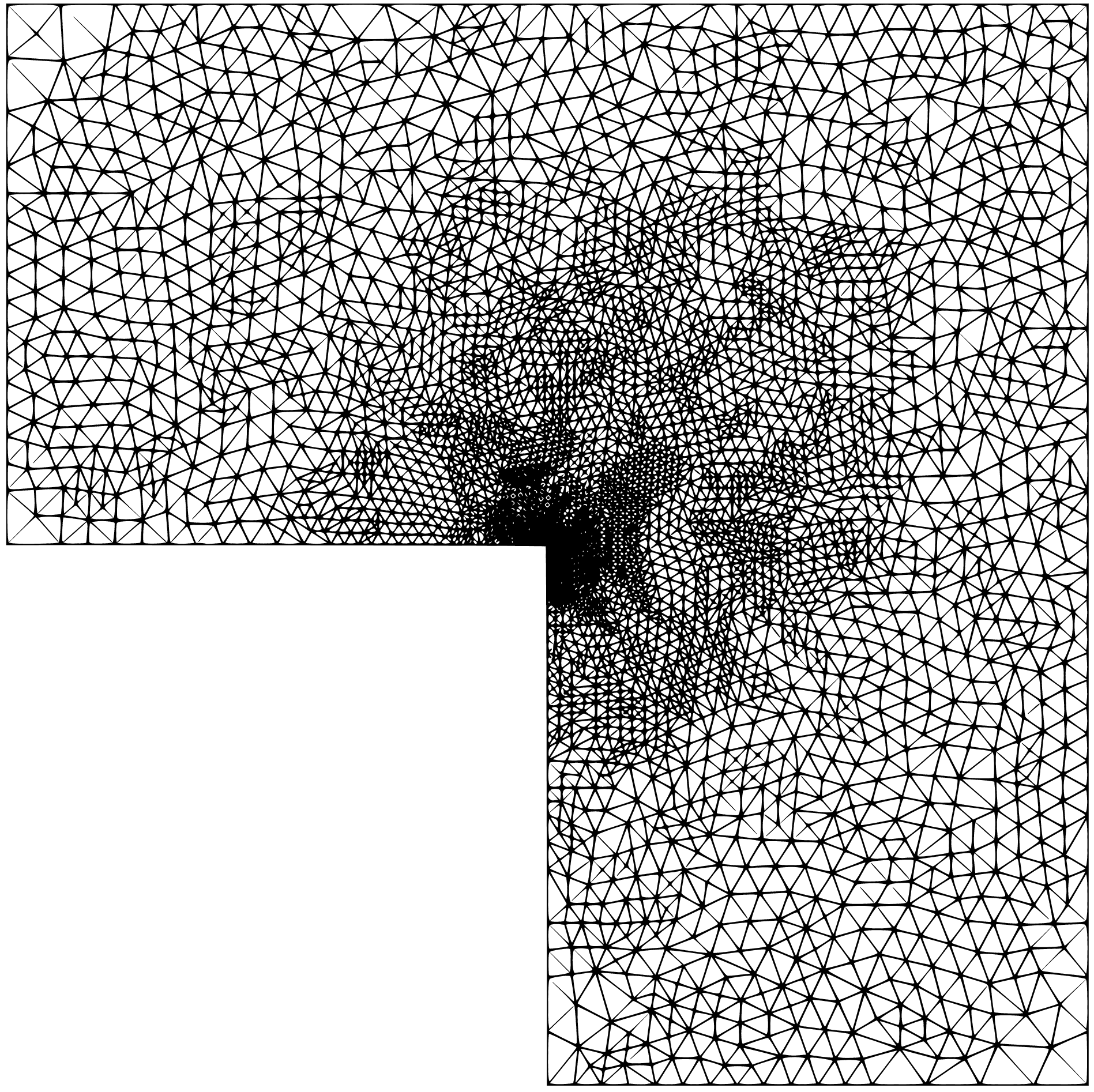}
            &
            \includegraphics[width=0.3\textwidth]{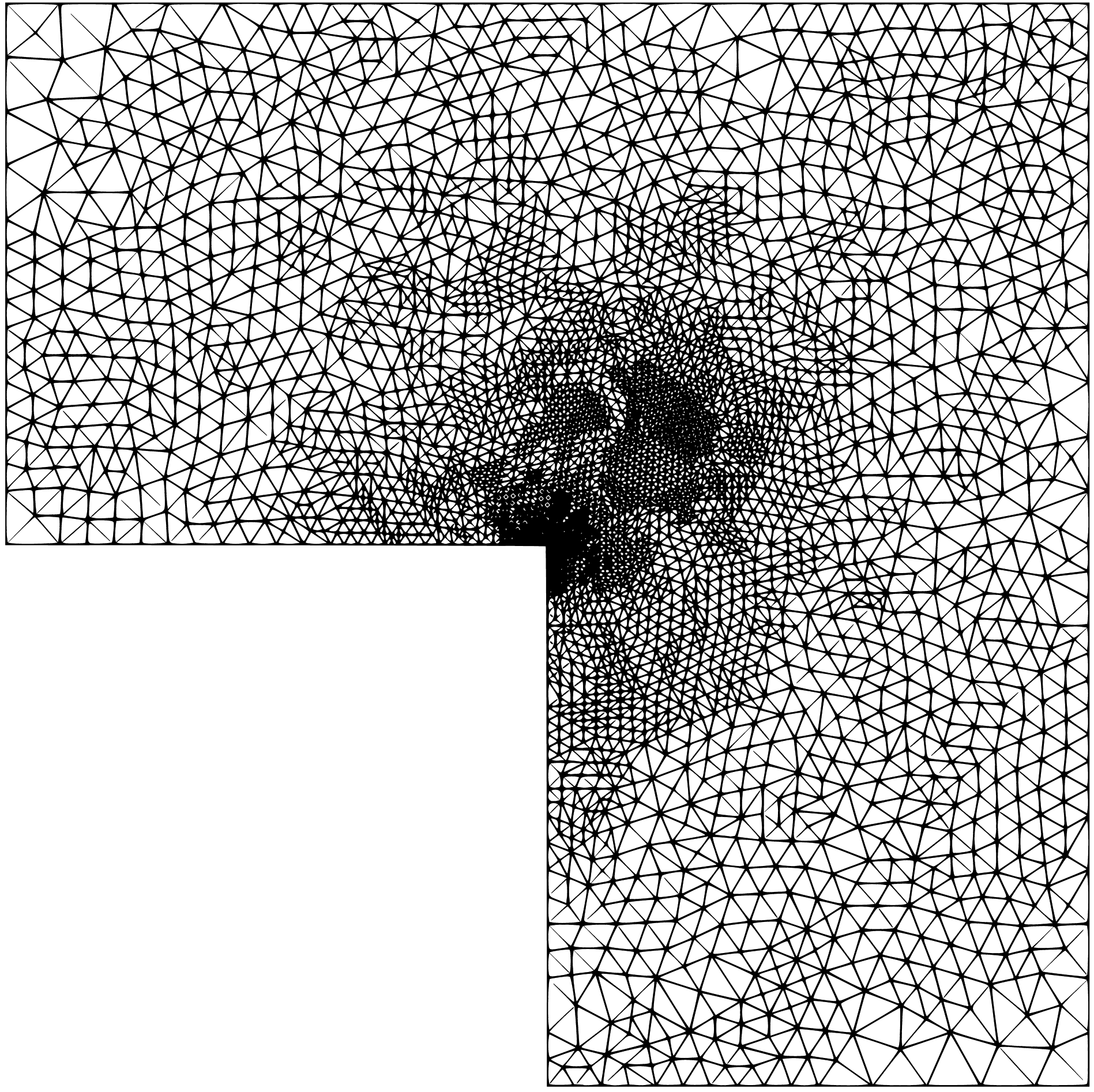}
            \\
            & \includegraphics[width=0.3\textwidth]{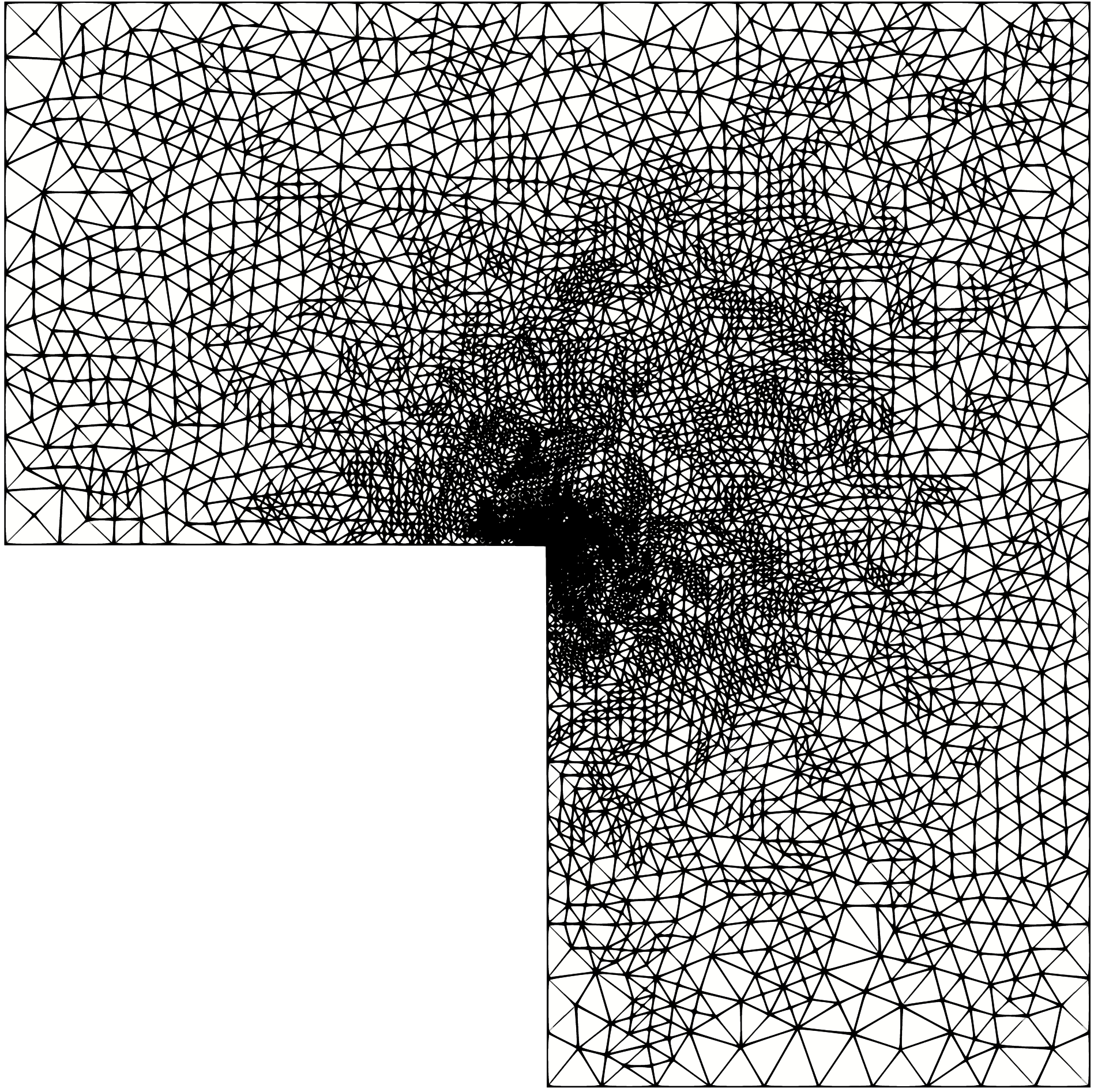}
            &
            \includegraphics[width=0.3\textwidth]{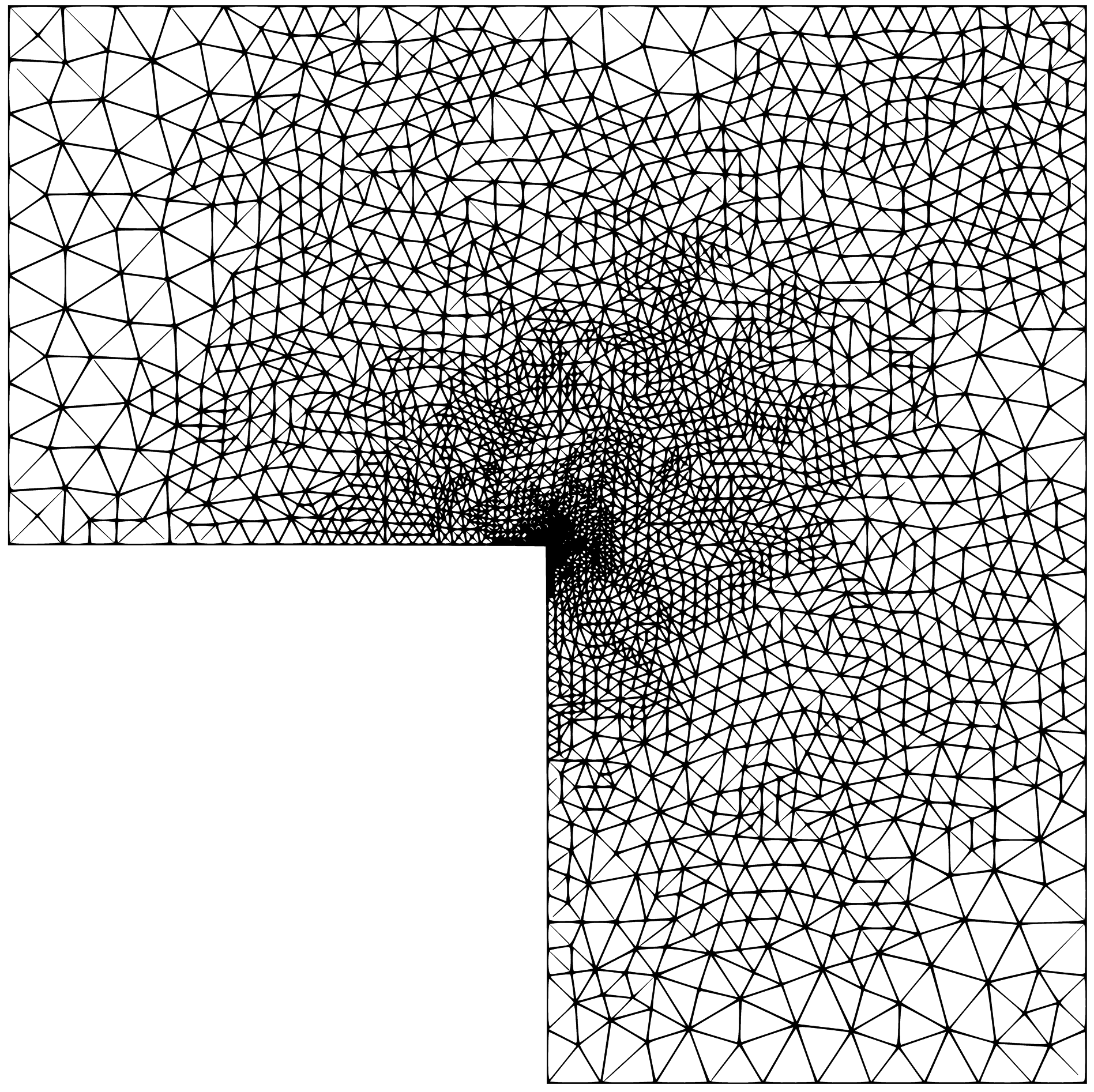}
        \end{tabular}
        \caption{L-shaped goal-oriented Poisson problem with linear elements: On top left the initial mesh used to start all the adaptive
            strategies. From top middle to bottom right, the adaptive meshes
        obtained after seven iterations of refinement strategies steered by
        weighted estimators derived respectively from
        $\eta_{\bw}^{4,2}$, $\eta_{\mathrm{res}}$, $\eta_{\bw}^{\mathrm{b}}$ and
        $\eta_{\mathrm{zz}}$ for both primal and dual problems.}
        \label{fig:goal-oriented-meshes}
    \end{figure}

\begin{figure}
        \includegraphics[trim=0 6 0 0, clip, width=\textwidth]{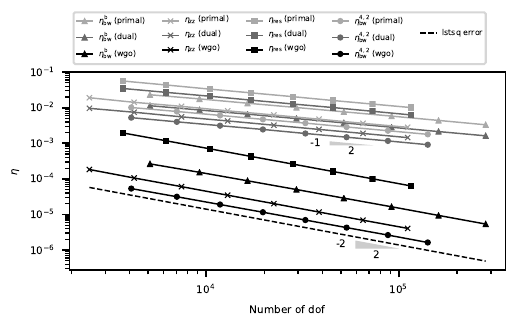}
        \caption{L-shaped goal-oriented Poisson problem with linear elements: plot comparing
            convergence of some goal-oriented adaptive strategies driven by four
            different estimators.
            Expected rates for primal and dual problems ($-0.5$) and goal
            functional ($-1$) shown by triangle markers.
            Comparison with an indicative line representing the higher order approximation
            of the errors of each strategy and obtained using least squares
            method.} 
        \label{fig:graph_boundary_singularity}
    \end{figure}

    \cref{fig:goal-oriented-meshes} shows refined meshes after seven
    iterations of the weighted goal oriented method.
    We can see that the meshes are mainly refined in the re-entrant corner as
    well as in the region on the right top of it where the goal functional
    focuses.
    In \cref{fig:graph_boundary_singularity} we show the convergence curves of
    some of these adaptive strategies. For each strategy, $\eta_u = \eta_z$ is
    the estimator specified in the legend and $\eta_w = \eta_u \eta_z$.
    All the strategies we have tried led to very similar higher-order
    approximate errors. So for the sake of clarity we have replaced the
    approximate errors by an indicative line computed using a regression from
    the least squares method (lstsq error), leading to the line that fits the best
    the values of the different approximate errors.
    As we can see, these adaptive strategies are reaching an optimal convergence
    rate.
    Although it is also the case for all the other strategies we have tried, we
    do not show the other results for the sake of concision.
    In the left table of \cref{fig:table_bws_goal_oriented} we show the
    efficiencies of the estimators $\eta_{w}$ where $\eta_u = \eta_z =
    \eta_{\bw}^{k_+,k_-}$.
    On the right table of \cref{fig:table_bws_goal_oriented} we take $\eta_u =
    \eta_z$ to be the estimators in the left column.
    As we can see on the efficiencies are not as good as in \cref{sec:lshaped}.
    The two best estimators are those derived from $\eta_{\bw}^{4,0}$ and
    $\eta_{\bw}^{4,2}$.
    With $\eta_{\bw}^{3,0}$, they are the only cases where the goal-oriented
    estimator $\eta_w$ is performing better than the goal-oriented estimator
    derived from the Zienkiewicz--Zhu estimator.
    The estimators $\eta_w$ derived from the bubble Bank--Weiser estimator as
    well as from the residual estimator are poorly overestimating the error.

    \begin{figure}
    \begin{center}
    \begin{tikzpicture}[x=0.8cm,y=0.6cm] 
    \def\ab{6.22};
    \def\ac{7.24};
    \def\bc{4.94};
    \def\ad{8.06};
    \def\bd{7.75};
    \def\cd{0.0};
    \def\ae{8.8};
    \def\be{10.02};
    \def\ce{1.66};
    \def\de{0.25};
    \def\oa{50.22};
    \def\ob{11.14};
    \def\oc{3.04};
    \draw (0,-1) -- (5, -1) -- (5, 4) -- (0, 4) -- (0, -1);
    \draw (0, 1) -- (5, 1);
    \draw (0, 2) -- (5, 2);
    \draw (0, 3) -- (5, 3);
    \draw (1, -1) -- (1, 4);
    \draw (2, -1) -- (2, 4);
    \draw (3, -1) -- (3, 4);
    \draw (4, -1) -- (4, 4);
    \draw (0, 4) -- (1, 3);
    \draw (0, 0) -- (5, 0);
    \draw (6.01, -1) -- (8.01, -1) -- (8.01, 2) -- (6.01, 2) -- (6.01, -1);
    \draw (6.01, -1) -- (8.01, -1);
    \draw (6.01, 0) -- (8.01, 0);
    \draw (6.01, 1) -- (8.01, 1);
    \draw (7.01, -1) -- (7.01, 2);

    \node at (0.33, 3.25) [text width = 3cm, align=center] {$k_-$};
    \node at (0.75, 3.66) [text width = 3cm, align=center] {$k_+$};
    \node at (0.5, -0.5) [text width = 3cm, align=center] {$3$};
    \node at (0.5, 0.5) [text width = 3cm, align=center] {$2$};
    \node at (0.5, 1.5) [text width = 3cm, align=center] {$1$};
    \node at (0.5, 2.5) [text width = 3cm, align=center] {$0$};

    \node at (1.5, 3.5) [text width = 3cm, align=center] {$1$};
    \node at (2.5, 3.5) [text width = 3cm, align=center] {$2$};
    \node at (3.5, 3.5) [text width = 3cm, align=center] {$3$};
    \node at (4.5, 3.5) [text width = 3cm, align=center] {$4$};

    \node at (1.5, 2.5) [align=center] {\ab};
    \node at (2.5, 2.5) [align=center] {\ac};
    \node at (3.5, 2.5) [align=center] {\ad};
    \node at (4.5, 2.5) [align=center] {\ae};

    \node at (1.5, 1.5) [align=center] {$\varnothing$};
    \node at (2.5, 1.5) [align=center] {\bc};
    \node at (3.5, 1.5) [align=center] {\bd};
    \node at (4.5, 1.5) [align=center] {\be};

    \node at (1.5, 0.5) [align=center] {$\varnothing$};
    \node at (2.5, 0.5) [align=center] {$\varnothing$};
    \node at (3.5, 0.5) [align=center] {\cd};
    \node at (4.5, 0.5) [align=center] {\ce};

    \node at (1.5, -0.5) [align=center] {$\varnothing$};
    \node at (2.5, -0.5) [align=center] {$\varnothing$};
    \node at (3.5, -0.5) [align=center] {$\varnothing$};
    \node at (4.5, -0.5) [align=center] {\de};

    \node at (6.51, 1.5) [text width = 3cm, align=center] {$\eta_{\mathrm{res}}$};
    \node at (6.51, 0.5) [text width = 3cm, align=center] {$\eta^{\mathrm{b}}_{\bw}$};
    \node at (6.51, -0.5) [text width = 3cm, align=center] {$\eta_{\mathrm{zz}}$};
    \node at (7.51, 1.5) [text width = 3cm, align=center] {\oa};
    \node at (7.51, 0.5) [text width = 3cm, align=center] {\ob};
    \node at (7.51, -0.5) [text width = 3cm, align=center] {\oc};
    \end{tikzpicture}
    \end{center}
    \caption{L-shaped goal-oriented Poisson problem with linear elements: efficiencies of the dual weighted estimators derived
    from $\eta_{\bw}^{k_+,k_-}$ and other estimators on the last mesh of an adaptively refined hierarchy.}
    \label{fig:table_bws_goal_oriented}
    \end{figure}

    \subsection{Nearly-incompressible elasticity}
    We consider the linear elasticity problem from\ \cite{carstensen_robust_2016} on
    the centered unit square domain $\Omega$ with homogeneous Dirichlet boundary
    conditions on $\Gamma_D = \Gamma$ ($u_D = 0$).
    The first Lamé coefficient is set to $\mu = 100$ and the Poisson ratio to
    $\nu = 0.3$ and $\nu = 0.499$.
    The problem data $\bm{f}$ is given by $\bm f = (f_1, f_2)$ with
    \begin{equation}
        \begin{array}{r c l}
            f_1(x, y) = -2\mu \pi^3 \cos(\pi y) \sin(\pi y)\big(2\cos(2\pi x) - 1\big),\\
            f_2(x, y) = 2 \mu \pi^3 \cos(\pi x) \sin(\pi x)\big(2\cos(2\pi y) - 1 \big).
        \end{array}
    \end{equation}
    The corresponding exact solution of the linear elasticity problem reads
    $\bm u =(u_1, u_2)$ with
    \begin{equation}
            u_1(x,y) = \pi \cos(\pi y)\sin^2(\pi x)\sin(\pi y),\
            u_2(x, y) = -\pi \cos(\pi x) \sin(\pi x)\sin^2(\pi y),
    \end{equation}
    the Herrmann pressure is zero everywhere on $\Omega$.
    In each case we discretize this problem using the Taylor--Hood element and an
    initial cartesian mesh and we apply our adaptive procedure driven by the
    Poisson estimator described in \cref{sec:elasticity}.
    We compare the Poisson estimators derived from different Bank--Weiser
    estimators and the residual estimator.

    As before, all the refinement strategies are achieving an optimal
    convergence rate no matter the value of $\nu$.
    \cref{fig:table_bws_elasticity} shows the results for $\nu=0.3$.
    We notice that almost all the Poisson estimators derived from Bank--Weiser
    estimators have a very good efficiency.
    The best estimator in this case is $\eta_{\bw}^{2,0}$ closely followed by
    $\eta_{\bw}^{3,0}$, $\eta_{\bw}^{\mathrm{b}}$ and $\eta_{\bw}^{4,0}$.
    Although the residual estimator still performs the worst, it is sharper
    than in all the previous test cases.
    As we can notice on \cref{fig:table_bws_elasticity_incompressible}, all the
    estimators are robust with respect to the incompressibility constraint.
    All the efficiencies have slightly increased and some estimators
    ($\eta_{\bw}^{2,0}$ and $\eta_{\bw}^{3,0}$) that where a lower bound of the error
    previously are now an upper bound.

    \begin{figure}
    \begin{center}
    \begin{tikzpicture}[x=0.8cm,y=0.6cm] 
    \def\ab{0.87};
    \def\ac{0.98};
    \def\bc{0.68};
    \def\ad{0.97};
    \def\bd{1.09};
    \def\cd{0.57};
    \def\ae{1.02};
    \def\be{1.19};
    \def\ce{0.68};
    \def\de{0.36};
    \def\oa{2.44};
    \def\ob{1.1};
    \draw (0,-1) -- (5, -1) -- (5, 4) -- (0, 4) -- (0, -1);
    \draw (0, 1) -- (5, 1);
    \draw (0, 2) -- (5, 2);
    \draw (0, 3) -- (5, 3);
    \draw (1, -1) -- (1, 4);
    \draw (2, -1) -- (2, 4);
    \draw (3, -1) -- (3, 4);
    \draw (4, -1) -- (4, 4);
    \draw (0, 4) -- (1, 3);
    \draw (0, 0) -- (5, 0);
    \draw (6.01, -1) -- (8.01, -1) -- (8.01, 2) -- (6.01, 2) -- (6.01, -1);
    \draw (6.01, -1) -- (8.01, -1);
    \draw (6.01, 0) -- (8.01, 0);
    \draw (6.01, 1) -- (8.01, 1);
    \draw (7.01, -1) -- (7.01, 2);

    \node at (0.33, 3.25) [text width = 3cm, align=center] {$k_-$};
    \node at (0.75, 3.66) [text width = 3cm, align=center] {$k_+$};
    \node at (0.5, -0.5) [text width = 3cm, align=center] {$3$};
    \node at (0.5, 0.5) [text width = 3cm, align=center] {$2$};
    \node at (0.5, 1.5) [text width = 3cm, align=center] {$1$};
    \node at (0.5, 2.5) [text width = 3cm, align=center] {$0$};

    \node at (1.5, 3.5) [text width = 3cm, align=center] {$1$};
    \node at (2.5, 3.5) [text width = 3cm, align=center] {$2$};
    \node at (3.5, 3.5) [text width = 3cm, align=center] {$3$};
    \node at (4.5, 3.5) [text width = 3cm, align=center] {$4$};

    \node at (1.5, 2.5) [align=center] {\ab};
    \node at (2.5, 2.5) [align=center] {\ac};
    \node at (3.5, 2.5) [align=center] {\ad};
    \node at (4.5, 2.5) [align=center] {\ae};

    \node at (1.5, 1.5) [align=center] {$\varnothing$};
    \node at (2.5, 1.5) [align=center] {\bc};
    \node at (3.5, 1.5) [align=center] {\bd};
    \node at (4.5, 1.5) [align=center] {\be};

    \node at (1.5, 0.5) [align=center] {$\varnothing$};
    \node at (2.5, 0.5) [align=center] {$\varnothing$};
    \node at (3.5, 0.5) [align=center] {\cd};
    \node at (4.5, 0.5) [align=center] {\ce};

    \node at (1.5, -0.5) [align=center] {$\varnothing$};
    \node at (2.5, -0.5) [align=center] {$\varnothing$};
    \node at (3.5, -0.5) [align=center] {$\varnothing$};
    \node at (4.5, -0.5) [align=center] {\de};

    \node at (6.51, 1.5) [text width = 3cm, align=center] {$\eta_{\mathrm{res}}$};
    \node at (6.51, 0.5) [text width = 3cm, align=center] {$\eta^{\mathrm{b}}_{\bw}$};
    \node at (6.51, -0.5) [text width = 3cm, align=center] {$\eta_{\mathrm{zz}}$};
    \node at (7.51, 1.5) [text width = 3cm, align=center] {\oa};
    \node at (7.51, 0.5) [text width = 3cm, align=center] {\ob};
    \node at (7.51, -0.5) [text width = 3cm, align=center] {$\varnothing$};
    \end{tikzpicture}
    \end{center}
    \caption{Nearly-incompressible elasticity ($\nu=0.3$) problem with
        Taylor--Hood elements: efficiencies of the Poisson
    estimators derived from $\eta_{\bw}^{k_+,k_-}$ and other estimators on the
    last mesh of an adaptively refined hierarchy.}
    \label{fig:table_bws_elasticity}
    \end{figure}
    
    \begin{figure}
    \begin{center}
    \begin{tikzpicture}[x=0.8cm,y=0.6cm] 
    \def\ab{0.94};
    \def\ac{1.02};
    \def\bc{0.77};
    \def\ad{1.04};
    \def\bd{1.1};
    \def\cd{0.67};
    \def\ae{1.1};
    \def\be{1.22};
    \def\ce{0.77};
    \def\de{0.44};
    \def\oa{2.47};
    \def\ob{1.13};
    \draw (0,-1) -- (5, -1) -- (5, 4) -- (0, 4) -- (0, -1);
    \draw (0, 1) -- (5, 1);
    \draw (0, 2) -- (5, 2);
    \draw (0, 3) -- (5, 3);
    \draw (1, -1) -- (1, 4);
    \draw (2, -1) -- (2, 4);
    \draw (3, -1) -- (3, 4);
    \draw (4, -1) -- (4, 4);
    \draw (0, 4) -- (1, 3);
    \draw (0, 0) -- (5, 0);
    \draw (6.01, -1) -- (8.01, -1) -- (8.01, 2) -- (6.01, 2) -- (6.01, -1);
    \draw (6.01, -1) -- (8.01, -1);
    \draw (6.01, 0) -- (8.01, 0);
    \draw (6.01, 1) -- (8.01, 1);
    \draw (7.01, -1) -- (7.01, 2);

    \node at (0.33, 3.25) [text width = 3cm, align=center] {$k_-$};
    \node at (0.75, 3.66) [text width = 3cm, align=center] {$k_+$};
    \node at (0.5, -0.5) [text width = 3cm, align=center] {$3$};
    \node at (0.5, 0.5) [text width = 3cm, align=center] {$2$};
    \node at (0.5, 1.5) [text width = 3cm, align=center] {$1$};
    \node at (0.5, 2.5) [text width = 3cm, align=center] {$0$};

    \node at (1.5, 3.5) [text width = 3cm, align=center] {$1$};
    \node at (2.5, 3.5) [text width = 3cm, align=center] {$2$};
    \node at (3.5, 3.5) [text width = 3cm, align=center] {$3$};
    \node at (4.5, 3.5) [text width = 3cm, align=center] {$4$};

    \node at (1.5, 2.5) [align=center] {\ab};
    \node at (2.5, 2.5) [align=center] {\ac};
    \node at (3.5, 2.5) [align=center] {\ad};
    \node at (4.5, 2.5) [align=center] {\ae};

    \node at (1.5, 1.5) [align=center] {$\varnothing$};
    \node at (2.5, 1.5) [align=center] {\bc};
    \node at (3.5, 1.5) [align=center] {\bd};
    \node at (4.5, 1.5) [align=center] {\be};

    \node at (1.5, 0.5) [align=center] {$\varnothing$};
    \node at (2.5, 0.5) [align=center] {$\varnothing$};
    \node at (3.5, 0.5) [align=center] {\cd};
    \node at (4.5, 0.5) [align=center] {\ce};

    \node at (1.5, -0.5) [align=center] {$\varnothing$};
    \node at (2.5, -0.5) [align=center] {$\varnothing$};
    \node at (3.5, -0.5) [align=center] {$\varnothing$};
    \node at (4.5, -0.5) [align=center] {\de};

    \node at (6.51, 1.5) [text width = 3cm, align=center] {$\eta_{\mathrm{res}}$};
    \node at (6.51, 0.5) [text width = 3cm, align=center] {$\eta^{\mathrm{b}}_{\bw}$};
    \node at (6.51, -0.5) [text width = 3cm, align=center] {$\eta_{\mathrm{zz}}$};
    \node at (7.51, 1.5) [text width = 3cm, align=center] {\oa};
    \node at (7.51, 0.5) [text width = 3cm, align=center] {\ob};
    \node at (7.51, -0.5) [text width = 3cm, align=center] {$\varnothing$};
    \end{tikzpicture}
    \end{center}
    \caption{Nearly-incompressible elasticity ($\nu=0.499$) problem with
        Taylor--Hood elements: efficiencies of the Poisson
    estimators derived from $\eta_{\bw}^{k_+,k_-}$ and other estimators on the
    last mesh of an adaptively refined hierarchy.}
    \label{fig:table_bws_elasticity_incompressible}
    \end{figure}

    \subsection{Human femur modelled using linear elasticity}
    In this test case we consider a linear elasticity problem on a domain inspired
    by a human femur bone~\footnote{The STL model of the femur bone can be
    found at \url{https://3dprint.nih.gov/discover/3dpx-000168} under a Public
    Domain license.}.

    The goal of this test case is not to provide an accurate description of the
    behavior of the femur bone but to demonstrate the applicability of our
    implementation to 3D dimensional goal-oriented problem with large number of
    degrees of freedom: the linear elasticity problem to solve on the initial
    mesh, using Taylor-Hood element has 247,233 degrees of freedom while our
    last refinement step reaches 3,103,594 degrees of freedom.

\begin{figure}
        \begin{center}
            \includegraphics[trim= 0 500 0 500, clip, width=\textwidth]{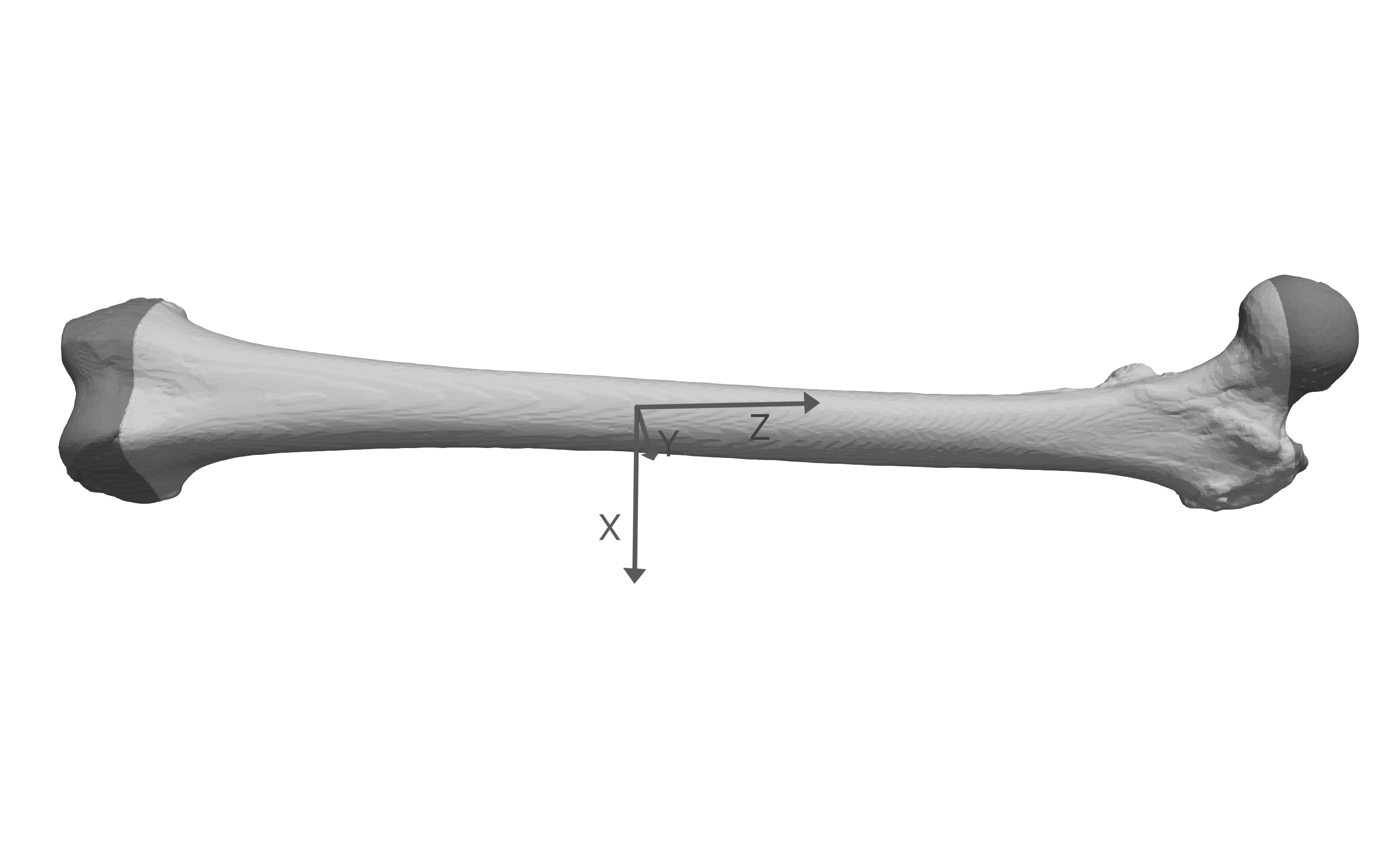}\\
            \vspace{-5pt}
            \includegraphics[trim= 0 600 0 600, clip, width=\textwidth]{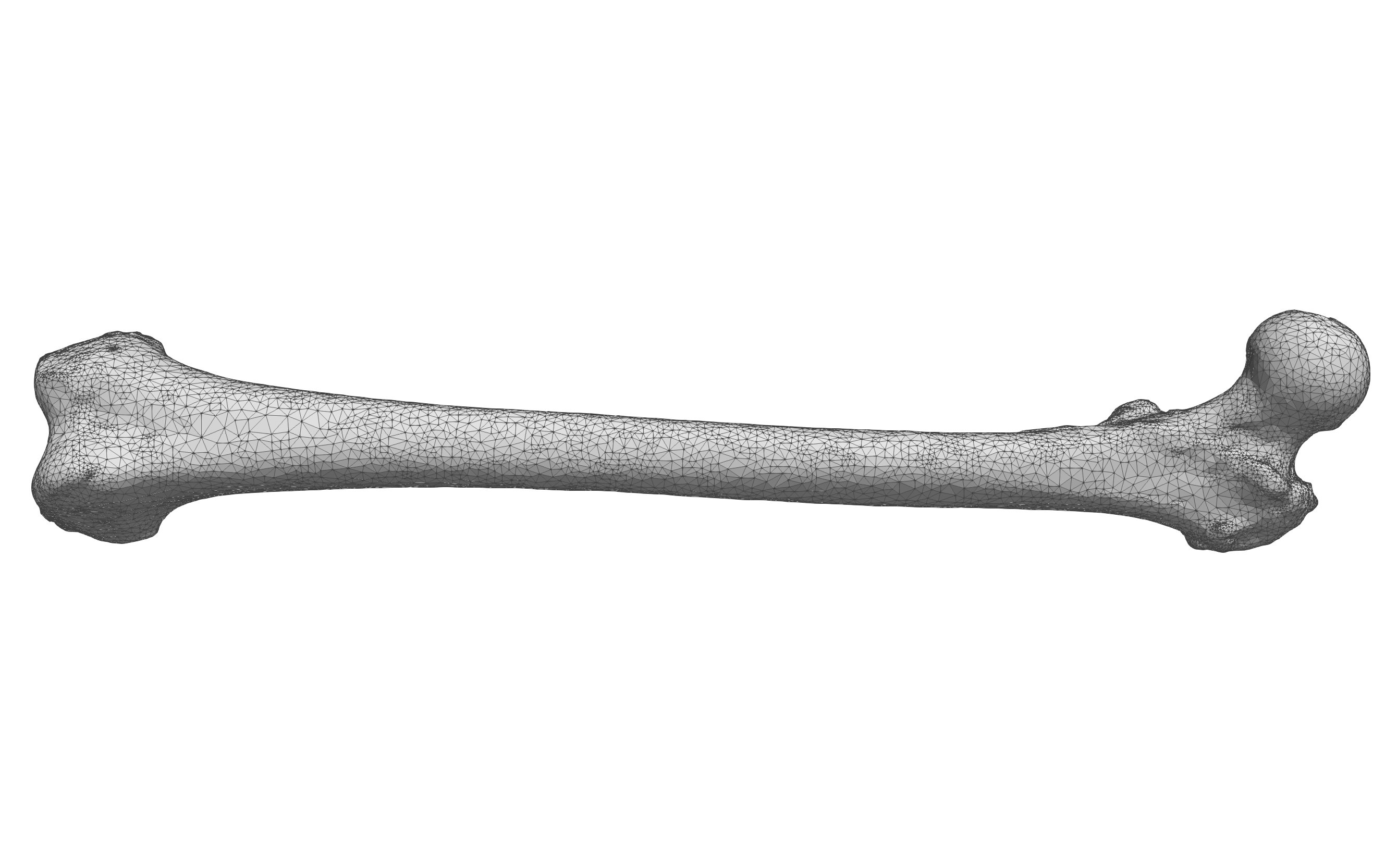}\\
            \includegraphics[trim= 0 600 0 600, clip, width=\textwidth]{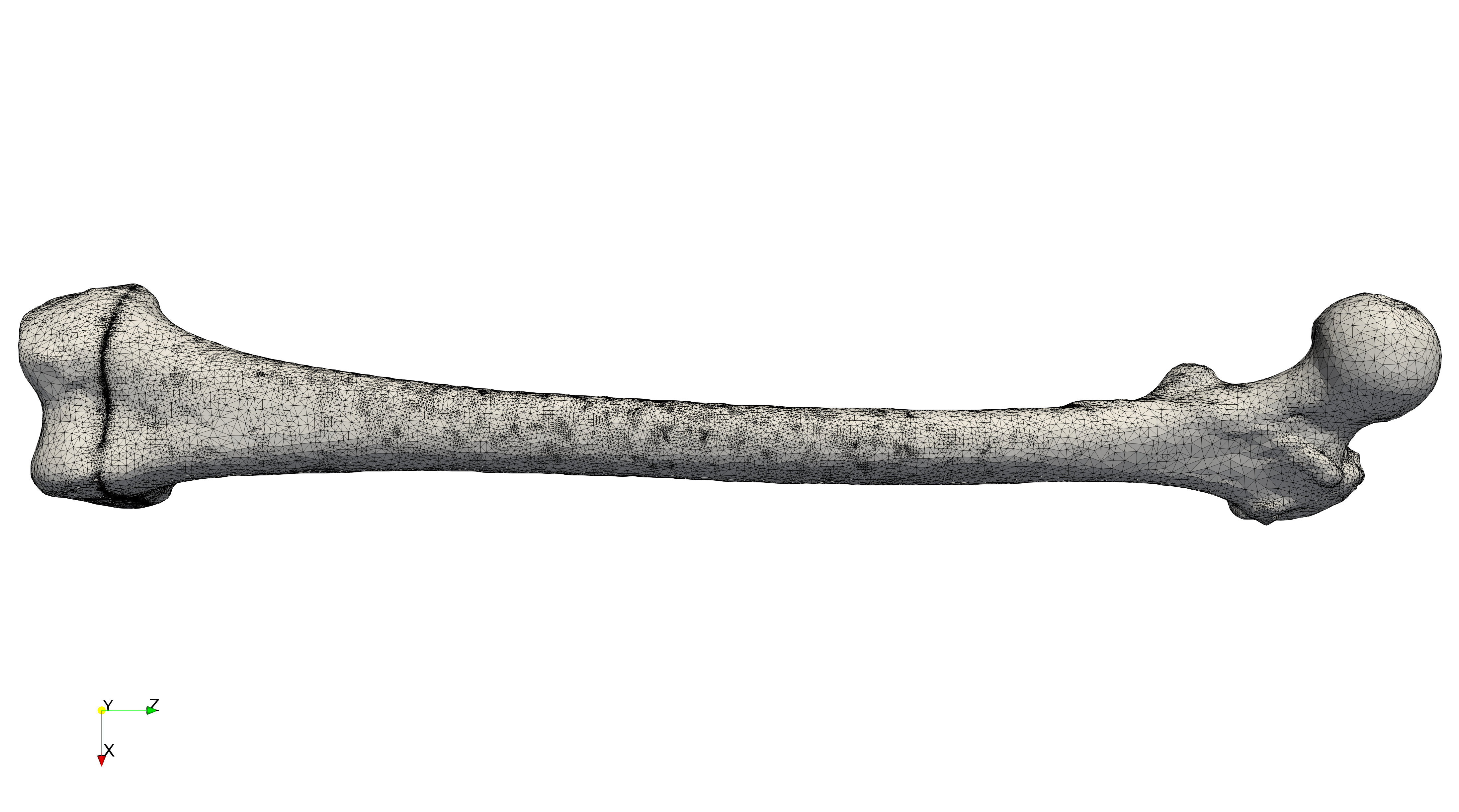}
        \end{center}
        \caption{Femur bone linear elasticity problem with Taylor--Hood
            elements: on the top, the three different regions of
        the boundary corresponding to different boundary conditions: the left dark
        grey region is the non-zero Neumann boundary, the middle
        light grey region is the zero Neumann boundary and the right dark grey
        region is the Dirichlet boundary.
        In the middle, the initial mesh.
        On the bottom, the last mesh after several steps of adaptive refinement.}
        \label{fig:femur_domain_bcs}
    \end{figure}

    The 3D mesh for analysis is build from the surface model using the C++ library CGAL
    \cite{cgal} via the Python front-end pygalmesh. 
    The material parameters, namely the Young's modulus is set to $20$ GPa and the Poisson's ratio
    to $0.42$ (see e.g.\ \cite{rupin_experimental_2008}).
    In addition, the load is given by $\mathbf f = (0, 0, 0)$, the Dirichlet
    data by $\bm u_D = 0$ on $\Gamma_D \subsetneq \Gamma$ represented as the
    left dark gray region of the boundary in \cref{fig:femur_domain_bcs} and $\bm g$
    the traction data is defined as $\bm g = (0, 0, 0)$ on the center light gray
    region of the boundary and is constant on the right dark gray region of the
    boundary $\bm g = (-10^{-7}, -10^{-7}, 10^{-6})$.
    The femur-shaped domain $\Omega$ as well as the initial and last meshes are
    shown in \cref{fig:femur_domain_bcs}.
    As we can see, the refinement occurs mainly in the central region of the
    femur, where the goal functional $J$ focus.
    Some artefacts can be seen as stains of refinement in the central region due
    to the fact that we use the initial mesh as our geometry and on the left due
    to the discontinuity in the boundary conditions.

\begin{figure}
        \includegraphics[trim= 0 5 0 0, clip, width=\textwidth]{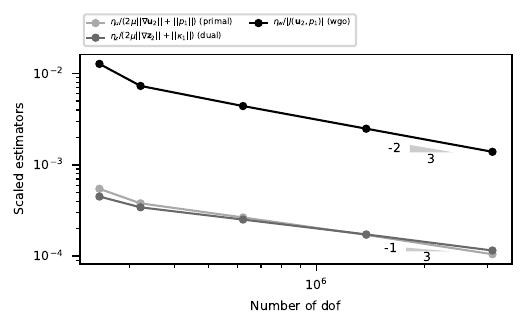}
        \caption{Femur bone linear elasticity problem with Taylor--Hood
            elements: convergence curves of the primal, dual and weighted estimators
            respectively scaled by the norm of the primal solution, dual
            solution and magnitude of the goal functional evaluated in the
        primal solution.}
        \label{fig:femur_conv}
    \end{figure}

    In \cref{fig:femur_conv} the primal solution is given by the couple
    $(\mathbf{u}_2, p_1)$ and the dual solution by $(\mathbf{z}_2, \kappa_1)$.
    As we can notice and as expected, the weighted estimator $\eta_w$ converges
    twice as fast as the primal and dual estimators.


    \subsection{Strong scaling study}\label{sec:strongscaling}
    Finally, we provide results showing that our implementation scales strongly
    in parallel and that for a large-scale three-dimensional problem this error
    estimation takes significantly less time than the solution of the primal
    problem.
    In this section we use the new DOLFINx solver\ \cite{Habera2020} with the
    matching implementation of our algorithm.
    
    We briefly discuss some aspects that are important for interpretation of the
    results. For a given cell the computation of the Bank--Weiser estimator
    requires geometry and solution data on the current cell and on all cells
    attached across its facets. So in a parallel computing context, cells
    located on the boundary of a partition require data from cells owned by
    another process. Both DOLFIN and DOLFINx support
    facet--mode ghosting where all data owned by cells on a partition boundary
    that share a facet are duplicated by the other process (ghost data). After
    the solution of the primal linear system the ghost data is updated between
    processes, which requires parallel communication. After this update, each
    process has a local copy of all of the data from the other rank needed to
    compute the Bank--Weiser estimator, and so the computation of the estimator
    is entirely local to a rank, i.e.\ without further parallel communication.
    
    Because of this locality a proper implementation of this algorithm
    should demonstrate strong scaling performance. Furthermore, it would
    be desirable that the error estimation takes significantly less time than
    the solution of the primal problem even when using state-of-the art linear
    solution strategies. The results in this section demonstrate that this is
    indeed the case.
    
    We solve \cref{eq:strong_form}
    where $\Omega$ is the unit cube $[0,1]^3$, $\Gamma_D = \partial \Omega$ and $\Gamma_N =
    \varnothing$.
    The data of this problem are given by $f(x,y,z) = 12 \pi^2 \sin(2\pi x)
    \sin(2\pi y) \sin(2\pi z)$ and $u_D(x, y, z) = 0$.
    Given these data the solution $u$ of \cref{eq:strong_form} is given by $u(x,
    y, z) = \sin(2\pi x) \sin(2\pi y) \sin(2\pi z)$.
    We use continuous quadratic Lagrange finite elements and the Bank--Weiser error
    estimation is performed using the pair $V_T^3$/$V_T^2$. The primal linear
    system matrix and right-hand side vector are assembled using standard
    routines in DOLFINx. The resulting linear system is solved with
    PETSc~\cite{balay_petsc_2016} using the conjugate gradient method preconditioned
    with Hypre BoomerAMG algebraic multigrid~\cite{falgout_hypre_2002}.
    
    The strong scaling study was carried out on the Aion cluster within the HPC
    facilities of the University of Luxembourg~\cite{varrette_management_2014}.
    The Aion cluster is a Atos/Bull/AMD supercomputer composed of 318 compute
    nodes each containing two AMD Epyc ROME 7H12 processors with 64 cores per
    processor (128 cores per node). The nodes are connected through a Fast
    InfiniBand (IB) HDR 100Gbps interconnect in a `fat-tree' topology. We invoke jobs using SLURM and ask for a
    contiguous allocation of nodes and exclusivity (no competing jobs) on each
    node. DOLFINx and PETSc are built using GCC 10.2.0 with Intel MPI and
    OpenBLAS. We use DOLFINx through its Python interface. The problem size is
    kept fixed at around 135 million degrees of freedom and the number of MPI
    ranks is increased from 128 (1 node, no interconnect communication) through to
    2048 (16 nodes, interconnect communication) by doubling the number of nodes
    and ranks used in the previous computation.
    
\begin{figure}
        \includegraphics[trim= 0 5 0 0, clip,width=\textwidth]{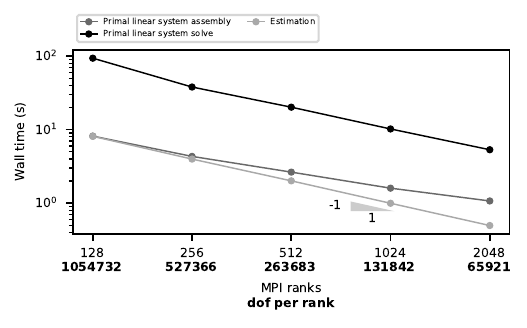}
        \caption{Strong scaling study on the University of Luxmebourg Aion HPC. Wall times for primal linear system assembly,
        primal linear system solve and error estimation of a three-dimensional
        Poisson finite element problem on the unit square, discretized with quadratic
        elements. $1/1$ triangle represents ideal strong scaling.}
        \label{fig:strong_scaling}
    \end{figure}

    In \cref{fig:strong_scaling} we show the results of the strong scaling
    study. We show wall time against MPI ranks and dof per rank for the primal
    linear system assembly, primal linear system solve, and the error
    estimation. For error estimation we are measuring steps 2 through 5 of
    \cref{subsec:outline}. Both the solve and estimation scale almost perfectly
    down to around 65 thousand dofs per rank. The primal system assembly does
    not scale as well as the estimation.
    This is because the primal system assembly is constrained by communication
    overheads and memory bandwidth, whereas the Bank--Weiser estimator
    computation is fully local and has much higher arithmetic intensity, so has
    not yet hit bandwidth limits of our system on the largest run.
    A further study (results not shown) using 96 MPI
    ranks per node yielded lower wall times and better strong scaling for primal
    linear system assembly, but the overall time for estimation and linear
    system solve increased and dominated any gains made in assembly. Comparing
    linear system assembly and solve with estimation time we can see that
    estimation is approximately one order of magnitude faster than solve time.
    
    \section{Conclusions}\label{sec:conclusions}
    In this paper we have shown how the error estimator of Bank--Weiser,
    involving the solution of a local problem on a special finite element
    space, can be mathematically reformulated and implemented
    straightforwardly in a modern finite element software with the aid of
    automatic code generation techniques.
    Through a series of numerical results we have shown that the estimator
    is highly competitive in accurately predicting the total global error and
    in driving an adaptive mesh refinement strategy.
    Furthermore, the basic methodology and implementation for the Poisson
    problem can be extended to tackle more complex mixed discretizations of
    PDEs including nearly-incompressible elasticity or Stokes problems.
    We have also shown the (strong) scalability of our method when implemented
    in parallel and that the error estimation time is significantly lower than the primal solution time on a large problem.

    \subsection*{Acknowledgements}
    The three-dimensional problems presented in this paper were carried out
    using the HPC facilities of the University of Luxembourg\
    \cite{varrette_management_2014}.

    We would like to thank Nathan Sime and Chris Richardson for helpful
    discussions on FEniCS's discontinuous Galerkin discretizations and PLAZA
    refinement strategy, respectively. We would also like to thank Prof.\
    Roland Becker for the motivating and fruitful discussions.

    \subsection*{Supplementary material}
    The DOLFINx version of the code can be found at \url{https://github.com/jhale/fenicsx-error-estimation} and the DOLFIN version
    at \url{https://github.com/rbulle/fenics-error-estimation}.
    
    The appendices 1 and 2 in of supplementary material contain two snippets
    showing the implementation of the Poisson estimator and the Poisson
    estimator for the nearly-incompressible elasticity problem using the DOLFIN
    version of the code.
    A simplified version of the code (LGPLv3) used to produce the results in
    this paper is archived at \url{https://doi.org/10.6084/m9.figshare.10732421}.
    A Docker image\ \cite{hale_containers_2017} is provided in which this code
    can be executed.

    \section*{Appendices}
    \renewcommand{\thesection}{A}
    \setcounter{equation}{0}

    \section{The residual estimator}\label{app:res}
    \subsection{Poisson equation}
    The class of residual estimators, the explicit residual estimator is part
    of, have been introduced for the first time in\
    \cite{babuska_posteriori_1978}.
    Let $h_T$ be the diameter (see e.g.\ \cite{scott_finite_1990}) of the cell
    $T$ and $h_E$ be the diameter of the facet $E$.
    The explicit residual estimator~\cite{ainsworth_posteriori_2011} on a cell
    $T$ for the Poisson problems \cref{eq:weak_form,eq:discrete_weak_form} is
    defined as
    \begin{equation}
        \label{eq:local_residual_estimator}
        \eta_{\mathrm{res}, T}^2 := h_T^2 \norm{f_T+\Delta u_k}_{T}^2
        + \sum_{E \in \E_I \cap \partial T} \frac{1}{2} h_E \norm{\jump{\partial_n u_k}_E}^2_{E}
        + \sum_{E \in \E_N \cap \partial T} h_E \norm{g_E-\partial_n u_k}^2_E,
    \end{equation}
    where $f_h$ and $g_h$ are the $L^2$ projections of $f$ and $g$ on $V^k$
    respectively.
    In order to take into account inhomogeneous Dirichlet boundary conditions,
    we define in addition the Dirichlet oscillations.
    If $E := \Gamma_D \cap T \neq \varnothing$, then
    \begin{equation}
        \mathrm{osc}_{D, E}^2 := h_E \norm{\nabla_{\Gamma} \left(g_E -
        u_k\right)}^2_{L^2(E)},
        \label{eq:dirichlet_oscillations}
    \end{equation}
    where $\nabla_{\Gamma}$ is the surface gradient and $g_E := \pi^+_T(g)$ is
    the $L^2$ projection of $g$ onto $V^{k+1}_T$ \cite{Aurada2013}.
    The global residual estimator reads
    \begin{equation}
        \label{eq:global_residual_estimator}
        \eta_{\mathrm{res}}^2 := \sum_{T \in \T} \eta_{\mathrm{res}, T}^2 +
        \mathrm{osc}_{D, \overline T \cap \Gamma_D}^2.
    \end{equation}
    \subsection{Linear elasticity equations}
    The residual estimator for the linear elasticity problem
    \cref{eq:elasticity:weak:equilibrium,eq:elasticity:weak:pressure,eq:elasticity:weak:discrete:equilibrium,eq:elasticity:weak:discrete:pressure}
    is given by
    \begin{equation}
        \eta_{\mathrm{res},T}^2 := \rho_T \norm{\bm R_T}_T^2 + \rho_d
        \norm{r_T}_T^2 + \sum_{E \in \partial T} \rho_E \norm{\bm R_E}_E^2,
        \label{eq:elasticity_local_residual}
    \end{equation}
    where the residuals $\bm R_T$, $r_T$ and $\bm R_E$ are respectively defined
    in
    \cref{eq:elasticity_vectorial_interior_res,eq:elasticity_scalar_interior_res,eq:elasticity_edge_res}
    and the constants $\rho_T$, $\rho_d$ and $\rho_E$ are given by
    \begin{equation}
        \rho_T := \frac{h_T (2\mu)^{-1/2}}{2},\quad \rho_d := \big(\lambda^{-1}
        + (2\mu)^{-1}\big)^{-1},\quad \rho_E := \frac{h_E(2\mu)^{-1}}{2},
    \end{equation}
    with $h_T$ the diameter of the cell $T$ and $h_E$ the length of the edge $E$.
    The global estimator reads
    \begin{equation}
        \eta_{\mathrm{res}}^2 := \sum_{T \in \T} \eta_{\mathrm{res},T}^2.
    \end{equation}

    \renewcommand{\thesection}{B}

    \section{The Zienkiewicz--Zhu estimator}\label{app:zz}
    The Zienkiewicz--Zhu estimator is a gradient recovery estimator based on an
    averaging technique introduced in\ \cite{zienkiewicz_simple_1987}.
    This estimator belongs to a general class of recovery estimators, see\
    \cite{cai_improved_2017,cai_recovery-based_2009,zhang_recovery_2001} for recent surveys
    and a reformulation of the recovery procedure in an $H(\mathrm{div})
    $-conforming space that has superior performance for problems with sharp
    interfaces.
    Despite the fact that some recovery estimators, especially when based on
    least squares fitting, are available for higher order finite elements (see for
    example \cite{zienkiewicz_superconvergent_1992}) we only consider the
    original estimator, defined for a piecewise linear finite element framework.

    Given the finite element solution $u_1 \in V^1$ the numerical flux
    $\rho_1 := \nabla u_1$ is a piecewise constant vector field.
    For each vertex $\chi \in \mathcal{N}$ in the mesh we denote
    $\omega_\chi$ the domain covered by the union of cells $T$ having common
    vertex $\chi$.
    The recovered flux $G(\rho_1) \in [V^1]^2$ has values at the degrees
    of freedom associated with the vertices $\mathcal{N}$ given by
    \begin{equation}
        G(\rho_1)(\chi) := \frac{1}{|\omega_\chi|} \int_{\omega_\chi} \rho_1\dx,
        \quad \forall \chi \in \mathcal{N}.
    \end{equation}
    The local Zienkiewicz--Zhu estimator is then defined as the discrepancy
    between the recovered flux and the numerical flux
    \begin{equation}
        \eta_{\mathrm{zz}, T} := \norm{G(\rho_1) - \rho_1}_{T},
        \quad \forall T \in \T.
    \end{equation}
    As for the residual estimator, we add Dirichlet oscillations (see
    \cref{eq:dirichlet_oscillations}) to take into
    account the Dirichlet boundary error.
    The global Zienkiewicz--Zhu estimator is given by
    \begin{equation}
        \eta_{\mathrm{zz}}^2 := \sum_{T \in \T} \eta_{\mathrm{zz}, T}^2 +
        \mathrm{osc}_{D, \overline T \cap \Gamma_D}^2.
    \end{equation}
    The code in the supplementary material contains a prototype
    implementation of the Zienkiewicz--Zhu estimator in FEniCS.
    We have implemented the local recovered flux calculation in Python rather
    than C\texttt{++}, so the runtime performance is far from optimal.

    \bibliographystyle{plain}
    \bibliography{ms}

\begin{thebibliography}{10}

\bibitem{ainsworth_influence_1996}
Mark Ainsworth.
\newblock The influence and selection of subspaces for a posteriori error
  estimators.
\newblock {\em Numerische Mathematik}, 73(4):399--418, June 1996.

\bibitem{ainsworth_posteriori_1997}
Mark Ainsworth and J.~Tinsley Oden.
\newblock A {Posteriori} {Error} {Estimators} for the {Stokes} and {Oseen}
  {Equations}.
\newblock {\em SIAM Journal on Numerical Analysis}, 34(1):228--245, February
  1997.

\bibitem{ainsworth_posteriori_2011}
Mark Ainsworth and J.~Tinsley Oden.
\newblock {\em A {Posteriori} {Error} {Estimation} in {Finite} {Element}
  {Analysis}}.
\newblock 2011.

\bibitem{cgal}
Pierre Alliez, Cl{\'e}ment Jamin, Laurent Rineau, St{\'e}phane Tayeb, Jane
  Tournois, and Mariette Yvinec.
\newblock {3D} mesh generation.
\newblock In {\em {CGAL} User and Reference Manual}. {CGAL Editorial Board},
  {5.1} edition, 2020.

\bibitem{alnaes_fenics_2015}
Martin Alnæs, Jan Blechta, Johan Hake, August Johansson, Benjamin Kehlet,
  Anders Logg, Chris Richardson, Johannes Ring, Marie~E Rognes, and Garth~N
  Wells.
\newblock The {FEniCS} {Project} {Version} 1.5.
\newblock {\em Archive of Numerical Software}, Vol 3, 2015.

\bibitem{alnaes_unified_2014}
Martin~S. Alnæs, Anders Logg, Kristian~B. Ølgaard, Marie~E. Rognes, and
  Garth~N. Wells.
\newblock Unified form language: {A} domain-specific language for weak
  formulations of partial differential equations.
\newblock {\em ACM Transactions on Mathematical Software}, 40(2):1--37,
  February 2014.

\bibitem{amestoy_fully_2001}
Patrick~R. Amestoy, Iain~S. Duff, Jean-Yves L'Excellent, and Jacko Koster.
\newblock A {Fully} {Asynchronous} {Multifrontal} {Solver} {Using}
  {Distributed} {Dynamic} {Scheduling}.
\newblock {\em SIAM Journal on Matrix Analysis and Applications}, 23(1):15--41,
  January 2001.

\bibitem{amestoy_hybrid_2006}
Patrick~R. Amestoy, Abdou Guermouche, Jean-Yves L’Excellent, and Stéphane
  Pralet.
\newblock Hybrid scheduling for the parallel solution of linear systems.
\newblock {\em Parallel Computing}, 32(2):136--156, February 2006.

\bibitem{Anciaux-Sedrakian2020}
A.~Anciaux-Sedrakian, L.~Grigori, Z.~Jorti, J.~Pape{\v{z}}, and S.~Yousef.
\newblock {Adaptive solution of linear systems of equations based on a
  posteriori error estimators}.
\newblock {\em Numer. Algorithms}, 84(1):331--364, 2020.

\bibitem{Arioli2013}
Mario Arioli, J{\"{o}}rg Liesen, Agnieszka Mi{\c{c}}dlar, and Zden{\v{e}}k
  Strako{\v{s}}.
\newblock {Interplay between discretization and algebraic computation in
  adaptive numerical solutionof elliptic PDE problems}.
\newblock {\em GAMM-Mitteilungen}, 36(1):102--129, aug 2013.

\bibitem{Aurada2013}
M.~Aurada, M.~Feischl, J.~Kemetm{\"{u}}ller, M.~Page, and D.~Praetorius.
\newblock {Each H 1/2 –stable projection yields convergence and
  quasi–optimality of adaptive FEM with inhomogeneous Dirichlet data in R d}.
\newblock {\em ESAIM Math. Model. Numer. Anal.}, 47(4):1207--1235, jul 2013.

\bibitem{bubuska_feedback_1984}
I.~Babuska and M.~Vogelius.
\newblock Feedback and adaptive finite element solution of one-dimensional
  boundary value problems.
\newblock {\em Numerische Mathematik}, 44(1):75--102, February 1984.

\bibitem{babuska_posteriori_1978}
I.~Babuška and W.~C. Rheinboldt.
\newblock A-posteriori error estimates for the finite element method.
\newblock {\em International Journal for Numerical Methods in Engineering},
  12(10):1597--1615, 1978.

\bibitem{balay_petsc_2016}
Satish Balay, Shrirang Abhyankar, Mark~F Adams, Jed Brown, Peter Brune, Kris
  Buschelman, Lisandro Dalcin, Victor Eijkhout, William~D Gropp, Dinesh
  Kaushik, Matthew~G Knepley, Lois~Curfman McInnes, Karl Rupp, Barry~F Smith,
  Stefano Zampini, Hong Zhang, and Hong Zhang.
\newblock {PETSc} {Users} {Manual}.
\newblock Technical Report ANL-95/11 - Revision 3.7, Argonne National
  Laboratory, 2016.

\bibitem{bank_posteriori_1985}
R.~E. Bank and A.~Weiser.
\newblock Some a posteriori error estimators for elliptic partial differential
  equations.
\newblock {\em Mathematics of Computation}, 44(170):283--283, May 1985.

\bibitem{bank_pltmg_1998}
Randolph~E. Bank.
\newblock {\em {PLTMG}: {A} {Software} {Package} for {Solving} {Elliptic}
  {Partial} {Differential} {Equations}: {Users}' {Guide} 8.0}.
\newblock Society for Industrial and Applied Mathematics, January 1998.

\bibitem{bank_superconvergent_2007}
Randolph~E Bank, Jinchao Xu, and Bin Zheng.
\newblock Superconvergent {Derivative} {Recovery} for {Lagrange} {Triangular}
  {Elements} of {Degree} p on {Unstructured} {Grids}.
\newblock {\em SIAM J. Numer. Anal.}, 45(5):2032--2046, 2007.

\bibitem{bartels_inhomogeneous_2004}
S.~Bartels, C.~Carstensen, and G.~Dolzmann.
\newblock Inhomogeneous {Dirichlet} conditions in a priori and a posteriori
  finite element error analysis.
\newblock {\em Numerische Mathematik}, 99(1):1--24, November 2004.

\bibitem{bartels_each_2002}
Sören Bartels and Carsten Carstensen.
\newblock Each averaging technique yields reliable a posteriori error control
  in {FEM} on unstructured grids. {Part} {II}: {Higher} order {FEM}.
\newblock {\em Mathematics of Computation}, 71(239):971--994, February 2002.

\bibitem{becker_weighted_2011}
Roland Becker, Elodie Estecahandy, and David Trujillo.
\newblock Weighted {Marking} for {Goal}-oriented {Adaptive} {Finite} {Element}
  {Methods}.
\newblock {\em SIAM Journal on Numerical Analysis}, 49(6):2451--2469, January
  2011.

\bibitem{becker_optimal_2001}
Roland Becker and Rolf Rannacher.
\newblock {An optimal control approach to a posteriori error estimation in
  finite element methods}.
\newblock {\em Acta Numer.}, 10:1--102, 2001.

\bibitem{beirao_da_veiga_priori_2008}
L.~Beirão~da Veiga, C.~Chinosi, C.~Lovadina, and R.~Stenberg.
\newblock A-priori and a-posteriori error analysis for a family of
  {Reissner}–{Mindlin} plate elements.
\newblock {\em BIT Numerical Mathematics}, 48(2):189--213, June 2008.

\bibitem{bespalov_goal-oriented_2019}
Alex Bespalov, Dirk Praetorius, Leonardo Rocchi, and Michele Ruggeri.
\newblock Goal-oriented error estimation and adaptivity for elliptic {PDEs}
  with parametric or uncertain inputs.
\newblock {\em Computer Methods in Applied Mechanics and Engineering},
  345:951--982, March 2019.

\bibitem{bespalov_t-ifiss_2019}
Alex Bespalov, Leonardo Rocchi, and David Silvester.
\newblock T-{IFISS}: a toolbox for adaptive {FEM} computation.
\newblock {\em Computers \& Mathematics with Applications}, 2020.

\bibitem{bulle_removing_2020}
Raphaël Bulle, Franz Chouly, Jack~S. Hale, and Alexei Lozinski.
\newblock Removing the saturation assumption in {Bank}–{Weiser} error
  estimator analysis in dimension three.
\newblock {\em Applied Mathematics Letters}, 107:106429, September 2020.

\bibitem{hale_implementation_2020}
Raphaël Bulle and Jack~S Hale.
\newblock An implementation of the {Bank}-{Weiser} error estimator in the
  {FEniCS} {Project} finite element software.
\newblock December 2020.

\bibitem{cai_improved_2017}
Zhiqiang Cai, Cuiyu He, and Shun Zhang.
\newblock Improved {ZZ} a posteriori error estimators for diffusion problems:
  {Conforming} linear elements.
\newblock {\em Computer Methods in Applied Mechanics and Engineering},
  313:433--449, January 2017.

\bibitem{cai_recovery-based_2009}
Zhiqiang Cai and Shun Zhang.
\newblock Recovery-{Based} {Error} {Estimator} for {Interface} {Problems}:
  {Conforming} {Linear} {Elements}.
\newblock {\em SIAM Journal on Numerical Analysis}, 47(3):2132--2156, January
  2009.

\bibitem{carstensen_axioms_2014}
C.~Carstensen, M.~Feischl, M.~Page, and D.~Praetorius.
\newblock Axioms of adaptivity.
\newblock {\em Computers \& Mathematics with Applications}, 67(6):1195--1253,
  April 2014.

\bibitem{carstensen_estimator_2010}
C.~Carstensen and C.~Merdon.
\newblock Estimator {Competition} for {Poisson} {Problems}.
\newblock {\em Journal of Computational Mathematics}, 28, 2010.

\bibitem{carstensen_each_2002}
Carsten Carstensen and Sören Bartels.
\newblock Each averaging technique yields reliable a posteriori error control
  in {FEM} on unstructured grids. {Part} {I}: {Low} order conforming,
  nonconforming, and mixed {FEM}.
\newblock {\em Mathematics of Computation}, 71(239):945--969, February 2002.

\bibitem{carstensen_robust_2016}
Carsten Carstensen and Joscha Gedicke.
\newblock Robust residual-based a posteriori {Arnold}--{Winther} mixed finite
  element analysis in elasticity.
\newblock {\em Computer Methods in Applied Mechanics and Engineering},
  300:245--264, March 2016.

\bibitem{concha2010}
{Concha INRIA Project-Team}.
\newblock {Complex Flow Simulation Codes based on High-order and Adaptive
  methods}.

\bibitem{duprez_quantifying_2020}
Michel Duprez, Stéphane Pierre~Alain Bordas, Marek Bucki, Huu~Phuoc Bui, Franz
  Chouly, Vanessa Lleras, Claudio Lobos, Alexei Lozinski, Pierre-Yves Rohan,
  and Satyendra Tomar.
\newblock Quantifying discretization errors for soft tissue simulation in
  computer assisted surgery: {A} preliminary study.
\newblock {\em Applied Mathematical Modelling}, 77:709--723, January 2020.

\bibitem{dorfler_convergent_1996}
Willy Dörfler.
\newblock A {Convergent} {Adaptive} {Algorithm} for {Poisson}’s {Equation}.
\newblock {\em SIAM Journal on Numerical Analysis}, 33(3):1106--1124, June
  1996.

\bibitem{dorfler_small_2002}
Willy Dörfler and Ricardo~H. Nochetto.
\newblock Small data oscillation implies the saturation assumption.
\newblock {\em Numerische Mathematik}, 91(1):1--12, March 2002.

\bibitem{elman_ifiss_2014}
Howard~C. Elman, Alison Ramage, and David~J. Silvester.
\newblock {IFISS}: {A} {Computational} {Laboratory} for {Investigating}
  {Incompressible} {Flow} {Problems}.
\newblock {\em SIAM Review}, 56(2):261--273, January 2014.

\bibitem{falgout_hypre_2002}
Robert~D Falgout and Ulrike~Meier Yang.
\newblock {hypre: A Library of High Performance Preconditioners}.
\newblock In Peter M~A Sloot, Alfons~G Hoekstra, C~J~Kenneth Tan, and Jack~J
  Dongarra, editors, {\em Comput. {Science} — {ICCS} 2002}, number 2331 in
  Lecture {Notes} in {Computer} {Science}, pages 632--641. Springer Berlin
  Heidelberg, apr 2002.

\bibitem{funken_efficient_2011}
Stefan Funken, Dirk Praetorius, and Philipp Wissgott.
\newblock Efficient implementation of adaptive {P1}-{FEM} in {Matlab}.
\newblock {\em Computational Methods in Applied Mathematics}, 11(4):460--490,
  2011.

\bibitem{gibson_slate_2019}
Thomas~H Gibson, Lawrence Mitchell, David~A Ham, and Colin~J Cotter.
\newblock {Slate: extending {Firedrake}'s domain-specific abstraction to
  hybridized solvers for geoscience and beyond}.
\newblock {\em Geosci. Model Dev. Discuss.}, pages 1--40, apr 2019.

\bibitem{giles_adjoint_2002}
Michael~B. Giles and Endre Süli.
\newblock Adjoint methods for {PDEs}: \textit{a posteriori} error analysis and
  postprocessing by duality.
\newblock {\em Acta Numerica}, 11:145--236, January 2002.

\bibitem{grisvard_elliptic}
P~Grisvard.
\newblock {\em {Elliptic problems in nonsmooth domains}}, volume~24 of {\em
  Monographs and Studies in Mathematics}.
\newblock Pitman (Advanced Publishing Program), Boston, MA, 1985.

\bibitem{guennebaud_eigen_2010}
Gaël Guennebaud, Benoît Jacob, and {Others}.
\newblock {\em Eigen v3}.
\newblock 2010.

\bibitem{Habera2020}
Michal Habera, Jack~S. Hale, Chris Richardson, Johannes Ring, Marie Rognes,
  Nate Sime, and Garth~N. Wells.
\newblock {FEniCSX: A sustainable future for the FEniCS Project}.
\newblock 2 2020.

\bibitem{hale_simple_2018}
Jack~S. Hale, Matteo Brunetti, Stéphane~P.A. Bordas, and Corrado Maurini.
\newblock Simple and extensible plate and shell finite element models through
  automatic code generation tools.
\newblock {\em Computers \& Structures}, 209:163--181, October 2018.

\bibitem{hale_containers_2017}
Jack~S. Hale, Lizao Li, Christopher~N. Richardson, and Garth~N. Wells.
\newblock Containers for {Portable}, {Productive}, and {Performant}
  {Scientific} {Computing}.
\newblock {\em Computing in Science \& Engineering}, 19(6):40--50, November
  2017.

\bibitem{hecht_freefem_2012}
Fr{\'{e}}d{\'{e}}ric Hecht.
\newblock {New development in freefem++}.
\newblock {\em J. Numer. Math.}, 20(3-4):1--14, jan 2012.

\bibitem{hoare_algorithm_1961}
C.~A.~R. Hoare.
\newblock Algorithm 65: find.
\newblock {\em Communications of the ACM}, 4(7):321--322, July 1961.

\bibitem{houston_sobolev}
P.~Houston, B.~Senior, and E.~S{\"{u}}li.
\newblock {Sobolev regularity estimation for hp-adaptive finite element
  methods}.
\newblock In {\em Numer. Math. Adv. Appl.}, pages 631--656. Springer Milan,
  Milano, 2003.

\bibitem{houston_automatic_2018}
Paul Houston and Nathan Sime.
\newblock Automatic {Symbolic} {Computation} for {Discontinuous} {Galerkin}
  {Finite} {Element} {Methods}.
\newblock {\em SIAM Journal on Scientific Computing}, 40(3):C327--C357, January
  2018.

\bibitem{khan_robust_2019}
Arbaz Khan, Catherine~E. Powell, and David~J. Silvester.
\newblock Robust a posteriori error estimators for mixed approximation of
  nearly incompressible elasticity.
\newblock {\em International Journal for Numerical Methods in Engineering},
  119(1):18--37, July 2019.

\bibitem{kirby_algorithm_2004}
Robert~C. Kirby.
\newblock Algorithm 839: {FIAT}, a new paradigm for computing finite element
  basis functions.
\newblock {\em ACM Transactions on Mathematical Software}, 30(4):502--516,
  December 2004.

\bibitem{kirby_compiler_2006}
Robert~C. Kirby and Anders Logg.
\newblock A compiler for variational forms.
\newblock {\em ACM Transactions on Mathematical Software}, 32(3):417--444,
  September 2006.

\bibitem{Liao_simple_2012}
Qifeng Liao and David Silvester.
\newblock A simple yet effective a posteriori estimator for classical mixed
  approximation of {Stokes} equations.
\newblock {\em Applied Numerical Mathematics}, 62(9):1242--1256, September
  2012.

\bibitem{logg_dolfin_2010}
Anders Logg and Garth~N. Wells.
\newblock {DOLFIN}: {Automated} finite element computing.
\newblock {\em ACM Transactions on Mathematical Software}, 37(2):1--28, April
  2010.

\bibitem{mitchell_collection_2017}
William~F. Mitchell.
\newblock {A collection of 2D elliptic problems for testing adaptive grid
  refinement algorithms}.
\newblock {\em Appl. Math. Comput.}, 220(February):350--364, sep 2013.

\bibitem{mitchell_collection_2013}
William~F. Mitchell.
\newblock A collection of {2D} elliptic problems for testing adaptive grid
  refinement algorithms.
\newblock {\em Applied Mathematics and Computation}, 220:350--364, September
  2013.

\bibitem{mommer_goal-oriented_2009}
Mario~S. Mommer and Rob Stevenson.
\newblock A {Goal}-{Oriented} {Adaptive} {Finite} {Element} {Method} with
  {Convergence} {Rates}.
\newblock {\em SIAM Journal on Numerical Analysis}, 47(2):861--886, January
  2009.

\bibitem{nochetto_removing_1993}
Ricardo~H. Nochetto.
\newblock Removing the saturation assumption in a posteriori error analysis.
\newblock {\em Istit. Lomb. Accad. Sci. Lett. Rend. A}, 127(1):67--82 (1994),
  1993.

\bibitem{nochetto_theory_2009}
Ricardo~H Nochetto, Kunibert~G Siebert, and Andreas Veeser.
\newblock {Theory of adaptive finite element methods: An introduction}.
\newblock In Ronald DeVore and Angela Kunoth, editors, {\em Multiscale,
  Nonlinear Adapt. Approx.}, pages 409--542, Berlin, Heidelberg, 2009. Springer
  Berlin Heidelberg.

\bibitem{Papez2018}
J.~Pape{\v{z}}, Z.~Strako{\v{s}}, and M.~Vohral{\'{i}}k.
\newblock {Estimating and localizing the algebraic and total numerical errors
  using flux reconstructions}.
\newblock {\em Numer. Math.}, 138(3):681--721, mar 2018.

\bibitem{pfeiler_dorfler_2019}
Carl-Martin Pfeiler and Dirk Praetorius.
\newblock Dörfler marking with minimal cardinality is a linear complexity
  problem.
\newblock {\em arXiv1907.13078 [cs, math]}, July 2019.

\bibitem{plaza_local_2000}
A.~Plaza and G.F. Carey.
\newblock Local refinement of simplicial grids based on the skeleton.
\newblock {\em Applied Numerical Mathematics}, 32(2):195--218, February 2000.

\bibitem{prudhomme_feel_2012}
Christophe Prud’homme, Vincent Chabannes, Vincent Doyeux, Mourad Ismail,
  Abdoulaye Samake, and Goncalo Pena.
\newblock Feel++ : {A} computational framework for {Galerkin} {Methods} and
  {Advanced} {Numerical} {Methods}.
\newblock {\em ESAIM: Proceedings}, 38:429--455, December 2012.

\bibitem{rathgeber_firedrake_2017}
Florian Rathgeber, David~A. Ham, Lawrence Mitchell, Michael Lange, Fabio
  Luporini, Andrew T.~T. Mcrae, Gheorghe-Teodor Bercea, Graham~R. Markall, and
  Paul H.~J. Kelly.
\newblock Firedrake: {Automating} the {Finite} {Element} {Method} by
  {Composing} {Abstractions}.
\newblock {\em ACM Transactions on Mathematical Software}, 43(3):1--27, January
  2017.

\bibitem{RenardGetFEM}
Yves Renard and K~Poulios.
\newblock {GetFEM: Automated FE modeling of multiphysics problems based on a
  generic weak form language}.
\newblock {\em (Submitted)}, 2020.

\bibitem{rodriguez_remarks_1994}
Rodolfo Rodríguez.
\newblock Some remarks on {Zienkiewicz}-{Zhu} estimator.
\newblock {\em Numerical Methods for Partial Differential Equations},
  10(5):625--635, September 1994.

\bibitem{rognes_automated_2013}
Marie~E. Rognes and Anders Logg.
\newblock Automated {Goal}-{Oriented} {Error} {Control} {I}: {Stationary}
  {Variational} {Problems}.
\newblock {\em SIAM Journal on Scientific Computing}, 35(3):C173--C193, January
  2013.

\bibitem{rupin_experimental_2008}
Fabienne Rupin, Amena Saied, Davy Dalmas, Françoise Peyrin, Sylvain Haupert,
  Etienne Barthel, Georges Boivin, and Pascal Laugier.
\newblock Experimental determination of {Young} modulus and {Poisson} ratio in
  cortical bone tissue using high resolution scanning acoustic microscopy and
  nanoindentation.
\newblock {\em J. Acoust. Soc. Am.}, 123(5):3785--3785, May 2008.

\bibitem{scott_finite_1990}
L.~Ridgway Scott and Shangyou Zhang.
\newblock Finite element interpolation of nonsmooth functions satisfying
  boundary conditions.
\newblock {\em Mathematics of Computation}, 54(190):483--483, May 1990.

\bibitem{van_der_walt_numpy_2011}
Stéfan van~der Walt, S~Chris Colbert, and Gaël Varoquaux.
\newblock The {NumPy} {Array}: {A} {Structure} for {Efficient} {Numerical}
  {Computation}.
\newblock {\em Computing in Science \& Engineering}, 13(2):22--30, March 2011.

\bibitem{varrette_management_2014}
Sebastien Varrette, Pascal Bouvry, Hyacinthe Cartiaux, and Fotis Georgatos.
\newblock Management of an academic {HPC} cluster: {The} {UL} experience.
\newblock In {\em 2014 {International} {Conference} on {High} {Performance}
  {Computing} \& {Simulation} ({HPCS})}, pages 959--967, Bologna, Italy, July
  2014. IEEE.

\bibitem{verfurth_posteriori_1994}
R.~Verfürth.
\newblock A posteriori error estimation and adaptive mesh-refinement
  techniques.
\newblock {\em Journal of Computational and Applied Mathematics},
  50(1-3):67--83, May 1994.

\bibitem{virtanen_scipy_2020}
Pauli Virtanen, Ralf Gommers, Travis~E. Oliphant, Matt Haberland, Tyler Reddy,
  David Cournapeau, Evgeni Burovski, Pearu Peterson, Warren Weckesser, Jonathan
  Bright, Stéfan~J. van~der Walt, Matthew Brett, Joshua Wilson, K.~Jarrod
  Millman, Nikolay Mayorov, Andrew R.~J. Nelson, Eric Jones, Robert Kern, Eric
  Larson, C.~J. Carey, Ilhan Polat, Yu~Feng, Eric~W. Moore, Jake VanderPlas,
  Denis Laxalde, Josef Perktold, Robert Cimrman, Ian Henriksen, E.~A. Quintero,
  Charles~R. Harris, Anne~M. Archibald, Antônio~H. Ribeiro, Fabian Pedregosa,
  Paul van Mulbregt, and Contributors.
\newblock {SciPy} 1.0--{Fundamental} {Algorithms} for {Scientific} {Computing}
  in {Python}.
\newblock {\em Nature Methods}, 17(3):261--272, March 2020.

\bibitem{zhang_recovery_2001}
Zhimin Zhang and Ningning Yan.
\newblock Recovery type a posteriori error estimates in finite element methods.
\newblock {\em J. Appl. Math. Comput.}, 8(2):235--251, 2001.

\bibitem{zienkiewicz_simple_1987}
O.~C. Zienkiewicz and J.~Z. Zhu.
\newblock A simple error estimator and adaptive procedure for practical
  engineering analysis.
\newblock {\em International Journal for Numerical Methods in Engineering},
  24(2):337--357, February 1987.

\bibitem{zienkiewicz_superconvergent_1992}
O.~C. Zienkiewicz and J.~Z. Zhu.
\newblock The superconvergent patch recovery ({SPR}) and adaptive finite
  element refinement.
\newblock {\em Computer Methods in Applied Mechanics and Engineering},
  101(1):207--224, December 1992.

\bibitem{olgaard_automated_2009}
Kristian~B. Ølgaard, Anders Logg, and Garth~N. Wells.
\newblock Automated {Code} {Generation} for {Discontinuous} {Galerkin}
  {Methods}.
\newblock {\em SIAM Journal on Scientific Computing}, 31(2):849--864, January
  2009.

\end{thebibliography}
\end{document}


\maketitle
\section{Indicative snippet of error estimation for Poisson
equation using Bank--Weiser estimator}\label{app:poisson}
We present here a snippet of DOLFIN Python code showing function to compute
the error of a Poisson problem using the Bank--Weiser estimator.

\begin{minted}[fontsize=\small, frame=lines, gobble=1]{python}
    from dolfin import *
    import fenics_error_estimation

    def estimate(u_h):
        """Bank-Weiser error estimation procedure for the Poisson problem.

        Parameters
        -----------
        u_h: dolfin.Function
        Solution of Poisson problem.

        Returns
        -------
        The error estimate on each cell of the mesh.
        """
        mesh = u_h.function_space().mesh()

        # Higher order space
        element_f = FiniteElement("DG", triangle, 2)
        # Low order space
        element_g = FiniteElement("DG", triangle, 1)

        # Construct the Bank-Weiser interpolation operator according to the
        # definition of the high and low order spaces.
        N = fenics_error_estimation.create_interpolation(element_f, element_g)

        V_f = FunctionSpace(mesh, element_f)
        e = TrialFunction(V_f)
        v = TestFunction(V_f)
        f = Constant(0.0)

        # Homogeneous zero Dirichlet boundary conditions
        bcs = DirichletBC(V_f, Constant(0.0), "on_boundary", "geometric")

        # Define the local Bank-Weiser problem on the full higher order space
        n = FacetNormal(mesh)
        a_e = inner(grad(e), grad(v))*dx
        # Residual
        L_e = inner(f + div(grad(u_h)), v)*dx + \
                inner(jump(grad(u_h), -n), avg(v))*dS

        # Local solves on the implied Bank-Weiser space. The solution is returned
        # on the full space.
        e_h = fenics_error_estimation.estimate(a_e, L_e, N, bcs)

        # Estimate of global error
        error = norm(e_h, "H10")

        # Computation of local error indicator.
        V_e = FunctionSpace(mesh, "DG", 0)
        v = TestFunction(V_e)

        eta_h = Function(V_e, name="eta_h")
        # By testing against v in DG_0 this effectively computes
        # the estimator on each cell.
        eta = assemble(inner(inner(grad(e_h), grad(e_h)), v)*dx)
        eta_h.vector()[:] = eta

        return eta_h
\end{minted}

\newpage
\section{Indicative snippet of error estimation for linear elasticity
equations using Poisson estimator}\label{app:stokes}
We give here a snippet of DOLFIN Python code showing function to compute the
error of a two-dimensional linear elasticity problem (discretized with
Taylor--Hood element) using the Poisson estimator, based on our implementation
of the Bank--Weiser estimator.

\begin{minted}[fontsize=\small, frame=lines, gobble=1]{python}
    import scipy.linalg as sp.linalg
    from dolfin import *
    import fenics_error_estimation

    def estimate(w_h, mu, lmbda):
        """
        Parameters
        -----------
        w_h: dolfin.Function
        Solution of the linear elasticity problem.
        mu: float
        First Lamé coefficient.
        lmbda: float
        Second Lamé coefficient.

        Returns
        -------
        The error estimate on each cell of the mesh.
        """

        mesh = w_h.function_space().mesh()

        u_h = w_h.sub(0)
        p_h = w_h.sub(1)

        # Vectorial high order space.
        X_element_f = VectorElement('DG', triangle, 3)

        # Scalar high order and low order spaces.
        S_element_f = FiniteElement('DG', triangle, 3)
        S_element_g = FiniteElement('DG', triangle, 2)

        # Construct the scalar projection matrix according to the definition
        # of the high and low order spaces.
        N_S = create_interpolation(S_element_f, S_element_g)

        # Construct the vectorial projection matrix as a block diagonal, each
        # block corresponding to a scalar problem.
        N_X = sp.linalg.block_diag(N_S, N_S)

        f = Constant((0., 0.))

        X_f = FunctionSpace(mesh, X_element_f)
        e_X = TrialFunction(X_f)
        v_X = TestFunction(X_f)

        # Homogeneous zero Dirichlet boundary conditions.
        bcs = DirichletBC(X_f, Constant((0., 0.)), 'on_boundary', 'geometric')

        # Cell residual.
        R_T = f + div(2.*mu*sym(grad(u_h))) - grad(p_h)

        # Facet residual.
        n = FacetNormal(mesh)
        R_E = (1./2.)*jump(p_h*Identity(2) - 2.*mu*sym(grad(u_h)), -n)

        # Local Poisson problem.
        a_X_e = 2.*mu*inner(grad(e_X), grad(v_X))*dx
        L_X_e = inner(R_K, v_X)*dx - inner(R_E, avg(v_X))*dS

        # Solve Poisson equation locally on implicit Bank--Weiser space.
        e_h = fenics_error_estimation.estimate(a_X_e, L_X_e, N_X, bcs)

        # Cell residual.
        rho_d = 1./(lmbda**(-1)+(2.*mu)**(-1))
        r_T = rho_d*(div(u_h) + lmbda**(-1)*p_h)

        # Computation of local error indicator.
        V_e = FunctionSpace(mesh, 'DG', 0)
        v = TestFunction(V_e)

        eta_h = Function(V_e)
        # By testing against v in DG_0 this effectively computes the estimator
        # on each cell.
        eta = assemble(2.*mu*inner(inner(grad(e_h), grad(e_h)), v)*dx + \
              rho_d**(-1)*inner(inner(eps_h, eps_h), v)*dx)
        eta_h.vector()[:] = eta

        return eta_h
\end{minted}